\documentclass[a4paper,11pt]{article}

\usepackage{url}
\usepackage{color}
\usepackage{common}
\usepackage{enumerate}
\usepackage{graphicx}
\usepackage[mathscr]{euscript}		%ermöglicht Umlaute
\usepackage{amsthm}
\usepackage{dsfont}
\usepackage[active]{srcltx}

\usepackage{psfrag}			%um Text in eps-Dokumenten mit latex zu kompilieren
\usepackage{hyperref}

\newcommand{\ind}{\mathds}
\newcommand{\floor}[1]{{\lfloor #1 \rfloor}}

\theoremstyle{plain}
\newtheorem{thm}{Theorem}[section]
\newtheorem{corollary}[thm]{Corollary}
\newtheorem{prop}[thm]{Proposition}
\newtheorem{lem}[thm]{Lemma}

\theoremstyle{definition}
\newtheorem{definition}[thm]{Definition}

\theoremstyle{remark}
\newtheorem{remark}[thm]{Remark}

\DeclareMathOperator{\esssup}{ess\,sup}

\DeclareMathOperator{\Prob}{Prob}

\numberwithin{equation}{section}

\title{LYAPUNOV EXPONENTS FOR THE ONE-DIMENSIONAL PARABOLIC ANDERSON MODEL WITH DRIFT}
\author{Alexander Drewitz\footnote{
{\it E-mail address: \href{mailto:drewitz@math.tu-berlin.de}{drewitz@math.tu-berlin.de}, Tel.: +49 30 314 23606.}} \vspace{1em}\\
{\it \footnotesize
Institut f\"ur Mathematik, Technische Universit\"at Berlin, MA 7-5, Str. des 17. Juni 136, 10623 Berlin, Germany
}
}
\date{\today}

\begin{document}

% FRAGEN:
% Explizite Abh\"angigkeit der $\lambda_p$ von $h$ (mon. non-increasing?) zeigen?
% 
% Annealed Lyap-Exp unstetig, oder doch $\beta_cr = \beta_cr^0$?
% 
% Proposition gross oder klein schreiben im Beweis derselben?

\maketitle

\begin{abstract}
We consider the solution $u$ to the one-dimensional parabolic Anderson model with homogeneous
initial condition $u(0, \cdot ) \equiv 1$, arbitrary drift and a time-independent potential bounded from above. Under
ergodicity and independence conditions
we derive representations for both the quenched Lyapunov exponent and, more importantly, the $p$-th annealed Lyapunov exponents
for {\it all} $p \in (0, \infty).$

These results enable us to
prove the heuristically plausible fact that the $p$-th annealed Lyapunov exponent converges
to the quenched Lyapunov exponent as $p \downarrow 0.$
Furthermore, we
show that $u$ is $p$-intermittent for $p$ large enough.

As a byproduct, we
compute the optimal quenched speed of the random walk appearing in the Feynman-Kac representation of $u$ under the corresponding
Gibbs measure. In this context, depending on the negativity of the potential, a phase transition from zero speed to positive speed
appears.

\vspace{.5em}
{\footnotesize

}
\end{abstract}

\section{Model and notation} \label{notation}

We consider the one-dimensional parabolic Anderson model with arbitrary drift and homogeneous initial
condition, i.e. the Cauchy problem
\begin{align}
\left. \begin{array}{rll}
\dfrac{\partial u}{\partial t} (t,x)
\hspace{-,5em} &= \kappa \Delta_h u(t,x) + \xi(x) u(t,x), & (t,x) \in \R_+ \times \Z, \vspace{.5em}\\
u(0,x) \hspace{-,5em}&= 1, & x \in \Z,
\end{array} \right\} \label{PAMArb}
\end{align}
where
$\kappa$ is a positive diffusion constant,
$h \in (0, 1]$ an arbitrary drift and
$\Delta_h$ denotes the discrete Laplace operator with drift $h$ given by
$$
\Delta_{h} u(t,x):= \frac{1+h}{2}(u(t,x+1) - u(t,x)) + \frac{1-h}{2}(u(t,x-1) - u(t,x)).
$$
Here and in the following,
$( \xi(x))_{x \in \Z} \in \Sigma := (-\infty,0]^{\Z}$ is a non-constant ergodic potential bounded from above.
The distribution of $\xi$ will be denoted by $\Prob$ and the corresponding expectation by $\langle \cdot \rangle.$
In our context, ergodicity is understood with respect to the left-shift $\theta$ acting on $\zeta \in \Sigma$ via
$
\theta((\zeta(x))_{x \in \Z}) := (\zeta(x+1))_{x \in \Z}.
$
Without further loss of generality, we will
assume
\begin{equation} \label{potEssSup}
\esssup \xi(0) = 0.
\end{equation}
The case $\esssup \xi(0) = c$ reduces to (\ref{potEssSup}) by the transformation $u \mapsto e^{ct}u.$
% Denoting $\Sigma:= (-\infty, 0]^\Z$ and assuming $\xi$ to be the identity on 
% a probability space $(\Sigma, {\cal B}(\Sigma), \Prob),$
We therefore have
\begin{equation} \label{potErg}
\Prob \in {\cal M}_1^e(\Sigma),
\end{equation}
where ${\cal M}_1(E)$ denotes the space of probability measures on a topological space $E.$
If we have a shift operator defined on $E$ (such as $\theta$ for $E=\Sigma$), then by ${\cal M}_1^s(E)$
and ${\cal M}_1^e(E)$ we denote the spaces of shift-invariant and ergodic probability measures on $E,$
respectively.
If not mentioned otherwise, we will always
assume the measures to be defined on the corresponding Borel $\sigma$-algebra
and the spaces of measures to be endowed with the topology of weak convergence.
We denote
$\Sigma_b := [b,0]^\Z$ and $\Sigma_b^+ := [b,0]^{\N_0}$ for $b \in (-\infty, 0)$.
Since the potential plays the role of a (random) medium,
we likewise refer to $\xi$ as the medium.
% In particular, for 
% $
% g: (\Sigma, {\cal B}(\Sigma)) \to (\R, {\cal B}(\R))
% $
% and
% $\nu \in {\cal M}_1(\Sigma)$, i.e. $\nu$ is a probabilty measures on $(\Sigma, {\cal B}(\Sigma))$,
% we define
% $
% \langle g \circ \xi \rangle_\nu := \int_X g \circ \xi \, d\Prob,
% $
% where the index $\nu$ means that we assume the potential $\xi$ to be distributed according to $\nu$, i. e. $\Prob \circ \xi^{-1} = \nu$.

Examples motivating the study of (\ref{PAMArb}) reach
from chemical kinetics (cf. \cite{GaMo-90} and \cite{CaMo-94}) to evolution theory (see \cite{EbEnEsFe-84}).
In particular, we may associate
to (\ref{PAMArb}) the following branching particle system: At time $t=0$ at each site $x \in \Z$ there starts
a particle moving independently from all others according to a continuous-time random walk with generator
$\kappa \Delta_{-h}.$ It is killed with rate $\xi^-$ and splits into two with rate $\xi^+.$
Each descendant moves independently from all other particles according to the same law as its ancestor. The expected
number of particles at time $t$ and site $x$ given the medium $\xi$ solves equation (\ref{PAMArb}).

\section{Main results} \label{mainResults}

Our central interest is in the quenched and $p$-th annealed Lyapunov exponents, which if they exist, are given by
\begin{equation} \label{quenchedLyapExpDef}
\lambda_0 := \lim_{t \to \infty} \frac1t \log u(t,0) \quad \text{a.s.}
\end{equation}
and
\begin{equation} \label{annealedLyapExpDef}
\lambda_p := \lim_{t \to \infty} \frac1t \log \langle u(t,0)^p \rangle^{1/p}, \quad p \in (0,\infty),
\end{equation}
respectively (cf. Theorems \ref{quenchedLyapExp} and \ref{annealedLyapExp}).

Of particular interest to us is the occurrence of intermittency of the solution $u,$
which heuristically
means that $u$ is irregular and exhibits pronounced spatial peaks.
The motivation for this is as follows. In equation (\ref{PAMArb}) two competing effects are present. On
one hand, the operator $\kappa \Delta_h$ induces a diffusion (in combination with a drift)
which tends to smooth the solution. On the
other hand, the influence of the random potential $\xi$ favours the spatial inhomogeneity of the solution $u.$
The occurrence of intermittency is therefore evidence that the influence of the potential
dominates the effect of diffusion.

A standard tool in the study of intermittency is in terms of the exponential growth of moments.\footnote{
An explicit geometric characterisation is more difficult and beyond the scope of this
article. For the zero drift case, see
 G\"artner, K\"onig and Molchanov \cite{GaKoMo-07}.}

\begin{definition}
 For $p \in (0, \infty),$ the solution $u$ to (\ref{PAMArb}) is called {\it p-intermittent} if
$
\lambda_{p+\varepsilon} > \lambda_p
$
for all $\varepsilon > 0$ sufficiently small.
\end{definition}

\begin{remark}
 It will turn out that $p$-intermittency implies $\lambda_{p+\varepsilon} > \lambda_p$
for all 
$\varepsilon > 0,$ cf. Proposition \ref{lyapExpGenProp} (a).
\end{remark}

Due to Jensen's inequality, $\lambda_{p + \varepsilon} \geq \lambda_p$ is always fulfilled;
$p$-intermittency manifests itself in the {\it strict} inequality.
In this case,
Chebyshev's inequality yields
$$
\Prob( u(t,0) > e^{\alpha t}) \leq e^{-\alpha p t} \langle u(t,0)^p \rangle \asymp e^{ (- \alpha +\lambda_p) p t }\to 0
$$
for $\alpha \in (\lambda_p, \lambda_{p+\varepsilon})$
and, at the same time,
$$
\langle u(t,0)^{p+\varepsilon} \ind{1}_{u(t,0) \leq e^{\alpha t}} \rangle
\leq e^{\alpha (p+\varepsilon) t} = o(\langle u(t,0)^{p+\varepsilon} \rangle)
$$
as $t \to \infty,$ which again implies
$$
\langle u(t,0)^{p+\varepsilon} \ind{1}_{u(t,0) > e^{\alpha t}} \rangle \sim \langle u(t,0)^{p+\varepsilon} \rangle.
$$
In particular, setting $\Gamma(t) := \{x \in \Z : u(t,x) > e^{\alpha t} \}$
and considering large centered intervals $I_t$ yields
$$
\vert I_t \vert^{-1} \sum_{x \in I_t} u(t,x)^{p+\varepsilon} \approx \vert I_t \vert^{-1} \sum_{x \in I_t \cap \Gamma(t)} u(t,x)^{p+ \varepsilon},
$$
due to Birkhoff's ergodic theorem. This justifies
the interpretation that for large times the solution $u$
develops (relatively) higher and higher peaks on fewer and fewer islands.
For further reading, see \cite{GaKo-05} and \cite{GaMo-90}.

In systems such as the one we are considering, one heuristically derives
\begin{align*}
 \lim_{p \downarrow 0} \frac{1}{p} \log \langle u(t,0)^p \rangle
= \big [\frac{d}{dp} \log \langle u(t,0)^p \rangle \big ] \big \vert_{p=0}
= \frac{\langle (\log u(t,0)) u(t,0)^p \rangle}{\langle u(t,0)^p \rangle} \Big\vert_{p=0} = \langle \log u(t,0) \rangle.
\end{align*}
It is therefore generally believed that the
$p$-th annealed
Lyapunov exponent converges to the quenched Lyapunov exponent
as $p \downarrow 0.$ Due to the fact that we are
able to compute the $p$-th annealed Lyapunov exponent for all $p \in (0, \infty),$ we can
prove this conjecture, cf. Theorem \ref{lyapExpCont}.

In order to formulate our results, we introduce some more notation.
Let
$Y= (Y_t)_{t \in \R_+}$ be a continuous-time random walk on $\Z$ with generator
$\kappa \Delta_{-h}.$
By $\P_x$ we denote the underlying probability measure
with $\P_x(Y_0 = x) = 1$ and we write
$\E_x$ for the expectation with respect to $\P_x.$
Let
$T_n$ be the first hitting time of $n \in \Z$ by $Y$ and define for $\beta \in \R,$
$$
L^+(\beta) := \Big \langle \Big( \log \E_1 \exp \Big\{ \int_0^{T_0} (\xi(Y_s) + \beta) \, ds \Big\} \Big)^+ \Big \rangle,
$$
$$
L^-(\beta) := \Big \langle \Big( \log \E_1 \exp \Big\{ \int_0^{T_0} (\xi(Y_s) + \beta) \, ds \Big\} \Big)^- \Big \rangle
$$
as well as
\begin{equation} \label{LDef}
L(\beta) := L^+ (\beta) - L^-(\beta)
\end{equation}
if this expression is well-defined, i.e. if at least one of the two terms on the right-hand side is finite.
We denote $\beta_{cr}$ the critical value such that
$L^+(\beta) = \infty$ for $\beta > \beta_{cr}$ and $L^+(\beta) < \infty$ for
$\beta < \beta_{cr}$. With this notation, we observe that $L(\beta)$ is well-defined for all $\beta \in (-\infty, \beta_{cr})$
at least, and taking into account (\ref{potEssSup}) one can show without much effort that
\begin{equation} \label{betaCrEst}
\beta_{cr} \in [0, \kappa].
\end{equation}

\subsection{Quenched regime} \label{subsecQuenchedRegime}

We start by considering the quenched Lyapunov exponent 
$
\lambda_0
$
and note that even the existence of the limit on the right-hand side of (\ref{quenchedLyapExpDef}) is not immediately obvious.
We will impose that either the random field $\xi$ is ergodic and bounded, i.e.
\begin{equation} \label{potErgBdd}
\Prob \in {\cal M}_1^e(\Sigma_b)
\end{equation}
for some $b < 0,$ or that $\xi$ consists of i.i.d. random variables, i.e.
\begin{equation} \label{potIID}
\Prob = \eta^\Z
\end{equation}
for some law $\eta \in {\cal M}_1 ((-\infty, 0])$.
Note that if (\ref{potIID}) is fulfilled, standard computations in combination with Lemma \ref{betaCrBeta}
yield $\beta_{cr} = \kappa (1-\sqrt{1-h^2}).$
In the context described above, we then have the following result.

\begin{thm} \label{quenchedLyapExp}
Assume (\ref{potEssSup}) and either (\ref{potErgBdd}) or (\ref{potIID}). Then
the quenched Lyapunov exponent $\lambda_0$
exists a.s. and is non-random. Furthermore, $\lambda_0$ equals the zero of
$
\beta \mapsto L(-\beta)
$
in $(-\beta_{cr}, 0)$ or, if such a zero does not exist, equals $-\beta_{cr}$.
% \begin{align}
% \lambda_0 
% &= 
% \left\{ \begin{array}{r@{\quad:\quad}l}
% %c^* & \text{ if } L^- (0) = \infty\\
% -\beta_{cr} & \text{ if } L^- (0) = \infty 
% \text{ or }\lim_{\beta \uparrow \beta_{cr}} L(\beta) \leq 0\\
% \text{zero of } \beta \mapsto L(-\beta) \text{ in } (-\beta_{cr}, 0) 
% & \text{ otherwise.}
% \end{array} \right\} \label{quenchedLyapExpZeroFormula}
% \end{align}
\end{thm}

For an outline of the proof of Theorem \ref{quenchedLyapExp} and in order to understand
the corollary below, we remark that
the unique bounded non-negative solution to (\ref{PAMArb}) is given
by the Feynman-Kac formula
\begin{align}\label{feynmanKac}
\begin{split}
u(t,x) = \E_0 \exp \Big\{ \int_0^t \xi(X_s) \, ds \Big\} &= \sum_{n \in \Z} \E_0 \exp \Big\{ \int_0^t \xi(X_s) \, ds \Big \} \ind{1}_{X_t = n}\\
&= \sum_{n \in \Z} \E_n \exp \Big \{ \int_0^t \xi(Y_s) \, ds \Big \} \ind{1}_{Y_t = 0}.
\end{split}
\end{align}
Here, in analogy to $Y$ we denote by $X = ( X_t )_{t \in \R_+}$ a continuous-time random walk
on $\Z$ with generator $\kappa \Delta_h.$
Note that $X$ and $Y$ may be regarded as time reversals of each other.

The proof of Theorem \ref{quenchedLyapExp} 
roughly proceeds as follows. Considering the Feynman-Kac representation
(\ref{feynmanKac}) of $u,$ the main contributions 
to (\ref{feynmanKac})
stem from summands with $n \approx \alpha t,$ i.e.
\begin{equation}\label{singleSummand}
u(t,0) \asymp \E_0 \exp \Big\{ \int_0^t \xi(X_s) \, ds \Big\} \ind{1}_{X_t \approx \alpha t}
\end{equation}
Here, $\alpha \geq 0$ denotes the optimal speed of the random walk $X$ within the random medium, cf.
Corollary \ref{optimalSpeed} below.
To show the desired behaviour we use large deviations 
for certain hitting times of the random walk, which then yield
a variational formula for $\lambda_0,$ cf. Corollary \ref{quenchedLyapExpVar}.
With this formula at hand, it is an easy task to complete the proof of Theorem \ref{quenchedLyapExp}.
See section \ref{proofsQuenched} for further details.

As a byproduct we obtain the following corollary on the optimal speed of the random walk
$X$ under the Gibbs measure
$$
\P_t^\xi (\cdot) := \frac{\E_0 \exp \{ \int_0^t \xi(X_s) \, ds \} \ind{1}_{X_t \in \cdot}}
{\E_0 \exp \{ \int_0^t \xi(X_s) \, ds \}}
$$
on $\R.$

\begin{corollary} \label{optimalSpeed}
Let the assumptions of Theorem \ref{quenchedLyapExp} be fulfilled.
\begin{enumerate}
 \item 
If $\lim_{\beta \uparrow \beta_{cr}} L(\beta) > 0,$ then for all $\varepsilon > 0,$
$$
\lambda_0 > \limsup_{t \to \infty} \frac1t \log \E_0 \exp \Big\{ \int_0^t \xi(X_s) \, ds \Big\}
\ind{1}_{X_t \notin t(\alpha^*-\varepsilon, \alpha^* +\varepsilon )}
$$
with
$$
\alpha^* := (L'(-\lambda_0))^{-1} = \Big \langle \frac{\E_1 T_0 \exp\{\int_0^{T_0} (\xi(Y_s) - \lambda_0) \, ds\}}
{\E_1 \exp\{\int_0^{T_0} (\xi(Y_s) -\lambda_0) \, ds \}} \Big \rangle^{-1} \in (0,\infty).
$$

\item

If $\lim_{\beta \uparrow \beta_{cr}} L(\beta) = 0,$ then for $m \in [0, (\lim_{\beta \uparrow \beta_{cr}} L'(\beta))^{-1}]$
and all $\varepsilon > 0,$
$$
\lambda_0 = \liminf_{t \to \infty} \frac1t \log \E_0 \exp \Big\{ \int_0^t \xi(X_s) \, ds \Big\}
\ind{1}_{X_t \in t(m-\varepsilon, m + \varepsilon)},
$$
while 
$$
\lambda_0 > \limsup_{t \to \infty} \frac1t \log \E_0 \exp \Big\{ \int_0^t \xi(X_s) \, ds \Big\}
\ind{1}_{X_t \notin t[-\varepsilon, (\lim_{\beta \uparrow \beta_{cr}} L'(\beta))^{-1} + \varepsilon]}.
$$

\item

If $\lim_{\beta \uparrow \beta_{cr}} L(\beta) < 0,$ then for all $\varepsilon > 0,$
$$
\lambda_0 > \limsup_{t \to \infty} \frac1t \log \E_0 \exp \Big\{ \int_0^t \xi(X_s) \, ds \Big\}
\ind{1}_{X_t \notin t(-\varepsilon, \varepsilon)}.
$$
\end{enumerate}
\end{corollary}

\begin{remark}
\begin{enumerate}
\item
The existence of $L'$ under the assumptions of part (b) from above
will be shown in Lemma \ref{lDiff} below.
Since $L$ is increasing and
convex on $(-\infty, \beta_{cr}),$
the limit $\lim_{\beta \uparrow \beta_{cr}} L'(\beta) > 0$ then
exists.

% \item
% In other words, the above result shows that
% typical paths of $X$ under the measures defined by
% $$
% A \mapsto \frac{\E_0 \exp \{ \int_0^t \xi(X_s) \, ds \} \ind{1}_{A}}{\E_0 \exp \{ \int_0^t \xi(X_s) \, ds \}}
% $$
% for $A \in \sigma \{ X_s : s \leq t\}$ have non-random speed given by $\alpha^*,$
% $m \in [0, \lim_{\beta \uparrow \beta_{cr}} L'(\beta) > 0]$ and $0,$ respectively. 

\item
Part (b) of the corollary can be viewed as a phase transition between the cases (a) and (c),
which correspond to the non-localised and localised regimes, respectively. Note that
the result of (c) can be considered a screening effect, where
the random walk is prevented from moving with positive speed due to the distribution of $\xi$ putting
much mass on very negative values,
cf. also \cite{BiKo-01}.

\item
Inspecting the proof of this corollary, one may observe that continuing the corresponding ideas
we would obtain a large deviations principle for the position of the random walk under the above Gibbs
measure. However, since our emphasis is rather on Lyapunov exponents and intermittency, we will
not carry out the necessary modifications.

\end{enumerate}
\end{remark}
% 
% \begin{remark} \label{geomInterpretLoc}
% Corollary \ref{optimalSpeed} is interesting also because it directly offers an (albeit coarse) insight into the geometrical structure
% of the solution to the following modification of (\ref{PAMArb}):
% \begin{align}
% \left. \begin{array}{rll}
% \dfrac{\partial v}{\partial t} (t,x)
% \hspace{-,5em} &= \kappa \Delta_{-h} v(t,x) + \xi(x) v(t,x), & (t,x) \in \R_+ \times \Z, \vspace{.5em}\\
% v(0,x) \hspace{-,5em}&= \delta_0(x), & x \in \Z,
% \end{array} \right\} \label{PAMArbLoc}
% \end{align}
% where drift and initial condition have been changed.
% Indeed, the Feynman-Kac formula and time reversal yield
% $$
% v(t,x) = \E_x \exp \Big\{ \int_0^t \xi(Y_s) \, ds \Big\} \delta_0(Y_t)
% =\E_0 \exp \Big\{ \int_0^t \xi(X_s) \, ds \Big\} \ind{1}_{X_t= x}.
% $$
% Corollary \ref{optimalSpeed} with the corresponding notations
% then yields that the essential peaks of $v$ are found at small islands around either
% $\alpha^*t,$ $mt$ with $m \in [0, (\lim_{\beta \uparrow \beta_{cr}} L'(\beta))^{-1}],$ or $0,$ respectively,
% depending on the behaviour of $L$ near $\beta_{cr}.$
% \end{remark}

\subsection{Annealed regime} \label{annealedRegimeSubSec}
In order to avoid technical difficulties,
we always assume in the annealed case that
\begin{equation} \label{potIIDBdd}
\Prob = \eta^\Z
\end{equation}
for some $\eta \in {\cal M}_1 ([b, 0])$ and $b \in (-\infty, 0).$
We are interested in
the existence of the annealed Lyapunov exponents
$
\lambda_p
$
for all $p > 0$
and will derive specific formulae for them.
The proof will use process level large deviations applied to the random medium $\xi.$
In order to be able to formulate our result, we have to introduce some further notation.
For $\zeta \in \Sigma_b$ we denote by $R_n(\zeta)$ the restriction of the empirical measure
$
n^{-1} \sum_{k=0}^{n-1} \delta_{\theta^k \circ \zeta} \in {\cal M}_1(\Sigma_b)
$
to ${\cal M}_1(\Sigma_b^+).$\footnote{If clear from the context,
we will interpret elements $ \nu$ of ${\cal M}_1 (\Sigma_b)$ as elements of ${\cal M}_1 (\Sigma_b^+)$
without further mentioning by considering $\nu \circ \pi_+^{-1}$ instead, where
$
\pi_+ : (x_n)_{n \in \Z} \mapsto (x_n)_{n \in \N_0}.
$
In the same fashion, we consider elements of $\Sigma_b$ as elements of $\Sigma_b^+.$
}
Using assumption (\ref{potIIDBdd}) we get that the uniformity condition (U) in section 6.3 of \cite{DeZe-98} is satisfied
for $(\xi(x))_{x \in \Z}.$
Hence, Corollaries 6.5.15 and 6.5.17 of the same reference
provide us with a full process level large deviations principle for the random sequence
of empirical measures $(R_n \circ \xi)_{n \in \N}$
on scale $n$ with rate function given by
\begin{align} \label{processLevelRateFunction}
\begin{split}
{\cal I}(\nu) := 
\left\{ \begin{array}{rl}
H(\nu_1^* \vert \nu_0^* \otimes \eta), & \text{if } \nu \text{ is shift-invariant,}\\
\infty, & \text{otherwise,}
\end{array} \right.
\end{split}
\end{align}
for $\nu \in {\cal M}_1(\Sigma_b^+)$.
In this expression, $H$ denotes relative entropy and 
writing $\pi_k$ for the projection mapping from
$\R^{\N_0}$ to $\R^k$ given by
$(x_n)_{n \in \N_0} \mapsto (x_0, \dots, x_{k-1}),$
measures $\nu_i^*$ are defined as follows:
For $i \in \{0,1\}$ and shift-invariant $\nu \in {\cal M}_1(\Sigma_b^+)$,
we denote by $\nu_i^*$ the unique probability measure on $[b,0]^{\Z \cap (-\infty, i]}$ such that, for each $k \in \N$
and each Borel set $\Gamma \subseteq [b, 0]^k$,
$$
\nu_i^*( \{ ( \dots, x_{i-k+1}, \dots, x_i) : (x_{i-k+1}, \dots, x_i) \in \Gamma \}) = \nu \circ \pi_k^{-1} (\Gamma).
$$
Note that $\nu_i^*$ is well-defined due to the shift-invariance of $\nu$.
Furthermore, set
\begin{equation*} %\label{LMeasureDef}
L (\beta, \nu) := \int_{\Sigma_b^+} \log \E_1 \exp \Big \{ \int_0^{T_0} (\zeta(Y_s) + \beta) \, ds \Big \}  \, \nu (d \zeta)
\end{equation*}
for all $\beta \in \R$ and $\nu \in {\cal M}_1(\Sigma^+_b)$.
% Note that $L(\beta, \nu)$ is well-defined for all $\beta \in \R$ in this case since the integrand
% on the right-hand side of (\ref{LDef}) is bounded from below on $\Sigma_b^+$ and hence
% its negative part is integrable.
In particular, we have $L(\beta) = L(\beta, \Prob).$
% For $\nu \in {\cal M}_1 (\Sigma_b^+)$ we denote by $\beta_{cr}(\nu)$ the
% critical value such that $L(\beta, \nu) = \infty$ for $\beta > \beta_{cr}(\nu)$ and $L(\beta, \nu) < \infty$ for
% $\beta < \beta_{cr}(\nu).$ Furthermore, we set
% $\beta_{cr}^0 := \beta_{cr}(\delta_0^{\N_0}).$
Employing the notation
\begin{equation} \label{lambdaPSupDef}
L_p^{sup}(\beta) := \sup_{\nu \in {\cal M}_1^s(\Sigma_b^+)}
\Big(
L(\beta, \nu) - \frac{{\cal I}(\nu)}{p}
\Big), \quad \beta \in \R,
\end{equation}
we are ready to formulate
our main result for the annealed setting.
\begin{thm} \label{annealedLyapExp}
Assume (\ref{potEssSup}) and (\ref{potIIDBdd}).
Then, for each $p \in (0,\infty),$ the $p$-th annealed Lyapunov exponent $\lambda_p$ exists. Furthermore, $\lambda_p$ equals
the zero of $\beta \mapsto L_p^{sup}(-\beta)$ in $(-\beta_{cr}, 0)$ or,
if such a zero does not exist, equals $-\beta_{cr}.$
\end{thm}

\begin{remark}
 In fact, $L_p^{sup}$ has at most one zero as it is strictly increasing, cf. Lemma \ref{LpsupProps} below.
\end{remark}
With respect to the proof of this theorem,
it turns out that the asymptotics of the $p$-th moment $\langle u(t,0)^p \rangle$ is the same as
the quenched behaviour of $u(t,0)^p$ but under a different distribution of the environment $\xi.$
This will be made precise by the use of the aforementioned process level large deviations for $R_n \circ \xi.$

The term $L_p^{sup}$ defined in (\ref{lambdaPSupDef}) and
appearing in the above characterisation of $\lambda_p$ admits a convenient interpretation as follows.
On the one hand,
distributions $\nu$ of our random medium which provide us with high values of
$L(\beta, \nu)$ 
can play an important role in attaining the supremum in the right-hand side of
(\ref{lambdaPSupDef}).
On the other hand, we have to pay a price for obtaining such (rare) distributions, which is given by
${\cal I} (\nu) / p.$ As is heuristically intuitive and evident from formula (\ref{lambdaPSupDef}),
this price in relation to the 
gain obtained by high values of $L(\beta, \nu)$ becomes smaller as $p$ gets larger.
Note that, heuristically, $L(\cdot, \nu)$ corresponds to the function appearing in Theorem
\ref{quenchedLyapExp} characterising the quenched Lyapnov exponent for a potential distributed according
to $\nu.$

% The principal steps in the proof of Theorem \ref{annealedLyapExp} are similar to those of the proof of Theorem
% \ref{quenchedLyapExp} in the sense that we look for the summand yielding the largest contribution in the
% Feynman-Kac formula (\ref{feynmanKac}). This time, however, we are taking in addition
% the expectation with respect to the medium.
% Therefore, to estimate the corresponding expression from above we use a process-level large deviations
% principle for $(R_n \circ \xi)_{n \in \N}$ and then apply Varadhan's lemma
% in order to reduce it to expressions already known from the quenched case.
% The lower bound consists of restricting the $R_n \circ \xi$ to approximate some $\nu \in {\cal M}_1^s(\Sigma_b^+)$
% and then taking the supremum over such $\nu.$

As mentioned before, we are interested in the 
intermittency of $u$ for which we have the following result:
\begin{prop} \label{pIntermittency}
For $p > 0$ large enough, the solution $u$ to (\ref{PAMArb}) is $p$-intermittent.
\end{prop}

Furthermore, as mentioned previously, one expects the $p$-th annealed Lyapunov exponent $\lambda_p$
to converge to the quenched Lyapunov exponent $\lambda_0$ as $p \downarrow 0.$

\begin{thm} \label{lyapExpCont}
 $$\lim_{p \downarrow 0} \lambda_p = \lambda_0.$$
\end{thm}

% Due to the fact that we have expressions for $\lambda_p$ for all $p \in [0, \infty)$ instead of only for $p \in \N_0$, we can prove
% the continuity of $\lambda_p$ in $p:$
% \begin{prop} \label{lyapExpCont}
% Under the same assumptions as in Theorem \ref{annealedLyapExp},
% $
% [0, \infty) \ni p \mapsto \lambda_p
% $
% is continuous.
% \end{prop}
% In particular, note that this includes the heuristically plausible continuity of $\lambda_p$ in $p=0.$

\subsection{Related work}
The parabolic Anderson model without drift and i.i.d. or Gaussian potential is well-understood,
see the survey of G\"artner and K\"onig \cite{GaKo-05} as well as the references therein.
As a common feature, these treatments take advantage of
the self-adjointness of the random Hamiltonian $\kappa \Delta + \xi$
which allows for a spectral theory approach to the respective problems.
In our setting, however, the (random) operators $\kappa \Delta_h + \xi$ are not self-adjoint, whence we
do not have the common functional calculus at our disposal.
As hinted at earlier, we therefore retreat to large deviations principles of certain hitting times connected to
(\ref{feynmanKac}) in dimension one.
Heuristically, another difference caused by the drift is that
the drift term of the Laplace operator makes it harder for the random walk $X$
appearing in the Feynman-Kac representation (\ref{feynmanKac}) of the solution to stay at islands of
values of $\xi$ close to its supremum $0.$
% Since the main contributions to $u$ in the Feynman-Kac representation stem from paths of $X$
% favouring sites with values of $\xi$ close to $0,$ it is heuristically plausible that the Lyapunov exponents should be
% nonincreasing as a function in $h \in (0,1].$
% At least in the case $\lambda_p = -\beta_{cr}$ we have a proof for
% this statement, see Proposition $\ref{critBetaIncr}.$

Our model without drift
has been dealt with by Biskup and K\"onig in \cite{BiKo-01a} (not restricted to one dimension)
and \cite{BiKo-01}. Here the authors found
formulae for the quenched and $p$-th annealed Lyapunov exponents for all $p \in (0, \infty)$ using
a spectral theory approach.
Furthermore, they
investigated the so-called
screening effect that can appear in dimension one.

A situation similar in spirit to ours has been examined in the seminal article \cite{GrHo-92}
by Greven and den Hollander,
motivated from the point of view of population dynamics.
The model treated there is a discrete-time branching model in random environment with drift and
corresponds to the case of a bounded i.i.d. potential.
In particular, their motivation stems from the discrete-time analogue described at the end
of section \ref{notation}.
The results are formulated by the use of nontrivial
variational formulae.
While the authors concentrate on the explicit dependence of the results on the drift parameter $h,$
an advantage of our approach
is that we may compute the $p$-th annealed Lyapunov exponents for {\it all} $p \in (0, \infty)$
and characterise them in a simpler way.

In the context of discrete time, 
it is well worth mentioning recent results of Flury \cite{Fl-07}. Departing from a different set of questions,
he computed
the quenched and first annealed Lyapunov exponents and obtains large deviations
for discrete-time random walks with drift under the influence of a random potential in
arbitrary dimensions. Using Varadhan's lemma,
he derives the result from a large deviations principle by Zerner \cite{Ze-98} for the case without drift.
However, it is not clear how to apply the corresponding techniques to our situation. Firstly, as pointed out
in \cite{Ze-98}, this large deviations principle does not carry over
to the continuous-time case automatically, which also involves large deviations
on the number of jumps. Secondly,
% the continuous-time analogue of Flury's model corresponds to symmetric simple
% random walk $S$ under a measure of the form
% $$
% \frac{1}{Z^h_{t,\xi}} \exp \Big\{ \int_0^t \xi(S_s) \, ds + h \cdot S_t \Big\} \P,
% $$
% which does not exactly describe our problem for random walk with drift $h.$
% Indeed, proceeding along this approach
% one would also have to consider large deviations for the number of jumps of $S.$
% Furthermore, 
and more importantly, it is not clear how to adapt the methods of
\cite{Ze-98} and \cite{Fl-07} to obtain $\lambda_p$ for general $p \in (0,\infty),$ which is the
main focus of this paper.

\subsection{Outline}

Section \ref{auxiliaryResults} contains auxiliary results both for the quenched and annealed context.
The proofs of Theorem \ref{quenchedLyapExp} and Corollary \ref{optimalSpeed} will be carried out in section
\ref{proofsQuenched}. In section \ref{auxiliaryAnnealed} we prove some results needed for the proof
of Theorem \ref{annealedLyapExp}. The latter is then the subject of section \ref{proofsAnnealed},
while section \ref{furtherResults} contains some further properties of the Lyapunov exponents
as well as the proofs of Proposition \ref{pIntermittency} and Theorem \ref{lyapExpCont}.

While the results we gave in section \ref{mainResults} are valid for arbitrary $h \in (0,1],$ the corresponding
proofs in sections \ref{auxiliaryResults} to \ref{proofsAnnealed} contain steps which a priori hold true
for $h \in (0,1)$ only. Section \ref{maxDrift} deals with the adaptations necessary to
obtain their validity for $h=1$ also. Finally, in section \ref{maxDrift} we will also give a more
convenient representation for $\lambda_p$ with $p \in \N,$ see Proposition
\ref{maxDriftRep}.

\section{Auxiliary results} \label{auxiliaryResults}

In this section we prove auxiliary results which will primarily facilitate the proof of the quenched
results given in section \ref{mainResults}, but will also play a role when deriving the annealed results.

All of the results hereafter implicitly assume (\ref{potEssSup}) and (\ref{potErg}) mentioned in section \ref{notation}.

Departing from (\ref{feynmanKac}), the strong Markov property supplies us with
\begin{align}\label{decomp}
\begin{split}
\E_n & \exp \Big\{ \int_0^t \xi(Y_s) \, ds \Big\} \ind{1}_{Y_t = 0}\\
&= \E_n \Big( \exp \Big \{ \int_0^{T_0} \xi(Y_s) \, ds \Big\} \ind{1}_{ T_0 \leq t}
\Big( \E_0 \exp \Big\{ \int_0^{t-r} \xi(Y_s) \, ds \Big \} \ind{1}_{Y_{t-r} = 0} \Big)_{ r = T_0 } \Big).
\end{split}
\end{align}
The advantage of considering the time reversal of (\ref{feynmanKac}) is now apparent: 
For a fixed realisation of the medium, the term
$
\E_0 \exp \{ \int_0^{t-r} \xi(Y_s) \, ds\} \ind{1}_{Y_{t-r} = 0}
$
sees the same part of the medium, independent of which $n \in \Z$ the random walk $Y$ is starting from.

The main results of this section are Proposition \ref{clambda}, which controls the aforementioned term, and the large
deviations principle of Theorem
\ref{arbDriftLDP}, which helps to control the remaining part of the right-hand side in (\ref{decomp}).
The remaining
statements of this section are of a more technical nature.

The following result is motivated in spirit by
section VII.6 in \cite{Fr-85}. Note that we exclude the case of absolute drift $h=1.$

\begin{prop} \label{clambda}
\begin{enumerate}
 \item 
For $h \in (0, 1)$ and $x,y \in \Z$, the finite limit
\begin{align} \label{cStarDef}
c^* := \lim_{t \to \infty} \frac1t \log \E_x \exp \Big\{ \int_0^t \xi(Y_s) \, ds \Big\} \ind{1}_{Y_t = y}
\end{align}
exists a.s., equals 
\begin{equation} \label{cStarSiteIndep}
\sup_{t \in (0, \infty)} t^{-1} \log \E_0 \exp \Big \{ \int_0^t \xi(Y_s) \, ds \Big \} \ind{1}_{Y_t = 0},
\end{equation}
and is non-random.
Furthermore, $c^* \leq -\beta_{cr}$.
If either (\ref{potErgBdd}) or (\ref{potIID}) hold true, then
\begin{equation} \label{cStarBetaEq}
c^* = -\beta_{cr}.
\end{equation}

\item
For $\beta > -c^*,$
we have
\begin{equation}
\E_1 \exp \Big\{ \int_0^{T_0} (\xi(Y_s) + \beta) \, ds \Big\} = \infty \quad a.s.
\end{equation}
If either (\ref{potErgBdd}) or (\ref{potIID}) hold true, then for each $\beta < \beta_{cr}$ there
exists a non-random constant $C_\beta < \infty$ such that
\begin{equation} \label{LBoundedIntegrand}
\E_1 \exp \Big\{ \int_0^{T_0} (\xi(Y_s) + \beta) \, ds \Big\} \leq C_\beta \quad a.s.
\end{equation}

\end{enumerate}
\end{prop}

\begin{remark}
 Identity (\ref{cStarBetaEq}) will prove useful for the simplification of the variational
problems arising in sections \ref{proofsQuenched} and \ref{proofsAnnealed}.
\end{remark}

\begin{proof}
We start with the proof of (a) and split it into four steps.

$(i)$
We first show that for all $x,y \in \Z,$ the limit in (\ref{cStarDef}) exists and equals the expression in
(\ref{cStarSiteIndep}).

For $t \geq 0$ and $x, y \in \Z,$ define
$$
p_{x,y} (t) := \log \E_x \exp \Big\{ \int_0^t \xi(Y_s) \, ds \Big\} \ind{1}_{Y_t =y }.
$$
Using the Markov property, we observe that $p_{0,0}$ is super-additive. Therefore, the limit $c^*$ of
$p_{0,0}(t)/t$ as $t \to \infty$ exists and
\begin{equation} \label{cStarSup}
c^* = \sup_{t \in (0,\infty)} p_{0,0}(t)/t \in (-\infty,0].
\end{equation}
For $x, y \in \Z,$
the Markov property applied at times $1$ and $t+1$ yields
\begin{align}\label{constEst}
\begin{split}
\E_x &\exp \Big\{ \int_0^{t+2} \xi(Y_s) \, ds \Big\} \ind{1}_{Y_{t+2} =y } \\
&\geq \E_x \Big( \ind{1}_{ Y_s \in \{0 \wedge x, \dots, 0 \vee x\} \, \forall s \in [0,1], Y_1 =0 }
\exp \Big\{ \int_0^{t+2} \xi(Y_s) \, ds \Big\}\\
&\quad \times \ind{1}_{Y_{t+1} = 0, Y_s \in \{0 \wedge y, \dots, 0 \vee y\} \forall s \in[t+1,t+2], Y_{t+2} =y } \Big) \\
% &= \min_{k \in \{ 0 \wedge x, \dots, 0 \vee x\}} e^{\xi(k)}
% \times \P_x(Y_s \in \{ 0\wedge x, \dots, 0 \vee x\} \, \forall s \in [0,1], Y_1=0) \nonumber \\
% & \quad \times \E_0 \Big( \exp \Big\{ \int_0^{t} \xi(Y_s) \, ds \Big\} \ind{1}_{ Y_{t} = 0} \Big) \nonumber \\
% &\quad \times \min_{k \in \{0 \wedge y, \dots, 0 \vee y\}} e^{\xi(k)} 
% \times \P_0(Y_s \in \{ 0\wedge y, \dots, 0 \vee y \} \, \forall s \in[0,1], Y_1 =y ) \nonumber \\
&\geq c_{x,y} \E_0 \exp \Big\{ \int_0^{t} \xi(Y_s) \, ds \Big\} \ind{1}_{ Y_{t} = 0},
\end{split}
\end{align}
where $c_{x,y}$ is an a.s. positive random variable given by
\begin{align} 
c_{x,y} := & \min_{k \in \{ 0 \wedge x, \dots, 0 \vee x\}} e^{\xi(k)} 
\times \P_x(Y_s \in \{ 0\wedge x, \dots, 0 \vee x\} \, \forall s \in [0,1], Y_1=0) \nonumber \\
&\quad \times  \min_{k \in \{ 0 \wedge y, \dots, 0 \vee y\}} e^{\xi(k)} 
\times \P_0(Y_s \in \{ 0\wedge y, \dots, 0 \vee y \} \, \forall s \in[0,1], Y_1 =y ). \label{cxyConstDef}
\end{align}
Similarly,
\begin{equation} \label{reverseConst}
\E_0 \exp \Big\{ \int_0^{t+2} \xi(Y_s) \, ds \Big\} \ind{1}_{Y_{t+2} =0} 
\geq c_{y,x} \E_x \exp \Big\{ \int_0^{t} \xi(Y_s) \, ds \Big\} \ind{1}_{ Y_{t} = y} .
\end{equation}
Now, combining (\ref{cStarSup}) to (\ref{reverseConst}) we conclude that
$
\lim_{t \to \infty} p_{x,y}(t)/t
$
exists and equals (\ref{cStarSiteIndep}).

$(ii)$
We next show that $c^*$ is non-random and (\ref{cMEquality}) holds.

Naming the dependence of $c^*$ on the realisation explicitly, we obtain
$$
c^* = c^*(\xi) = c^*(\theta \circ \xi)
$$
by the use of (i). Thus, $c^*$ is non-random by Birkhoff's ergodic theorem.

In order to derive (\ref{cMEquality}), observe that for $M \in \N$ the function
$$
p_M (t) := \log \E_0 \exp \Big\{ \int_0^t \xi(Y_s) \, ds \Big\} \ind{1}_{T_{-M} > t, Y_t = 0}
$$
is super-additive. Hence
$
c_M^* := \lim_{t \to \infty} p_M(t)/t 
$
is well-defined and equals
$
\sup_{t \in (0, \infty)} p_M(t)/t.
$
Obviously, $p_M(t)$ is nondecreasing in $M$ and $p_M(t) \leq p_{0,0}(t),$ whence
\begin{equation} \label{cMSmallerC}
c_M^* \leq c^*
\end{equation}
for all $M \in \N.$
On the other hand, since 
$c_M^* \geq p_M(t)/t$
and
$p_M(t) \uparrow p_{0,0}(t)$ as $M \to \infty,$ we get
$\lim_{M \to \infty} c_M^* \geq p_{0,0}(t)/t$ for all $t$ and, consequently,
$
\lim_{M \to \infty} c_M^* \geq c^*.
$
Together with (\ref{cMSmallerC}) it follows that
\begin{equation} \label{cStarApprox}
\lim_{M \to \infty} c_M^* = c^*.
\end{equation}
Similarly to the previous step we compute
\begin{align*}
\E_0 \exp \Big\{ & \int_0^{t+2} \xi(Y_s) \, ds \Big\} \ind{1}_{T_{-M} > t+2, Y_{t+2} = 0} \\
% & \geq \P_0 (Y_s \in \{0,1\} \forall s \in [0,1], Y_1 = 1) \exp\{ \xi(0) \wedge \xi(1)\} \\
% &\quad \times \E_1 \exp \Big\{ \int_0^t \xi(Y_s) \, ds \Big\} \ind{1}_{T_{-M} > t, Y_t = 1} \\
% &\quad \times \P_1 (Y_s \in \{0,1\} \forall s \in [0,1], Y_1 = 0}) \exp\{ \xi(0) \wedge \xi(1)\}
&\geq c_{1,1} \E_1 \exp \Big\{ \int_0^t \xi(Y_s) \, ds \Big \} \ind{1}_{T_{-M} > t, Y_t = 1}\\
&= c_{1,1} \E_0 \exp \Big \{ \int_0^t (\theta \circ \xi)(Y_s) \, ds \Big \} \ind{1}_{T_{-(M+1)} > t, Y_t = 0}.
\end{align*}
Taking logarithms, dividing both sides by $t$ and letting $t$ tend to infinity,
we obtain
$
c^*_M(\xi) \geq c^*_{M+1}(\theta \circ \xi).
$
Iterating this procedure and using the monotonicity of
$c_M^*$ in $M,$ we obtain
$$
c_M^*(\xi) \geq \frac1n \sum_{j=1}^{k-1} c_{M+j}^* (\theta^j \circ \xi) + \frac1n \sum_{j=k}^n c^*_{M+k}(\theta^j \circ \xi)
$$
for all $n \in \N$ and $k \leq n$.
Birkhoff's ergodic theorem now yields $c_M^*(\xi) \geq \langle c_{M+k}^*(\xi) \rangle$
a.s. for all $k \in \N$. Because we clearly have
$
c_M^*(\xi) \leq c_{M+k}^*(\xi),
$
this implies that $c_M^*$ is constant a.s. and independent of $M.$
Due to (\ref{cStarApprox}) this gives $c_M^* = c^*$ a.s. for all $M \in \N,$
and thus
\begin{equation} \label{cMEquality}
 c^* = \lim_{t \to \infty} \frac1t \log \E_0 \exp \Big\{ \int_0^t \xi(Y_s)\, ds \Big\} \ind{1}_{T_{-M} > t, Y_t = 0}
\end{equation}
for all $M \in \N.$ By a similar derivation as above we find that (\ref{cMEquality}) holds for all
$M \in \Z \backslash \{0\},$ a fact which will prove useful in section \ref{proofsAnnealed}.

$(iii)$
The next step is to prove that $\beta_{cr} \leq -c^*$.

Given $t, \varepsilon > 0,$
we apply the Markov property at time $t$ to obtain
\begin{align*}
\E_0 &\exp \Big\{ \int_0^{T_{-1}} (\xi(Y_s) - c^* + \varepsilon) \, ds \Big\}\\
&\geq \E_0 \exp \Big\{ \int_0^{T_{-1}} (\xi(Y_s) - c^* + \varepsilon) \, ds \Big\} \ind{1}_{T_{-1} > t, Y_t = 0}\\
&= \E_0 \exp \Big\{ \int_0^t (\xi(Y_s) - c^* + \varepsilon) \, ds \Big\} \ind{1}_{T_{-1} > t, Y_t = 0}\\
&\quad \times \E_0 \exp \Big\{ \int_0^{T_{-1}} (\xi(Y_s) -c^* + \varepsilon) \, ds \Big\}.
\end{align*}
The second factor on the right-hand side is positive a.s. while, as we infer from (\ref{cMEquality}) for $M=1,$
the first one is logarithmically
equivalent to
$
e^{t\varepsilon}.
$
Thus, we deduce
\begin{equation} \label{bFirstPart}
\E_0 \exp \Big\{ \int_0^{T_{-1}} (\xi(Y_s) - c^* + \varepsilon) \, ds \Big\} = \infty \quad \text{a.s.,}
\end{equation}
and, using the shift invariance of $\xi,$ we get
\begin{equation} \label{bFirstPartShift}
\E_1 \exp \Big\{ \int_0^{T_{0}} (\xi(Y_s) - c^* + \varepsilon) \, ds \Big\} = \infty \quad \text{a.s.}
\end{equation}
In particular, this implies
$
L^+(-c^* + \varepsilon) = \infty
$
and thus $\beta_{cr} \leq -c^*.$

% PROBLEM, falls $L(0) = -\infty.$!!!!!!!!!!!11
% Now the second factor in the last equation is positive $\eta_0$-a.s. and due to our penultimate step in the proof,
% the first factor is exponentially equivalent to $\exp\{t(c^* - c^* + \varepsilon)\}$ which tends to $\infty$ as 
% $t \to \infty.$ It ensues that
% $$
% \Big \langle \Big( \log \Big( \E_0 \exp \Big\{ \int_0^{T_{-1}} (\xi(Y_s) - c_\varepsilon^*)ds \Big\} \Big) \Big \rangle = \infty.
% $$
% Since $\varepsilon > 0$ was arbitrarily chosen, this completes the proof of the Proposition.

$(iv)$
This part consists of showing that $\beta_{cr} \geq -c^*$ if either
(\ref{potErgBdd}) or (\ref{potIID}) is fulfilled.

% We first consider the case (\ref{potErgBdd}).
Note that the shift invariance of $\xi$ yields
$$
L^+(\beta) =
\Big \langle \Big(  \log \E_0 \exp \Big\{ \int_0^{T_{-1}} (\xi(Y_s) + \beta) \, ds\Big \} \Big)^+ \Big \rangle.
$$
Using (\ref{cStarSup}) as well as (\ref{reverseConst}) and taking into account that $c^* \leq 0,$ we get
\begin{equation} \label{upperSupAddEst}
\E_0 \exp \Big\{ \int_0^t \xi(Y_s) \, ds \Big\} \ind{1}_{Y_t =-1} \leq e^{c^* t}/c_{-1,0}
\end{equation}
for all $t \in (0,\infty).$ 
% Here,
% since $\xi$ is bounded from below by $b,$
% $c_{-1,0}$ of (\ref{cxyConstDef})
% can be chosen not to depend on the
% realisation of the medium anymore.
Consequently, we compute for $n \in \N$ and $\varepsilon > 0$:
\begin{align}
\E_0 &\exp \Big\{ \int_0^{T_{-1}} (\xi(Y_s) -c^* - \varepsilon) \, ds \Big\}
\ind{1}_{ T_{-1} \in (n-1, n], Y_s =-1 \, \forall s \in [T_{-1}, n] } \nonumber \\
&\leq \E_0 \exp \Big\{ \int_0^{n} \xi(Y_s) \, ds \Big\}
\ind{1}_{ Y_{n} = -1 } \exp\{-\xi(-1)\}
\exp\{ -c^*n - \varepsilon (n-1)\} \nonumber\\
&\leq \exp\{ - \varepsilon (n-1) - \xi(-1)\} /c_{-1,0} \quad \text{a.s.},\label{upperEstSplit}
\end{align}
where we have used $c^* \leq 0$ to deduce the first inequality and (\ref{upperSupAddEst}) to obtain the last one.
Analogously, the strong Markov property at time $T_{-1}$ supplies us with the lower bound
\begin{align}
\E_0 &\exp \Big\{ \int_0^{T_{-1}} (\xi(Y_s) -c^* - \varepsilon) \, ds \Big\}
\ind{1}_{ T_{-1} \in (n-1, n], Y_{s} =-1 \, \forall s \in [T_{-1}, n]} \nonumber\\
&\geq \E_0 \exp \Big\{ \int_0^{T_{-1}} (\xi(Y_s) - c^* - \varepsilon) \, ds \Big\}
\ind{1}_{T_{-1} \in (n-1, n]} \P_{-1} (Y_{s} =-1 \, \forall s \in [0,1]).\label{lowerEstSplit}
\end{align}
Since $\P_0(T_{-1} < \infty) = 1$,
combining (\ref{upperEstSplit}) with (\ref{lowerEstSplit}) and summing over $n \in \N,$ we get
\begin{equation} \label{partBSecondFirst}
\E_0 \exp \Big\{ \int_0^{T_{-1}} (\xi(Y_s) - c^* - \varepsilon) \, ds \Big\} 
\leq C \sum_{n \in \N} \exp \{ - \varepsilon n\} < \infty \quad \text{a.s.},
\end{equation}
where
\begin{equation} \label{univConst}
C := \big( \P_{-1} (Y_s =-1 \, \forall s \in [0,1]) \big)^{-1} \exp\{\varepsilon - \xi(-1)\} /c_{-1,0}.
\end{equation}

We now distinguish cases and first assume (\ref{potErgBdd}). In this case, $c_{-1,0}$ can
be bounded from below by some constant $\underline{c}_{-1,0} > 0$ a.s., whence
$C$ can be bounded from above by the non-random constant
\begin{equation} \label{overlineConst}
\overline{C} := \big( \P_{-1} (Y_s =-1 \, \forall s \in [0,1]) \big)^{-1} \exp\{\varepsilon - b \} /\underline{c}_{-1,0}.
\end{equation}
In particular, using (\ref{partBSecondFirst}) this implies $L^+ (-c^* - \varepsilon) < \infty,$
whence we deduce $\beta_{cr} \geq -c^*$.

To treat the second case assume (\ref{potIID}).
Due to (\ref{potEssSup}) and (\ref{potIID}), we infer
$
\Prob ( \xi(-1) \geq b) > 0
$
for any $b \in (-\infty, 0);$ fix one such $b.$
On
$\{ \xi(-1) \geq b\},$
as before, $C$ may be bounded from above by the corresponding constant $\overline{C}$ of (\ref{overlineConst})
and therefore 
\begin{equation} \label{partBSecondFirstUnderline}
\E_0 \exp \Big\{ \int_0^{T_{-1}} (\xi(Y_s) - c^* - \varepsilon) \, ds \Big\} 
\leq \overline{C} \sum_{n \in \N} \exp \{ - \varepsilon n\} < \infty
\end{equation}
$
\Prob (\cdot \vert \xi(-1) \geq b)
${-a.s.}
Since
the left-hand side of (\ref{partBSecondFirstUnderline}) does not depend on the actual value of $\xi(-1),$
the independence property
(\ref{potIID}) yields that (\ref{partBSecondFirstUnderline}) even holds $\Prob$-a.s., which
finishes the proof of part (a).

It remains to prove (b).
The first part was already established in (\ref{bFirstPartShift}).
Under assumption (\ref{potErgBdd}), the upper bound is a consequence of (\ref{partBSecondFirst})
with $C$ replaced by $\overline{C}$ of (\ref{overlineConst}); otherwise,
if (\ref{potIID}) is fulfilled, the upper bound follows from the last conclusion in the proof of part (a) (iv).
\end{proof}

Proposition \ref{clambda} enables us to control the asymptotics of the second expectation on the right
of (\ref{decomp}). To deal with the first expression, we define for $n \in \N$ and $\zeta \in \Sigma_b^+$
the probability measures
$$
\P^\zeta_n(A) := (Z^\zeta_n)^{-1} \E_n \exp \Big \{ \int_0^{T_0} \zeta(Y_s) \, ds \Big\} \ind{1}_A
$$
with $A \in {\cal F}$ and the normalising constant
$$
Z_n^\zeta := \E_n \exp \Big\{ \int_0^{T_0} \zeta(Y_s) \, ds \Big\}.
$$
The expectation with respect to $\P_n^\zeta$ will be denoted $\E_n^\zeta.$
By considering
$
\P_n^\xi \circ (T_0/n)^{-1},
$
we obtain a random sequence of probability measures on $\R_+$
for which we aim to prove a large deviations principle (see Theorem \ref{arbDriftLDP} below).
As common in the context of large deviations,
we define the moment generating function
$$
\Lambda(\beta) := \lim_{n \to \infty} \frac1n \log \E_n^\xi \exp\{ \beta T_0\} = \langle \log \E_1^\xi \exp \{ \beta T_0 \} \rangle,
\quad \beta \in \R,
$$
where the last equality stems from Birkhoff's ergodic theorem. Note that
\begin{equation} \label{LambdaLSum}
\Lambda (\beta) = L(\beta) - L(0)
\end{equation}
whenever the right-hand side is well-defined.

The following lemma tells us that the critical value $\beta_{cr}$ of $L^+$ also applies to $\Lambda$
and is positive.

\begin{lem} \label{lambdaCrit}
Assume $h \in (0, 1).$ Then
\begin{enumerate}
 \item 
$\Lambda(\beta) < \infty$ for $\beta < \beta_{cr}$, while $\Lambda(\beta) = \infty$ for $\beta > \beta_{cr};$

\item
$\beta_{cr}$ is positive.

\end{enumerate}
\end{lem}

\begin{remark}
Note that for $h = 1$ we can explicitly compute
$
c^* = \xi(0) - \kappa
$
as well as
$
\beta_{cr} = \kappa.
$
\end{remark}

\begin{proof}
$(a)$
 Since $Z_1^\xi \leq 1,$ we get
\begin{equation} \label{obvLambdaIneq}
\Big( \log \E_1 \exp \Big\{ \int_0^{T_0} (\xi(Y_s) + \beta) \, ds \Big \} \Big)^+
\leq \log \E_1^\xi \exp \{ \beta T_0 \}
\end{equation}
for $\beta \geq 0.$
Consequently, since $\beta_{cr} \geq 0,$
it is evident that $\Lambda(\beta) = \infty$ for $\beta > \beta_{cr}.$ 
For the remaining part of the statement, we estimate with $\beta \in [0, \kappa):$
\begin{align}
\E_1^\xi \exp \{ \beta T_0 \}
&\leq \frac{\E_1 \exp \big\{ \int_0^{T_0} (\xi(Y_s) + \beta) \, ds \big\} \ind{1}_{T_0 \leq T_2}}
{\E_1 \exp \big\{ \int_0^{T_0} \xi(Y_s) \, ds \big\} \ind{1}_{T_0 \leq T_2}} \nonumber\\
&\quad + \frac{\E_1 \exp \big\{ \int_0^{T_0} (\xi(Y_s) + \beta) \, ds \big\} \ind{1}_{T_2 \leq T_0}}
{\E_1 \exp \big\{ \int_0^{T_0} \xi(Y_s) \, ds \big\} \ind{1}_{T_0 \leq T_2}} \nonumber\\
&= \frac{\frac{1+h}{2} \frac{\kappa}{\kappa - \xi(1) - \beta} + \frac{1-h}{2} \frac{\kappa}{\kappa - \xi(1) - \beta}
\E_2 \exp \big\{ \int_0^{T_0} (\xi(Y_s) + \beta) \, ds \big\} }
{\frac{1+h}{2} \frac{\kappa}{\kappa - \xi(1)}} \nonumber \\
&\leq \frac{\kappa - \xi(1)}{\kappa - \xi(1) - \beta}
\Big( 1 + \E_2 \exp \Big\{\int_0^{T_0} (\xi(Y_s) + \beta) \, ds \Big\} \Big)\nonumber \\
&\leq \frac{\kappa}{\kappa - \beta}
\Big( 1 + \E_1 \exp \Big\{ \int_0^{T_0} ( (\theta \circ \xi)(Y_s) + \beta) \, ds \Big\}\nonumber \\
&\quad \times \E_1 \exp \Big\{ \int_0^{T_0} (\xi(Y_s) + \beta) \, ds \Big\} \Big) \quad \text{a.s.} \label{MVTBound}
\end{align}
Taking logarithms on both sides and using the inequality
$
\log(1+xy) \leq \log 2 + \log^+ x + \log^+ y
$
for $x,y> 0,$
we arrive at
\begin{align}
\log \E_1^\xi \exp \{\beta T_0\}
&\leq \log \frac{\kappa}{\kappa - \beta} + \log 2 +
\Big( \log \E_1 \exp \Big\{ \int_0^{T_0} ( (\theta \circ \xi)(Y_s) + \beta) \, ds \Big\} \Big)^+\nonumber \\
&\quad + \Big( \log \E_1 \exp \Big\{ \int_0^{T_0} (\xi(Y_s) + \beta) \, ds \Big\} \Big)^+. \label{inobvLambdaIneq}
\end{align}
We observe that for $\beta < \beta_{cr} (\leq \kappa$ due to (\ref{betaCrEst})$)$
the right-hand side is integrable with respect to $\Prob$ and hence
so is the left-hand side; thus,
$\Lambda(\beta) < \infty$ for $\beta < \beta_{cr}$.

$(b)$
It is sufficient to prove $\beta_{cr} > 0$ for the vanishing potential $\xi(x) \equiv 0$.
But in this case we have $\beta_{cr} = -c^*$ due to part (a) of Proposition \ref{clambda}. Therefore, using
the definition of $c^*,$ the proof reduces to a standard
large deviations bound and will be omitted. 

\end{proof}

We next prove the following properties of $L$ defined by (\ref{LDef}). Recall that $L$ is well-defined
for $\beta < \beta_{cr}.$

\begin{lem} \label{lDiff}
\begin{enumerate}
\item

If $L(\beta) > -\infty$ for some $\beta \in (-\infty, \beta_{cr}),$ then the same is true for
all $\beta \in (-\infty, \beta_{cr}).$

% \item
% $\Prob$-a.s.,
% the function
% $$
% (-\infty, \beta_{cr}) \ni \beta \mapsto \log \E_1 \exp \Big\{ \int_0^{T_0} (\xi(Y_s) + \beta) \, ds \Big\}
% $$
% is continuously differentiable on $(-\infty, \beta_{cr}).$

\item 
If the function $L$ is finite on $(-\infty, \beta_{cr}),$ then it
is continuously differentiable on this interval. Its
derivative is given by
\begin{equation} \label{LDerivative}
L'(\beta) = \Big \langle \frac{\E_1 T_0 \exp \{ \int_0^{T_0} (\xi(Y_s) + \beta) \, ds \}}
{\E_1 \exp \{ \int_0^{T_0} (\xi(Y_s) + \beta) \,ds \}} \Big \rangle.
\end{equation}

\item
If the assumptions of (b) apply, then
$$
\lim_{\beta \to -\infty} L'(\beta) = 0.
$$
\end{enumerate}
\end{lem}

\begin{remark}
In the situation that $L$ is finite, this lemma yields that $\Lambda'(\beta)$ is also given
by the expression in (\ref{LDerivative}) (cf. (\ref{LambdaLSum})).
\end{remark}

\begin{proof}
$(a)$
Assume $L(\beta) > -\infty$ for some $\beta \in (-\infty, \beta_{cr}).$
Due to the monotonocity of $L,$ it suffices to show $L(\beta-c) > -\infty$ for all $c > 0.$
We apply a reverse H\"older inequality for $q < 0 < r < 1$ with $q^{-1} + r^{-1} = 1$ to obtain
\begin{align*}
\E_1 &\exp \Big \{ \int_0^{T_0}(\xi(Y_s) + \beta -c) \, ds \Big\}\\
&= \E_1 \exp \Big\{ \int_0^{T_0} (\xi(Y_s) + \beta)/r \, ds \Big\} 
\exp \Big\{ \int_0^{T_0} ((\xi(Y_s) + \beta)/q - c) \, ds \Big\}\\
&\geq \Big( \E_1 \exp \Big\{ r \int_0^{T_0} (\xi(Y_s) + \beta)/r \, ds \Big\} \Big)^\frac1r
\Big( \E_1 \exp \Big\{ q \int_0^{T_0} ((\xi(Y_s) + \beta)/q -c) \, ds \Big\} \Big)^\frac1q\\
&= \Big( \E_1 \exp \Big\{ \int_0^{T_0} (\xi(Y_s) + \beta) \, ds \Big\} \Big)^\frac1r
\Big( \E_1 \exp \Big\{ \int_0^{T_0} (\xi(Y_s) + \beta - qc) \, ds \Big \} \Big)^\frac1q.
\end{align*}
Using the definition of $L,$ we obtain
\begin{equation*}
 L(\beta-c)
\geq \frac1r L(\beta) + \frac1q L(\beta - qc)
\end{equation*}
and for $\vert q \vert > 0$ small enough such that $\beta - qc < \beta_{cr},$ the second
summand is finite.
The first summand is finite by assumption, whence the claim follows.

$(b)$
The proof is standard and uses assertion (\ref{LBoundedIntegrand}) of Proposition \ref{clambda}. The details are left
to the reader.

$(c)$
Due to $(b)$ it is sufficient to show
that the integrand converges to $0$ pointwise and then apply dominated convergence\footnote{
Indeed, dominated convergence is applicable since the integrand is increasing in $\beta$ as one may check
by considering its derivative, and it is integrable for $\beta = 0,$ cf. e.g. (\ref{MVTBound}), Lemma \ref{lambdaCrit} (b)
and Proposition \ref{clambda} (b).
}
to infer the desired
result.

For this purpose fix a realisation of the medium and observe
$
\P_{1}^\xi \circ T_0^{-1} \ll \lambda
$
with $\lambda$ denoting the Lebesgue measure on $\R_+.$
Then the random density
$$
f:= \frac{ d \P_{1}^\xi \circ T_0^{-1} } {d \lambda}
$$
is well-defined.
It follows that
$$
\E_{1}^\xi T_0 \exp \{ \beta T_0\} = \int_{\R_+} x \exp\{ \beta x\} f(x) \, dx
$$
and splitting this integral we compute for $\varepsilon > 0$ and $\beta < 0$:
\begin{align} \label{densityIntEst}
\begin{split}
\int_0^\varepsilon x \exp \{ \beta x \} f(x) \, dx &\geq c \int_0^\varepsilon x \exp \{ \beta x \} \, dx
= c \Big( \frac{1}{\beta} x \exp\{ \beta x\} \vert_{x=0}^\varepsilon - \frac{1}{\beta^2} (\exp \{ \beta \varepsilon \} - 1 ) \Big)\\
&= c \Big( \frac{\varepsilon}{\beta} \exp \{ \beta \varepsilon \} + \frac{1}{\beta^2} - \frac{1}{\beta^2} \exp \{ \beta \varepsilon \} \Big),
\end{split}
\end{align}
where $c > 0$ is chosen such that $f \geq c$ holds $\lambda_{[0,\varepsilon]}$-a.s.
Similarly, the remaining part is estimated by
$$
\int_\varepsilon^\infty x \exp \{ \beta x \} f(x) \, dx \leq \varepsilon \exp \{\beta \varepsilon\} \int_\varepsilon^\infty f(x) \, dx
\leq \varepsilon \exp \{ \beta \varepsilon \}
$$
for $\beta < -\varepsilon^{-1}.$
Thus, for each $\varepsilon > 0$ we can choose $\beta$ small enough such that
$
\E_{1}^\xi T_0 \exp \{ \beta T_0\} \ind{1}_{T_0 \geq \varepsilon}
\leq \E_{1}^\xi T_0 \exp \{ \beta T_0\} \ind{1}_{T_0 \leq \varepsilon}
$
whence it follows for such $\beta$ that
$$
\E_{1}^\xi T_0 \exp \{ \beta T_0\} \leq 2 \varepsilon 
\E_{1}^\xi \exp \{ \beta T_0\}.
$$
This proves
that the above integrand converges to $0$ a.s. for $\beta \to -\infty$ and the result follows.
\end{proof}

In contrast to $L,$ the function $\Lambda$ may never take the value $-\infty,$ as is seen in the following lemma.
\begin{lem} \label{LambdaProp}
$\Lambda(\beta) > -\infty$ for all $\beta \in (-\infty, \beta_{cr}).$
\end{lem}

\begin{proof}
Due to monotonicity, it is sufficient to show $\Lambda(\beta) > -\infty$ for all $\beta \in (-\infty, 0).$ For this purpose
we choose such $\beta$
and estimate
\begin{align*}
 \frac{\E_1 \exp \{ \int_0^{T_0} (\xi(Y_s) + \beta) \, ds \}}{\E_1 \exp \{ \int_0^{T_0} \xi(Y_s) \, ds \}}
&\geq \frac{\E_1 \exp \{ \int_0^{T_0} (\xi(Y_s) + \beta) \, ds \} \ind{1}_{T_0 \leq T_2}}
{\E_1 \exp \{ \int_0^{T_0} \xi(Y_s) \, ds \} \ind{1}_{T_0 \leq T_2} (\P_1 (T_0 \leq T_2))^{-1}}\\
&= \frac{1+h}{2} \frac{\kappa - \xi(1)}{\kappa - \xi(1) - \beta}.
\end{align*}
Taking logarithms and expectations, we see that
$\Lambda(\beta) > -\infty$ for all $\beta \in (-\infty, 0).$
\end{proof}

We now have the necessary tools available to tackle the desired large deviations principle.
Let $\Lambda^*$ denote the Fenchel-Legendre transform of $\Lambda$ given by
$$
\Lambda^*(\alpha) := \sup_{\beta \in \R} (\beta \alpha - \Lambda(\beta)) = \sup_{\beta < \beta_{cr}} (\beta \alpha - \Lambda(\beta)),
\quad \alpha \in \R,
$$
where the second equality is due to Lemma \ref{lambdaCrit} (a).
Furthermore, for $M > 0$ and $n \in \N$, define
$$
\P^\xi_{M,n} = \P^\xi_n( \cdot \vert \{\tau_k \leq M \; \forall k \in \{1, \dots, n \}\})
$$
where 
$\tau_k := T_{k-1} -T_k,$
$k \in \N.$
The corresponding expectation is denoted by $\E_{M,n}^\xi$ 

\begin{thm} \label{arbDriftLDP}
For almost all realisations of $\xi,$
the sequence of probability measures $(\P_n^\xi \circ (T_0/n)^{-1} )_{n \in \N}$ on $\R_+$ satisfies
a large deviations principle on scale $n$ with
deterministic, convex
good rate function $\Lambda^*.$
\end{thm}

\begin{proof}
Being the supremum of affine functions, $\Lambda^*$ is lower semi-continuous and convex. Furthermore,
since $\Lambda(0) = 0,$ it follows that
$\Lambda^*(u) \geq 0$ for all $u \in \R$.
Choosing $\beta \in (0, \beta_{cr}),$ which is possible due to Lemma \ref{lambdaCrit} (b),
we find that for any $M \geq 0$ the set
$$
\{\alpha \in \R : \beta \alpha - \Lambda(\beta) \leq M \}
\cap
\{\alpha \in \R : -\beta \alpha - \Lambda(-\beta) \leq M \}
$$
is compact and, in particular, $\Lambda^*$ has compact level sets; thus, $\Lambda^*$ is a good convex rate function.

The
upper large deviations bound for closed sets is a direct consequence of the
G\"artner-Ellis theorem (cf. Theorem 2.3.6 in \cite{DeZe-98}).

To prove the lower large deviations bound for open sets,
we cannot directly apply the G\"artner-Ellis theorem since the steepness assumption (cf. Definition 2.3.5 (c) in \cite{DeZe-98})
is possibly not fulfilled.
Indeed, if $h = 1$ it may occur that
$
\lim_{\beta \uparrow \beta_{cr}} \vert \Lambda^{'}(\beta) \vert < \infty
$
since in this case $\beta_{cr} = \kappa,$
\begin{equation} \label{lambdaPrimeExplicit}
\Lambda'(\beta) = \Big \langle \frac{\E_1^\xi T_0 \exp \{ \beta T_0\}}{\E_1^\xi \exp\{ \beta T_0\}} \Big \rangle
= \Big \langle \frac{1}{\kappa - \beta -\xi(0)} \Big \rangle
\end{equation}
and $\Lambda$ is steep if and only if
$
-1/\xi(0)
$
is not integrable.

To circumvent this problem, we retreat to the measures $\P_{M,n}^\xi$ and for the corresponding
logarithmic moment generating function
we write
$$
\Lambda_M (\beta) := \lim_{n \to \infty} \frac1n \log \E_{M,n}^\xi \exp \{ \beta T_0\}
= \langle \log \E_{M,1}^\xi \exp \{ \beta T_0\} \rangle, \quad \beta \in \R,
$$
where the equality is due to Birkhoff's ergodic theorem. Using dominated convergence, one checks that
$\Lambda_M$ is essentially smooth (cf. Definition 2.3.5 in \cite{DeZe-98}). We may therefore apply the
G\"artner-Ellis theorem to the sequence $(\P_{M,n}^\xi \circ (T_0/n)^{-1})_{n \in \N}$ to obtain for any
$G \subseteq \R_+$ open and $x \in G$ the estimate
\begin{equation} \label{BddGELowerBd}
\liminf_{n \to \infty} \frac1n \log \P_{M,n}^\xi \circ (T_0 /n)^{-1} (G)
\geq -\Lambda_M^*(x),
\end{equation}
where $\Lambda_M^*(\alpha) := \sup_{\beta \in \R} (\beta \alpha - \Lambda_M (\beta))$
denotes the Fenchel-Legendre transform for $\alpha \in \R.$

In order to use this result for our original problem, we recall that a sequence of functions
$(f_n)_{n \in \N}$ from $\R$ to $\overline{\R}$
{\it epi-converges} to a function $f: \R \to \overline{\R}$ at $x_0 \in \R$
if and only if
$$
\liminf_{n \to \infty} f_n(x_n) \geq f(x_0)
$$
for all sequences $(x_n)_{n \in \N} \subset \R$ converging to $x_0$ and
$$
\limsup_{n \to \infty} f_n(x_n) \leq f(x_0)
$$
for some sequence $(x_n)_{n \in \N} \subset \R$ converging to $x_0.$

Using the facts that
\begin{enumerate}
\item $\Lambda$ is continuous on $(-\infty, \beta_{cr}),$

\item $\Lambda_M \to \Lambda$ pointwise as $M \to \infty,$

\item $\Lambda_M$ is a monotone function and the sequence $(\Lambda_M)_{M \in \N}$ is monotone when
restricted either to $(-\infty, 0]$ or $[0, \infty),$
\end{enumerate}
we deduce that $\Lambda_M$ epi-converges towards $\Lambda$ as $M \to \infty.$
Therefore, since we note that $\Lambda$ (cf. Lemma \ref{LambdaProp}) and the $(\Lambda_M)_{M \in \N}$ are proper,
lower semi-continuous and convex functions, we conclude using Theorem 11.34 in \cite{RoWe-98} that
$\Lambda_M^*$ epi-converges towards $\Lambda^*$ as $M \to \infty$ along $\N.$ 
Choosing $G$ and $x$ as above we therefore find a sequence $(x_M)_{M \in \N} \subset G$ with
$\lim_{M \to \infty} \Lambda^*_M(x_M) = \Lambda^*(x).$ Employing (\ref{BddGELowerBd}) we thus obtain
$$
\limsup_{M \to \infty} \liminf_{n \to \infty} \frac1n \log \P_{M,n}^\xi (T_0 / n)^{-1} (G) \geq -\Lambda^*(x),
$$
which in combination with 
$$
\P_n^\xi \circ (T_0 / n )^{-1} (G) \geq \P_{M,n}^\xi (T_0 / n)^{-1} (G) \cdot \P_M^\xi (\tau_k \leq M \, \forall k \in \{1, \dots, n\})
$$
yields
$$
\liminf_{n \to \infty} \frac1n \log \P_n^\xi \circ (T_0 /n)^{-1} (G) \geq - \Lambda^*(x).
$$
This finishes the proof of the lower bound.

\end{proof}

\section{Proofs for the quenched regime} \label{proofsQuenched}
The outline of this section is as follows.
We will use the large deviations principle of Theorem \ref{arbDriftLDP}
to derive a variational formula for the lower logarithmic
bound of $u$ (cf. Lemma \ref{quenchedLowerBd}). In combination with further estimates, the large deviations
principle will also prove valuable in establishing a similar estimate for the upper bound, see Lemma
\ref{quenchedUpperBd}.
Combining Lemmas \ref{quenchedLowerBd} and \ref{quenchedUpperBd} we obtain Corollary
\ref{quenchedLyapExpVar} and can then complete the proof of Theorem \ref{quenchedLyapExp}.
Both the lower and upper bounds will essentially depend on the fact that
(\ref{feynmanKac}) can be used to obtain
(\ref{singleSummand}).
This result will then be shown explicitly to yield the proof of Corollary \ref{optimalSpeed}.
Here, $\asymp$ means exponential equivalence.

We start with the proof of the lower bound.
\begin{lem} \label{quenchedLowerBd}
Let $I \subset [0, \infty)$ be an interval. Then for almost all realisations of $\xi,$
\begin{equation} \label{quenchedLowerBdFormula}
 \liminf_{t \to \infty} \frac1t \log \sum_{n \in tI} \E_n \exp \Big\{ \int_0^t \xi(Y_s) \, ds \Big\} \ind{1}_{Y_t = 0}
\geq \sup_{\alpha \in \stackrel{\circ}{I} \cup (I\cap \{0\})} \inf_{\beta < \beta_{cr}}
(-\beta + \alpha L(\beta)).
\end{equation}
\end{lem}

\begin{proof}
For $\delta > 0$ and $\alpha \in \stackrel{\circ}{I}$ we obtain using (\ref{decomp}):
\begin{align*}
&\E_\floor{\alpha t} \exp \Big\{ \int_0^t \xi(Y_s) \, ds \Big\} \ind{1}_{Y_t = 0}\\
&\geq \E_\floor{\alpha t} \Big( \exp \Big\{ \int_0^{T_0} \xi(Y_s) \, ds \Big\} \ind{1}_{ \frac{T_0}{\floor{\alpha t}} \leq \frac{1}{\alpha} }
\Big( \E_0 \exp \Big\{ \int_0^{t-r} \xi(Y_s) \, ds \Big\} \ind{1}_{Y_{t-r} = 0} \Big)_{r = T_0} \Big)\\
&\geq \sup_{m \in (0,1/\alpha)} \E_\floor{\alpha t} \Big( \exp \Big\{ \int_0^{T_0} \xi(Y_s) \, ds \Big\}
\ind{1}_{\frac{T_0}{\floor{\alpha t}} \in (0, 1/\alpha) \cap (m-\delta, m + \delta)}\\
&\quad \times\Big( \E_0 \exp \Big \{ \int_0^{t-r} \xi(Y_s) \, ds \Big\} \ind{1}_{Y_{t-r} = 0} \Big)_{r = T_0} \Big).
\end{align*}
Applying
Proposition \ref{clambda} (a) to the inner expectation and
Theorem \ref{arbDriftLDP} to the remaining part of the right-hand side, we have
\begin{align}
\liminf_{t \to \infty} & \frac1t \log \E_\floor{\alpha t} \exp \Big\{ \int_0^t \xi(Y_s) \, ds \Big\} \ind{1}_{Y_t = 0}\nonumber \\ 
&\geq  \sup_{m \in (0,1/\alpha)} \Big( \alpha \big(  -\inf_{x \in (0, 1/ \alpha) \cap (m-\delta, m + \delta)} \Lambda^*(x) + L(0)\big)
+ (1- \alpha (m - \delta)) c^* \Big).\label{quenchedLowerFirstIEst}
\end{align}
Note here that for a potential unbounded from below, $L(0)=-\infty$ is possible;
% however, reading through the proof
% of Proposition \ref{arbDriftLDP}, we can compensate $L(0)$ with the normalising factors $Z_{M,n}^\xi$ in the limit
% to obtain for $\delta \downarrow 0$ using the lower semi-continuity of $I$:
nevertheless, observe that in this case also $L(\beta) = -\infty$ for all $\beta < \beta_{cr},$ 
see Lemma \ref{lDiff} (a).
% Indeed,
% \begin{align*}
%  L(\beta) &= \Big \langle \E_1 \exp \Big\{ \int_0^{T_0} (\xi(Y_s) + \beta) \, ds \Big\} \Big \rangle\\
% & \leq \Big \langle \Big( \E_1 \exp \Big\{ \int_0^{T_0} p\xi(Y_s) \, ds \Big\} \Big)^\frac1p \Big\rangle\\
% & \quad + \log \Big( \E_1 \exp \Big\{q \beta T_0 \Big\} \Big)\frac1q
% \end{align*}
% with $p^{-1} + q^{-1} = 1.$ Choosing $q > 0$ small enough, the second summand on the right-hand side is finite, while
% the first one equals $-\infty$ since $p > 1.$
Therefore, the following computations hold true even if $L(0) = - \infty.$
The lower semi-continuity of $\Lambda^*$ supplies us with 
$
\lim_{\delta \downarrow 0} \inf_{x \in (0, 1/\alpha) \cap (m-\delta, m + \delta)} \Lambda^*(x) = \Lambda^*(m).
$
Hence, taking $\delta \downarrow 0$ in (\ref{quenchedLowerFirstIEst}) yields
\begin{align}
\liminf_{t \to \infty} & \frac1t \log \E_\floor{\alpha t} \exp \Big\{ \int_0^t \xi(Y_s) \, ds \Big\} \ind{1}_{Y_t = 0}\nonumber \\
&\geq  \sup_{m \in (0,1/\alpha)}
\big ( \alpha (-\Lambda^*(m) + L(0)) + (1-\alpha m)c^* \big ) \nonumber \\
&= c^* + \alpha \sup_{m \in (0, 1/\alpha)} \inf_{\beta < \beta_{cr}}
( m(\beta_{cr}-\beta)+ L(\beta)),
\label{artificialLabel}
\end{align}
where we used (\ref{LambdaLSum}) and (\ref{cStarBetaEq}) to obtain the equality.
Thus, the supremum in $m$ is taken for $m=1/\alpha.$ Hence,
\begin{align}
\liminf_{t \to \infty} & \frac1t \log \E_\floor{\alpha t} \exp \Big\{ \int_0^t \xi(Y_s) \, ds \Big\} \ind{1}_{Y_t = 0}
\geq \inf_{\beta < \beta_{cr}} ( -\beta + \alpha L(\beta) ). \label{lowerBDVariationalFormulaI}
\end{align}
Now for the case $\alpha = 0$ we observe
$$
\liminf_{t \to \infty} \frac1t \log \E_0 \exp \Big\{ \int_0^t \xi(Y_s) \, ds \Big\} \ind{1}_{Y_t = 0} = - \beta_{cr}
$$
due to Proposition \ref{clambda} (a). Since
$
\inf_{\beta < \beta_{cr}} ( -\beta + \alpha L(\beta) )
$
evaluates to $-\beta_{cr}$ for $\alpha = 0,$ in combination with (\ref{lowerBDVariationalFormulaI}) this
finishes the proof of (\ref{quenchedLowerBdFormula}).
\end{proof}

Next, we turn to the upper bound which is slightly more involved.
\begin{lem} \label{quenchedUpperBd}
Let $I \subset [0, \infty)$ be a compact interval. Then for almost all realisations of $\xi,$
\begin{equation}\label{quenchedUpperBdFormula}
 \limsup_{t \to \infty} \frac1t \log \sum_{n \in tI} \E_n \exp \Big\{ \int_0^t \xi(Y_s) \, ds \Big\} \ind{1}_{Y_t = 0}
\leq \sup_{\alpha \in I} \inf_{\beta < \beta_{cr}}
(-\beta + \alpha L(\beta)).
\end{equation}
\end{lem}

\begin{proof}
$(i)$
First, assume that $\inf I > 0$ and write $I = [\varepsilon, \gamma].$ 
Then for
$\delta > 0$ we choose numbers $(\alpha^\delta_k)_{k=1}^{n}$ such that
$\varepsilon  = \alpha^\delta_1 < \alpha^\delta_2 < \dots < \alpha^\delta_{n} = \gamma$
and $\max_{k=1, \dots, n-1} (\alpha^\delta_{k+1} - \alpha^\delta_k) < \delta.$
Then
\begin{align*}
\limsup_{t \to \infty} &\frac1t \log \sum_{n \in tI}
\E_n \exp \Big\{ \int_0^t \xi(Y_s) \, ds \Big\} \ind{1}_{Y_t = 0}\\
 & = \max_{k=1, \dots, n-1} \limsup_{t \to \infty} \frac1t \log
\sum_{n \in t[\alpha_k^\delta, \alpha_{k+1}^\delta]} \E_n \exp \Big\{\int_0^t \xi(Y_s) \, ds \Big \} \ind{1}_{Y_t = 0}.
\end{align*}
Using (\ref{decomp}) and Proposition \ref{clambda} (a) we get
\begin{align}\label{quenchedFirstUpper}
\begin{split}
&\sum_{n \in t[\alpha_k^\delta, \alpha_{k+1}^\delta]}
\E_n \exp \Big\{\int_0^t \xi(Y_s) \, ds \Big \} \ind{1}_{Y_t = 0} \\
&=\sum_{n \in t[\alpha_k^\delta, \alpha_{k+1}^\delta]}
\E_n \Big( \exp \Big\{ \int_0^{T_0} \xi(Y_s) \, ds \Big\} \ind{1}_{T_0 \leq t} \Big( \E_0 \exp \Big\{ \int_0^{t-r} \xi(Y_s) \, ds \Big\}
\ind{1}_{Y_{t-r} = 0} \Big)_{r=T_0} \Big) \\
& \leq \sum_{n \in t[\alpha_k^\delta, \alpha_{k+1}^\delta]}
\E_n \exp \Big\{ \int_0^{T_0} \xi(Y_s) \, ds \Big\} \ind{1}_{\frac{T_0}{n} \leq \frac{1}{\alpha_k^\delta}}
\exp\{ c^*(t-T_0)\},
%\\
% &\asymp e^{c^*t} \E_\floor{\alpha_j^\delta t} \exp \Big\{ \int_0^{T_0} (\xi(Y_s) - c^*) \, ds \Big\}
% \ind{1}_{\frac{T_0}{\floor{\alpha_j^\delta t}} \leq \frac{1}{\alpha_k^\delta}},
\end{split}
\end{align}
which by (\ref{cStarBetaEq}) and the exponential Chebyshev inequality can be bounded from above by
\begin{align} \label{quenchedSecondUpper}
\begin{split}
&\inf_{\beta > 0} \sum_{n \in t[\alpha_k^\delta, \alpha_{k+1}^\delta]}
\E_n \exp \Big\{ \int_0^{T_0} (\xi(Y_s) - c^* - \beta) \, ds \Big\} \exp\{ c^*t\} \exp \{ \beta n / \alpha_k^\delta\}\\
&\leq
\exp\{ c^*t\} \inf_{\beta < \beta_{cr}} \sum_{n \in t[\alpha_k^\delta, \alpha_{k+1}^\delta]}
\E_n \exp \Big\{ \int_0^{T_0} (\xi(Y_s) +\beta) \, ds \Big\} \exp \{ (-\beta -c^*) n/ \alpha_k^\delta\}.
\end{split}
\end{align}
Therefore, combining (\ref{quenchedFirstUpper}) and (\ref{quenchedSecondUpper}) we arrive at
\begin{align} \label{quenchedThirdUpper}
 \begin{split}
\limsup_{t \to \infty} &\frac1t \log \sum_{n \in t[\alpha_k^\delta, \alpha_{k+1}^\delta]}
\E_n \exp \Big\{\int_0^t \xi(Y_s) \, ds \Big \} \ind{1}_{Y_t = 0}\\
& \leq c^* \Big(1 - \frac{\alpha_{j_k}^\delta}{\alpha_k^\delta} \Big)
+\alpha_{j_k}^\delta \inf_{\beta < \beta_{cr}} (-\beta/\alpha_k^\delta + (\Lambda(\beta) + L(0))),
\end{split}
\end{align}
where $j_k= k$ if the summands on the right-hand side of (\ref{quenchedSecondUpper})
have nonpositive exponential rates in $n$ for some $\beta>0$ and $j_k=k+1$ otherwise.
Now if $L(0) = -\infty$, then obviously the right-hand side of (\ref{quenchedThirdUpper}) equals $-\infty$ and
(\ref{quenchedUpperBdFormula}) holds true.
Therefore, we assume
$% \label{LBiggerMinusInf}
 L(0) > -\infty
$
from now on, which due to Lemma \ref{lDiff} (a) implies $L(\beta) > -\infty$ for all $\beta \in (-\infty, \beta_{cr}).$
By (\ref{LambdaLSum}),
the right-hand side of (\ref{quenchedThirdUpper})
evaluates to 
\begin{align}
% VORHERIGE VERSION
% c^* + \alpha_j^\delta & \sup_{m \in [0, 1/\alpha_k^\delta ]} \inf_{\beta < \beta_{cr}} 
% ( -c^*m -\beta m + L(\beta) )\nonumber \\
c^* \Big( 1 - \frac{\alpha_{j_k}^\delta}{\alpha_k^\delta} \Big)
+ \frac{\alpha_{j_k}^\delta}{\alpha_k^\delta} \inf_{\beta < \beta_{cr}} ( -\beta + \alpha_k^\delta L(\beta) ). \label{piecewiseUpperBd}
\end{align}
% where to obtain the equality we observed that due to (\ref{cStarBetaEq}), the supremum in $m$ is taken for $m = 1/\alpha_k^\delta.$
Using Lemma \ref{lDiff} (c)
one can show that the family of functions indexed by $\delta,$ which are defined piecewise constant by the
right-hand side of (\ref{piecewiseUpperBd}) for $\alpha \in [\alpha_k^\delta, \alpha_{k+1}^\delta),$
$k \in \{1, \dots, n-2\},$ and $\alpha \in [\alpha_{n-1}^\delta, \alpha_n^\delta],$ converges
uniformly in $\alpha \in [\varepsilon, \gamma]$ to
\begin{equation} \label{infBeta}
\inf_{\beta < \beta_{cr}} (-\beta + \alpha L(\beta))
\end{equation}
as $\delta \downarrow 0.$
% IN DISS IN DISS IN DISS IN DISS IN DISS IN DISS IN DISS IN DISS
% Indeed, this is not obvious only for the term $\inf_{\beta < \beta_{cr}} ( -\beta + \alpha_k^\delta L(\beta) ).$
% But differentiating 
% \begin{equation} \label{betaDiffTerm}
% -\beta + \alpha L(\beta)
% \end{equation}
% with respect to $\beta \in (-\infty, \beta_{cr}),$
% Lemma \ref{lDiff} $(b)$ (applicable due to Lemma \ref{lDiff} $(a)$ and since we assumed $L(0) > -\infty$ above)
% supplies us with
% $$
% -1 + \alpha \Big \langle \frac{ \E_1^\xi T_0 \exp \{ \beta T_0 \}}{ \E_1^\xi \exp \{ \beta T_0 \}} \Big \rangle.
% $$
% Lemma \ref{lambdaChoice} $(a)$ yields that this derivative is negative for some $\beta = \beta_*$ small enough
% and all $\alpha \in [\varepsilon, \gamma].$
% Hence, since (\ref{betaDiffTerm}) is convex in $\beta,$ this implies that the infimum
% in (\ref{infBeta}) can be restricted to $[\beta_*, \beta_{cr})$
% for all $\alpha \in [\varepsilon, \gamma]$ simultaneously.
% Setting $s: = \sup \{ L \leq \beta_{cr}/\varepsilon\},$ we may restrict the infimum in $\beta$ even to
% $[\beta_*, s).$ 
% The desired uniform convergence now follows, since
% $$
% (\alpha, \beta) \mapsto -\beta + \alpha_k^\delta L(\beta)
% $$
% for $\alpha \in [\alpha_k^\delta, \alpha_{k+1}^\delta]$ converges uniformly on $[\varepsilon, \gamma] \times [\beta_*, s).$
% 
% IN DISS IN DISS IN DISS IN DISS IN DISS IN DISS IN DISS IN DISS
Taking $\delta \downarrow 0$ and the supremum over $\alpha \in [\varepsilon, \gamma],$ we therefore
obtain from the previous relations:
\begin{align}
\limsup_{t \to \infty} \frac1t \log \sum_{n \in t[\varepsilon, \gamma]} \E_n \exp \Big\{ \int_0^t \xi(Y_s) ds \Big\}
\ind{1}_{Y_t = 0}
&\leq \sup_{\alpha \in [\varepsilon, \gamma]} \inf_{\beta < \beta_{cr}} 
(-\beta + \alpha L(\beta) ). \label{quenchedPosVel}
\end{align}

$(ii)$
It remains to consider the case that $\inf I = 0.$
Then we either find $\varepsilon > 0$ such that
\begin{align}
\limsup_{t \to \infty} & \frac1t \log
\sum_{n \in tI} \E_n \exp \Big \{ \int_0^t \xi(Y_s) \, ds \Big\} \ind{1}_{Y_t = 0} \nonumber \\
& = \limsup_{t \to \infty} \frac1t \log
\sum_{n \in t(I \cap [\varepsilon, \infty))} \E_n \exp \Big \{ \int_0^t \xi(Y_s) \, ds \, \Big\} \ind{1}_{Y_t = 0} \label{linMass},
\end{align}
in which case the problem reduces to the previous case
and (\ref{quenchedUpperBdFormula}) holds true in particular.
Otherwise, for each $\varepsilon >0$
we have ``$>$'' in (\ref{linMass}) instead of ``$=$''.
We would then
find a function $\varphi: [0,\infty) \to \N_0$ such that $\varphi(t)/t \to 0$ as $t \to \infty$ and
which satisfies the first equality in
\begin{align}
\limsup_{t \to \infty} & \frac1t \log
\sum_{n \in tI} \E_n \exp \Big \{ \int_0^t \xi(Y_s) \, ds \Big\} \ind{1}_{Y_t = 0} \nonumber \\
&= \limsup_{t \to \infty} \frac1t \log
\E_{\varphi(t)} \exp \Big \{ \int_0^t \xi(Y_s) \, ds \Big\} \ind{1}_{Y_t = 0} \nonumber \\
&= \limsup_{t \to \infty} \frac1t \log
\E_{\varphi(t)} \Big( \exp \Big\{ \int_0^{T_0} \xi(Y_s) \, ds \Big\} \ind{1}_{T_0 \leq t} \nonumber \\
&\quad \times \Big( \E_0 \exp \Big\{ \int_0^{t-r} \xi(Y_s) \, ds \Big\} \ind{1}_{Y_{t-r} = 0} \Big)_{r = T_0} \Big) \nonumber \\
&\leq \sup_{\alpha \in [\delta ,1 - \delta]} \Big(
\limsup_{t \to \infty} \frac{1}{t} \log
\E_{\varphi(t)} \Big( \exp \Big\{ \int_0^{T_0} \xi(Y_s) \, ds \Big\} \ind{1}_{T_0 \in 
t[\alpha-\delta, \alpha + \delta]} \Big) \nonumber \\
& \quad \times \sup_{r \in t[\alpha-\delta, \alpha + \delta]}
\E_0 \exp \Big\{ \int_0^{t-r} \xi(Y_s) \, ds \Big\} \ind{1}_{Y_{t-r} = 0}
\Big) \label{subLinearMassGrowth}
\end{align}
with $\delta> 0$ small.
The exponential Chebyshev inequality for $\beta \in (0, \beta_{cr})$ supplies us with
\begin{align} \label{quenchedDegenUpperBd}
\begin{split}
\E_{\varphi(t)} \exp \Big \{ & \int_0^{T_0} \xi(Y_s) \, ds \Big\} \ind{1}_{T_0 \geq (\alpha - \delta)t}\\
&\leq \E_{\varphi(t)} \exp \Big\{ \int_0^{T_0} (\xi(Y_s) + \beta) \, ds \Big\} \exp \{ -\beta (\alpha-\delta)t \}.
\end{split}
\end{align}
Taking $\delta \downarrow 0$ and $\beta \uparrow \beta_{cr}$ in this inequality, we deduce
\begin{align*}
\limsup_{t \to \infty} \frac1t \log
\E_{\varphi(t)} \exp \Big \{ \int_0^{T_0} \xi(Y_s) \, ds \Big\} \ind{1}_{T_0 \geq (\alpha - \delta)t}
\leq - \alpha \beta_{cr}. 
\end{align*}
Now taking $\delta \downarrow 0$ in (\ref{subLinearMassGrowth}) we obtain in combination with
Proposition \ref{clambda} (a):
\begin{align}
\limsup_{t \to \infty} & \frac1t \log
\sum_{n \in tI} \E_n \exp \Big \{ \int_0^t \xi(Y_s) \, ds \Big\} \ind{1}_{Y_t = 0} \nonumber \\
&\leq \sup_{\alpha \in [0,1]} (-\alpha \beta_{cr} + (1-\alpha)c^*) = -\beta_{cr}. \label{lambdaTilde}
\end{align}
But $-\beta_{cr}$ is just the result when replacing $\sup_{\alpha \in I}$ by
$\alpha = 0$ in (\ref{quenchedUpperBdFormula}). This finishes
the proof of the lemma.
\end{proof}
The next result establishes the intuitively plausible fact that only summands in the direction of the drift
are relevant on an exponential scale.

\begin{lem} \label{quenchedLeftSmallerRightExp}
For all $\delta \geq 0$ and $\gamma \geq \delta$ we have
 \begin{align}
\limsup_{t \to \infty} &\frac1t \log \sum_{n \in t[-\gamma , -\delta ]} \E_n \exp \Big\{ \int_0^t \xi(Y_s) \,ds\Big\}
\ind{1}_{Y_t = 0} \nonumber \\
&\leq \delta \log \frac{1-h}{1+h} + 
\liminf_{t \to \infty} \frac1t \log \sum_{n \in t[\delta, \gamma]} \E_n \exp \Big\{ \int_0^t \xi(Y_s) \, ds \Big\}
\ind{1}_{Y_t = 0} \quad \text{a.s.}\label{quenchedOneSidedConc}
\end{align}
% If the right-hand side is bigger than $-\beta_{cr},$ then (\ref{oneSidedConc}) holds with a strict inequality.
\end{lem}
\begin{proof}
% For $n \in \N$ we denote by ${\cal W}_n(1,0)$ the set of discrete-time $n$-step simple random walks on $\Z$
% starting in $1$ at time $0,$ terminating in $0$ and time $n$ and being strictly larger than $0$ at all times
% in between.
% Using the decomposition 
% \begin{align*} 
% \E_1 \exp & \Big\{ \int_0^{T_0} \xi(Y_s) \, ds \Big\}
% = \sum_{n \in \N_0} \Big(\frac{1+h}{2}\Big)^{n+1} \Big( \frac{1-h}{2} \Big)^n
% \sum_{s \in {\cal W}_{2n+1} (1,0)} \prod_{k=0}^{2n} \frac{\kappa}{\kappa - \xi(s_k)},
% \end{align*}
% we get the second equality in

First we observe
that the function
$$
(0,\infty) \ni \alpha \mapsto \inf_{\beta < \beta_{cr}} (-\beta + \alpha L(\beta))
$$
is either constant $-\infty$ (if and only if $L \equiv -\infty$ on $(-\infty, \beta_{cr}),$ cf. Lemma \ref{lDiff} (a))
or continuous, since it is concave.
% \footnote{
% Indeed, as an infimum over continuous functions, the function (denoted $g$ now)
% is upper semi-continuous.
% Now if there existed $\alpha^* \in (0,\infty)$ at which it was not lower semi-continuous, then we
% would find $\varepsilon > 0$ and a monotone sequence $\alpha_n \to \alpha^*$ (say increasing) such that
% \begin{equation} \label{nonLowerSemicont}
% g(\alpha_n) \leq g(\alpha^*) - \varepsilon
% \end{equation}
% for all $n \in \N.$ Thus, there would exist $\beta_n \in (-\infty, \beta_{cr})$
% such that $-\beta_n + \alpha_n L(\beta_n) \leq g(\alpha_n) + \varepsilon/2$ as well as
% $-\beta_n + \alpha L(\beta_n) \geq g(\alpha^*).$ But taking $n \to \infty$ and using
% (\ref{nonLowerSemicont}) this would imply $g(\alpha) = -\infty$ for $\alpha < \alpha^*.$ But since
% $\beta \mapsto -\beta + \alpha L(\beta)$ is bounded from below on $(-\infty, \beta_{cr})$
% (use Lemma \ref{lDiff} (c) and convexity), this yields a contradiction.
% }
In combination with the proofs of Lemmas \ref{quenchedLowerBd} and \ref{quenchedUpperBd} we therefore infer the existence of
$\alpha_+ \in [\delta, \gamma]$ such that
\begin{align*}
\lim_{t \to \infty} \frac1t \log \sum_{n \in t[\delta, \gamma]} \E_n \exp \Big\{ \int_0^t \xi(Y_s) \, ds \Big\}
\ind{1}_{Y_t = 0}
&=\lim_{t \to \infty} \frac1t \log \E_\floor{\alpha_+ t} \exp \Big\{ \int_0^t \xi(Y_s) \, ds \Big\}
\ind{1}_{Y_t = 0}\\
&= \sup_{\alpha \in [\delta,\gamma]} \inf_{\beta < \beta_{cr}} (-\beta + \alpha L(\beta)) \quad \text{a.s.}
\end{align*}
Employing similar arguments, one may show that the analogues of
Lemmas \ref{quenchedLowerBd} and \ref{quenchedUpperBd} for $I \subset (-\infty,0]$
also hold and we infer the existence of $\alpha_- \in [-\gamma, -\delta]$ such
that 
\begin{equation}\label{negSumEq}
\lim_{t \to \infty} \frac1t \log \sum_{n \in t[-\gamma , -\delta ]} \E_n \exp \Big\{ \int_0^t \xi(Y_s) \,ds\Big\}
\ind{1}_{Y_t = 0}
=\lim_{t \to \infty} \frac1t \log \E_\floor{\alpha_- t} \exp \Big\{ \int_0^t \xi(Y_s) \, ds \Big\}
\ind{1}_{Y_t = 0}
\end{equation}
exists and is deterministic also.
Next we observe that for $n \in \N$
\begin{align} \label{eqInDistribution}
\begin{split}
\E_{-n} \exp \Big\{ \int_0^t \xi(Y_s) \, ds \Big\} \ind{1}_{Y_t = 0}
&= \E_0 \exp \Big\{ \int_0^t \xi(X_s) \, ds \Big\} \ind{1}_{X_t = -n}\\
&=\Big( \frac{1-h}{1+h} \Big)^n \E_0 \exp \Big\{ \int_0^t \xi(Y_s) \, ds \Big\} \ind{1}_{Y_t = -n}\\
&\stackrel{d}{=} \Big( \frac{1-h}{1+h} \Big)^n \E_n \exp \Big\{ \int_0^t \xi(Y_s) \, ds \Big\} \ind{1}_{Y_t = 0},
\end{split}
\end{align}
where the first equality follows from time reversal, the second by comparing the transition probabilities
of $X$ and $Y,$ and
the last equality follows from the shift invariance of $\xi.$

Employing (\ref{eqInDistribution}) in combination with (\ref{negSumEq}) we conclude
\begin{align*}
 \lim_{t \to \infty} \frac1t \log &\sum_{n \in t[-\gamma , -\delta ]} \E_n \exp \Big\{ \int_0^t \xi(Y_s) \,ds\Big\}
\ind{1}_{Y_t = 0}\\
&=\lim_{t \to \infty} \frac1t \log \E_\floor{\alpha_- t} \exp \Big\{ \int_0^t \xi(Y_s) \, ds \Big\}
\ind{1}_{Y_t = 0}\\
&= \vert \alpha_- \vert \frac{1-h}{1+h} +
\lim_{t \to \infty} \frac1t \log \E_{-\floor{\alpha_- t}} \exp \Big\{ \int_0^t \xi(Y_s) \, ds \Big\}
\ind{1}_{Y_t = 0}\\
&\leq  \vert \alpha_- \vert \frac{1-h}{1+h} +
\lim_{t \to \infty} \frac1t \log \sum_{n \in t[\delta, \gamma]} \E_n \exp \Big\{ \int_0^t \xi(Y_s) \, ds \Big\}
\ind{1}_{Y_t = 0}.
\end{align*}
Indeed, to justify the last equality note that the limits on both sides exist and are constant
a.s.; (\ref{eqInDistribution}) then yields the equality in question.
This finishes the proof.
\end{proof}

\begin{lem} \label{quenchedFarAwaySummands}
We have
\begin{equation} \label{crucialGamma}
\limsup_{t \to \infty} \frac1t \log
\sum_{n \notin t[-\gamma, \gamma]} \E_n \exp \Big \{ \int_0^t \xi(Y_s) \, ds \Big\} \ind{1}_{Y_t = 0}
\to -\infty
\end{equation}
as $\gamma \to \infty$.
\end{lem}

\begin{proof}
Note that by the use of Stirling's formula we obtain for $\gamma > \kappa e$:
\begin{align*}
\sum_{n \geq \gamma t} \E_n \exp \Big \{ \int_0^t \xi(Y_s) \, ds \Big\} \ind{1}_{Y_t = 0}
&\leq \sum_{n \geq \gamma t} \P_n (T_0 \leq t)
\leq \sum_{n \geq \gamma t} \sum_{k \geq n} e^{-\kappa t} \frac{(\kappa t)^k}{k!}\\
&\leq e^{-\kappa t} \sum_{n \geq \gamma t} \sum_{k \geq n} \Big( \frac{\kappa t e}{\gamma t} \Big)^k\\
&=Ce^{-\kappa t} \sum_{n \geq \gamma t} \Big( \frac{\kappa e}{\gamma} \Big)^n
= Ce^{-\kappa t} \Big( \frac{\kappa e}{\gamma} \Big)^\floor{\gamma t}
\end{align*}
where $C$ is a generic constant depending on $\kappa$ and $\gamma,$ swallowing all sums appearing in the geometric series. 
Since an analogous result is valid for 
$$
\sum_{n \leq -\gamma t} \E_n \exp \Big \{ \int_0^t \xi(Y_s) \, ds \Big\} \ind{1}_{Y_t = 0},
$$
we infer that
(\ref{crucialGamma})
holds
as $\gamma \to \infty$.
\end{proof}

\begin{corollary} \label{quenchedLyapExpVar}
 The quenched Lyapunov exponent $\lambda_0$ exists and is given by
\begin{equation} \label{quenchedLyapExpVarFormula}
\lambda_0 = \sup_{\alpha \in [0,\gamma]} \inf_{\beta < \beta_{cr}}
(-\beta + \alpha L(\beta))
\end{equation}
for all $\gamma > 0$ large enough.
\end{corollary}
\begin{proof}
We take advantage of (\ref{feynmanKac}) to split for $\gamma > 0:$
\begin{align}
 u(t,0) &\leq \sum_{ n \notin t[-\gamma , \gamma ]} \E_n \exp \Big\{ \int_0^t \xi(Y_s) \, ds \Big\} \ind{1}_{Y_t = n} \nonumber \\
&\quad + \sum_{ n \in [0, \gamma t]} \E_n \exp \Big\{ \int_0^t \xi(Y_s) \, ds \Big\} \ind{1}_{Y_t = n} \label{quenchedFeynmanKacSplit}\\
&\quad + \sum_{ n \in [-\gamma t,0]} \E_n \exp \Big\{ \int_0^t \xi(Y_s) \, ds \Big\} \ind{1}_{Y_t = n}. \nonumber
\end{align}
Lemma \ref{quenchedLeftSmallerRightExp} yields that the third summand is logarithmically negligible
when compared to the second. Since the first summand can be made arbitrarily small for $\gamma$ large,
according to Lemma \ref{quenchedFarAwaySummands}, we obtain in combination with Lemma \ref{quenchedUpperBd}:
$$
\limsup_{t \to \infty} \frac1t \log u(t,0) \leq \sup_{\alpha \in [0,\gamma]} \inf_{\beta < \beta_{cr}}
(-\beta + \alpha L(\beta))
$$
for $\gamma$ large enough.

With respect to the lower bound,
Lemma \ref{quenchedLowerBd} in combination with (\ref{feynmanKac}) supplies us with 
$$
\liminf_{t \to \infty} \frac1t \log u(t,0) \geq \sup_{\alpha \geq 0} \inf_{\beta < \beta_{cr}}
(-\beta + \alpha L(\beta)).
$$
Combining these two estimates we infer the existence of $\lambda_0$ and the variational formula
(\ref{quenchedLyapExpVarFormula}).
\end{proof}

We are now ready to prove the results of subsection \ref{subsecQuenchedRegime}.
\begin{proof}[{\bf Proof of Theorem \ref{quenchedLyapExp}}]
Corollary \ref{quenchedLyapExpVar} supplies us with the existence of $\lambda_0$ and the variational formula
(\ref{quenchedLyapExpVarFormula}).
If $L$ does not have a zero in $(0, \beta_{cr}),$ then we have
$L(\beta) < 0$ for all $\beta < \beta_{cr}.$ Thus, the supremum in $\alpha$ is taken in $\alpha = 0$
and the right-hand side of (\ref{quenchedLyapExpVarFormula})
evaluates to $-\beta_{cr}.$ If $L$ does have a zero in $(0, \beta_{cr}),$
then inspecting (\ref{quenchedLyapExpVarFormula}) and differentiating with respect to $\beta,$
% we observe that $\lambda_0$ equals the zero of $L(-\cdot),$
% which finishes the proof of the theorem.
we observe that the supremum over $\alpha$ is a maximum taken in $\alpha = (L'(\beta_z))^{-1},$ 
with $\beta_z$ denoting the zero of $L$ in $(0, \beta_{cr}).$ Consequently, we deduce
that $\lambda_0$ equals $-\beta_z,$ which finishes the proof.
\end{proof}

\begin{proof} [\bf Proof of Corollary \ref{optimalSpeed}]

$(a)$
Note that $L$ has a zero in $(0, \beta_{cr})$ by assumption and thus Theorem \ref{quenchedLyapExp} implies
$
-\lambda_0 < \beta_{cr}.
$
Therefore, by
Lemma \ref{lDiff} (b) we may deduce
 $(L'(-\lambda_0))^{-1} \in (0, \infty).$
% NICHT LOESCHEN!
% IN DISS IN DISS IN DISS IN DISS IN DISS IN DISS IN DISS IN DISS IN DISS
% 
% The fact that
% $$
% \Big \langle \frac{\E_1 T_0 \exp \{ \int_0^{T_0} (\xi(Y_s) - \lambda_0) \, ds \}}
% {\E_1 \exp \{ \int_0^{T_0} (\xi(Y_s) -\lambda_0) \, ds\} } \Big \rangle > 0
% $$
% (IS OBVIOUS since the integrand is nonnegative and positive on a set of positive measure OR)
% follows, since by assumption we have $-\lambda_0 < \beta_{cr}$ and using \ref{partBSecondFirst} it ensues
% that the denominator is bounded from above a.s. for $h \in (0,1)$ by a constant.
% For $h=1,$ a direct computation yields the same result, cf. (\ref{lambdaPrimeExplicit}).
% 
% The finiteness follow from Lemma \ref{lDiff}.
% Thus, $\alpha^* \in (0, \infty).$
% 
% NICHT LOESCHEN!
% IN DISS IN DISS IN DISS IN DISS IN DISS IN DISS IN DISS IN DISS IN DISS

Using the time reversal of (\ref{feynmanKac}) and Lemma \ref{quenchedFarAwaySummands}, it suffices to show
\begin{align}
\lambda_0 > \limsup_{t \to \infty} \frac1t \log 
\sum_{n \in t ([-\gamma, \gamma] \backslash (\alpha^*-\varepsilon, \alpha^* +\varepsilon))}
\E_n \exp \Big\{ \int_0^t \xi(Y_s) \, ds \Big\}
\ind{1}_{Y_t =0} \label{firstOptSpeedToShow}
\end{align}
for $\gamma$ large enough.
First, observe that due to Lemma \ref{quenchedUpperBd} we have
\begin{align}
\limsup_{t \to \infty} \frac1t &\log \sum_{n \in t([0,\gamma] \backslash (\alpha^* - \varepsilon, \alpha^* + \varepsilon))}
\E_n \exp \Big\{ \int_0^t \xi(Y_s) \, ds \Big\}
\ind{1}_{Y_t = 0} \nonumber \\
& \leq
\sup_{\alpha \in [0,\gamma] \backslash (\alpha^* - \varepsilon, \alpha^* + \varepsilon)} \inf_{\beta < \beta_{cr}}
(-\beta + \alpha L(\beta)). \label{posSummandsOptSpeedEst}
\end{align}
Differentiating the expression
\begin{equation} \label{varFormula}
-\beta + \alpha L(\beta)
\end{equation}
with respect to $\beta$ we
obtain
\begin{equation} \label{varFormulaDerivative}
-1 + \alpha \Big \langle \frac{\E_1^\xi T_0 \exp\{\beta T_0\}}{\E_1^\xi \exp\{ \beta T_0\}} \Big \rangle,
\end{equation}
cf. Lemma \ref{lDiff} (b).
% Setting this expression equal $0$ for $\beta = -\lambda_0$ and keeping in mind the convexity of 
% (\ref{lowerBDVariationalTermI}) in $\beta$ we observe that this expression has a minimum in $\beta = -\lambda_0$
% if and only if
% $$
% \alpha = \Big \langle \frac{\E_1^\xi T_0 \exp\{ -\lambda_0 T_0\}}{\E_1^\xi \exp\{ -\lambda_0 T_0\}} \Big \rangle^{-1}.
% $$
% In particular, Theorem (\ref{quenchedLyapExp}) in combination with (\ref{lowerBDVariationalFormulaI}) yields
% \begin{align*}
% \lim_{t \to \infty} \frac1t \log u(t,0) &= \sup_{\alpha > 0} \inf_{\beta < \beta_{cr}} (-\beta + \alpha L(\beta))
% = \inf_{\beta < \beta_{cr}} \Big( -\beta + 
% \Big \langle \frac{\E_1^\xi T_0 \exp\{ -\lambda_0 T_0 \}}{\E_1^\xi \exp\{-\lambda_0 T_0\}} \Big \rangle^{-1}
% L(\beta) \Big)\\
% &= \lim_{t \to \infty} \frac1t \log \E_0 \exp \Big\{ \int_0^t \xi(X_s) \, ds \Big\} \ind{1}_{X_t \in (\floor{\alpha t}- \delta,\floor{\alpha t}+\delta)}
% \end{align*}
Now (\ref{varFormulaDerivative}) as a function of $\beta$ is continuous at $-\lambda_0$ and
inserting $\beta = -\lambda_0$ as well as $\alpha = \alpha^*,$ the term in (\ref{varFormulaDerivative}) evaluates to 0.
% Indeed, since the numerator is bounded from above by a finite constant a.s. for all $\beta \leq \beta_{cr}$ due to Proposition \ref{clambda} $(c)$
% and the denominator is larger than $1$, the Dominated Convergence Theorem supplies us with the desired continuity.
Therefore, for $\varepsilon \in (0, \alpha^*)$ there exists $\delta >0$ such that for all
$\alpha$ with $\vert \alpha - \alpha^* \vert \geq \varepsilon$ 
and $\beta$ with $\vert \beta - (-\lambda_0) \vert < \delta,$
the derivative (\ref{varFormulaDerivative}) is bounded away
from $0.$
Since,
according to Theorem \ref{quenchedLyapExp}, setting $\beta = -\lambda_0$ in (\ref{varFormula})
evaluates to $\lambda_0$ independently of the value of $\alpha,$
this boundedness yields
$$
\inf_{\beta < \beta_{cr}} (-\beta + \alpha L(\beta)) \leq \lambda_0 - \delta^*
$$
for some $\delta^*>0$ and all $\alpha \notin (\alpha^*-\varepsilon, \alpha^* + \varepsilon).$
Consequently, we get
$$
\sup_{\alpha \in [0,\gamma] \backslash (\alpha^*-\varepsilon, \alpha^* + \varepsilon)} \inf_{\beta < \beta_{cr}} (-\beta + \alpha L(\beta))
\leq \lambda_0 - \delta^* < \lambda_0.
$$
Therefore, using (\ref{posSummandsOptSpeedEst}) we have
\begin{align*}
\limsup_{t \to \infty} \frac1t &\log \sum_{n \in t([0,\gamma] \backslash (\alpha^* - \varepsilon, \alpha^* + \varepsilon))}
\E_n \exp \Big\{ \int_0^t \xi(Y_s) \, ds \Big\}
\ind{1}_{Y_t = 0} < \lambda_0.
\end{align*}
Combining this estimate with Lemma \ref{quenchedLeftSmallerRightExp} thus yields
\begin{align}
\limsup_{t \to \infty} \frac1t \log \sum_{n \in t( [-\gamma, 0]\backslash (-\alpha^* -\varepsilon, -\alpha^* + \varepsilon))}
\E_n \exp \Big\{ \int_0^{T_0} \xi(Y_s) \, ds \Big\} \ind{1}_{Y_t = 0} < \lambda_0. \label{negSummandsOptSpeedEst}
\end{align}
Furthermore, the same lemma supplies us with the first inequality in
\begin{align} \label{negatedOptSpeedEst}
\begin{split}
\limsup_{t \to \infty} \frac1t &\log \sum_{n \in t( -\alpha^* -\varepsilon, -\alpha^* + \varepsilon)}
\E_n \exp \Big\{ \int_0^{T_0} \xi(Y_s) \, ds \Big\} \ind{1}_{Y_t = 0}\\
&\leq (\alpha^* - \varepsilon) \log \frac{1-h}{1+h} 
+ \sup_{\alpha \in [0,\gamma] \backslash (\alpha^* - \varepsilon, \alpha^* + \varepsilon)} \inf_{\beta < \beta_{cr}}
(-\beta + \alpha L(\beta)) 
< \lambda_0.
\end{split}
\end{align}
Combining (\ref{posSummandsOptSpeedEst}), (\ref{negSummandsOptSpeedEst}) and (\ref{negatedOptSpeedEst}) gives
(\ref{firstOptSpeedToShow}) and hence finishes the proof of part (a).

$(b)$
To show the equality, observe that the assumption $\lim_{\beta \uparrow \beta_{cr}} L(\beta) = 0$ implies
that $L$ has no zero in $(-\beta_{cr}, 0)$ and hence
$\lambda_0 = - \beta_{cr}$ due to Theorem \ref{quenchedLyapExp}.
Now choose $m \in [0, (\lim_{\beta \uparrow \beta_{cr}} L'(\beta))^{-1}]$ and $\varepsilon \in (0,m).$
Employing time reversal and Lemma \ref{quenchedLowerBd} we arrive at
\begin{align}\label{constrainedLowerBd}
\begin{split}
 \liminf_{t \to \infty} &\frac1t \log \E_0 \exp \Big\{ \int_0^t \xi(X_s) \, ds \Big\} \ind{1}_{X_t \in t(m-\varepsilon, m+\varepsilon)}\\
&= \liminf_{t \to \infty} \frac1t \log \sum_{n \in t(m-\varepsilon,m+\varepsilon)}
\E_n \exp \Big\{ \int_0^t \xi(Y_s) \, ds \Big\} \ind{1}_{Y_t = 0}\\
&\geq \sup_{\alpha \in (m-\varepsilon, m+\varepsilon)} \inf_{\beta < \beta_{cr}} (-\beta + \alpha L(\beta))
= \inf_{\beta < \beta_{cr}} ( -\beta + (m-\varepsilon) L(\beta)),
\end{split}
\end{align}
where the last equality follows since $L(\beta) < 0$ for $\beta < \beta_{cr}.$
Differentiating the inner term of the right-hand side with respect to $\beta$ yields
$-1 + (m-\varepsilon) L'(\beta)$ which due to our choice of $m$ is smaller than $0$ for all $\beta < \beta_{cr}.$
Thus, (\ref{constrainedLowerBd}) evaluates to $-\beta_{cr},$ and this proves the desired equality.

To prove the inequality, observe that
Lemmas
\ref{quenchedUpperBd}, \ref{quenchedLeftSmallerRightExp} and \ref{quenchedFarAwaySummands} yield
\begin{align*}
\limsup_{t \to \infty} &\frac1t \log \sum_{n \leq -t\varepsilon}
\E_n \exp \Big\{ \int_0^t \xi(Y_s) \, ds \Big\} \ind{1}_{Y_t = 0}
\leq 
\varepsilon \log \frac{1-h}{1+h} +
\sup_{\alpha \geq \varepsilon} \inf_{\beta < \beta_{cr}} (-\beta + \alpha L(\beta)).
\end{align*}
Using Corollary \ref{quenchedLyapExpVar} we have that the right-hand side is strictly 
smaller than $\lambda_0.$

For the remaining summands,
Lemma \ref{quenchedUpperBd} in combination with Lemma \ref{quenchedFarAwaySummands} yields
\begin{align}\label{optSpeedBPosSumUpper}
\begin{split}
\limsup_{t \to \infty} &\frac1t \log \sum_{n \geq t((\lim_{\beta \uparrow \beta_{cr}} L'(\beta))^{-1} + \varepsilon)}
\E_n \exp \Big\{ \int_0^t \xi(Y_s) \, ds \Big\} \ind{1}_{Y_t = 0}\\
&\leq 
\sup_{\alpha \geq (\lim_{\beta \uparrow \beta_{cr}} L'(\beta))^{-1} + \varepsilon} \inf_{\beta < \beta_{cr}} (-\beta + \alpha L(\beta)),
\end{split}
\end{align}
and the derivative $-1 + \alpha L'(\beta)$ of the inner term on the right-hand side
with respect to $\beta$ is positive and bounded away from $0$ for all
$\beta < \beta_{cr}$ large enough and all $\alpha \geq (\lim_{\beta \uparrow \beta_{cr}} L'(\beta))^{-1} +\varepsilon.$
Thus, we conclude that
the right-hand side of (\ref{optSpeedBPosSumUpper}) is strictly smaller than $\lambda_0= -\beta_{cr},$ which
finishes the proof.

$(c)$
Using time reversal we get
for $\gamma > \varepsilon:$
\begin{align*}
\limsup_{t \to \infty} &\frac1t \log \E_0 \exp \Big\{ \int_0^t \xi(X_s) \, ds \Big\} \ind{1}_{X_t \notin t(-\varepsilon,\varepsilon)}\\
& = \limsup_{t \to \infty} \frac1t \log \sum_{n \notin t(-\varepsilon,\varepsilon)}
\E_n \exp \Big\{ \int_0^t \xi(Y_s) \, ds \Big\} \ind{1}_{Y_t = 0}\\
& = \limsup_{t \to \infty} \frac1t \log \sum_{n \notin t(-\gamma,\gamma)}
\E_n \exp \Big\{ \int_0^t \xi(Y_s) \, ds \Big\} \ind{1}_{Y_t = 0}\\
& \quad \vee \limsup_{t \to \infty} \frac1t \log 
\sum_{n \in t((-\gamma,\gamma) \backslash (-\varepsilon,\varepsilon))}\E_n \exp \Big\{ \int_0^t \xi(Y_s) \, ds \Big\}\ind{1}_{Y_t= 0}
\end{align*}
According to Lemma \ref{quenchedFarAwaySummands}, the first term on the right-hand side
tends to $-\infty$ as $\gamma \to \infty,$
while Lemma \ref{quenchedLeftSmallerRightExp} combined with
Lemma \ref{quenchedUpperBd} implies that the second can be estimated from above by
$$
\sup_{\alpha \in [\varepsilon, \gamma]} \inf_{\beta < \beta_{cr}} (-\beta + \alpha L(\beta))
= \inf_{\beta < \beta_{cr}} (-\beta + \varepsilon L(\beta)),
$$
where the equality follows since $L(\beta) < 0$ for all $\beta < \beta_{cr}$ by assumption.
Thus, this expression evaluates to $-\beta_{cr} + \varepsilon \lim_{\beta \uparrow \beta_{cr}} L(\beta) < -\beta_{cr} = \lambda_0,$
and the statement follows.
\end{proof}

\section{Auxiliary results particular to the annealed regime} \label{auxiliaryAnnealed}
This section contains mainly technical results, which will be employed in the proof of Theorem
\ref{annealedLyapExp}. The results given here are in parts
generalisations of corresponding results for a finite state space given in section $IX.2$ and $A.9$
of \cite{El-85}.
% 
% The previous lemma now tells us that we may apply Varadhan's Lemma (see e.g. Theorem 4.3.1 in \cite{DeZe-98}) to obtain in combination with
% taking $\delta \downarrow 0$ in (\ref{upperChebyParamEst}):
% \begin{align*}
% \limsup_{t \to \infty} \frac1t \log \langle u(t,0)^p \rangle
% \leq \sup_{\alpha \geq 0} \inf_{\beta \geq 0} \Big(p(\beta + c^*) + \alpha \sup_{\nu \in {\cal M}_1(\Sigma)} 
% (p\Lambda_\nu (-\beta - c^*) - {\cal I}(\nu)) \Big);
% \end{align*}
% considering that ${\cal I}$ is a infinite for non-stationary $\nu$, it ensues that
% \begin{equation} \label{statSup}
% \limsup_{t \to \infty} \frac1t \log \langle u(t,0)^p \rangle 
% \leq \sup_{\alpha \geq 0} \inf_{\beta \geq -\tilde{\beta_{cr}}} \Big(p\beta + \alpha \sup_{\nu \in {\cal M}^s_1(\Sigma_b)}
% (p \Lambda_\nu (-\beta) - I(\nu) ) \Big).
% \end{equation}
% It is our aim to restrict the supremum in $\nu \in {\cal M}^s_1(\Sigma_b)$ to $\nu \in {\cal M}_1^e(\Sigma_b)$, since our treatment of the lower 
% bound needs ergodicity.
% For this purpose we now show that $\nu \mapsto \phi_\beta(\nu) - {\cal I}(\nu)$ is an affine function, which follows if we prove that
% $
% \nu \mapsto {\cal I}(\nu)
% $
% defines an affine function.

\begin{lem} \label{relEntropLem}
Let $n \in \N$, $\rho \in {\cal M}_1(\R)$ and $\nu \in {\cal M}_1^s (\R^{\N_0}).$
Then
\begin{equation} \label{relEntropRep}
\sum_{i=1}^n H(\pi_i \nu \vert \pi_{i-1}\nu \otimes \rho)
= H(\pi_n \nu \vert \rho^n).
\end{equation}
\end{lem}
\begin{proof}
 The result is a consequence of the decomposition of relative entropy given for example
in Theorem D.13, \cite{DeZe-98}.
Indeed, from this theorem it follows that
\begin{equation} \label{relEntropDec}
H(\pi_n \nu \vert \rho^n ) = H(\pi_{n-1} \nu \vert \rho^{n-1})
+ \int_{\R^{n-1}} H \big(\pi_n \nu ^{(x_1, \dots, x_{n-1})} \vert {(\rho^n)}^{(x_1, \dots, x_{n-1})} \big)
\pi_{n-1} \nu(d(x_1, \dots,x_{n-1}))
\end{equation}
where for a measure $\mu$ on ${\cal B}(\R^n)$ the regular conditional probability distribution of $\mu$ given $\pi_{n-1}$
is denoted by
$\R^{n-1} \ni (x_1, \dots, x_{n-1}) \mapsto \mu^{(x_1, \dots, x_{n-1})} \in {\cal M}_1(\R^n).$
Thus, to establish (\ref{relEntropRep})
it suffices to show
\begin{equation} \label{substRelEntrEq}
H(\pi_n \nu \vert \pi_{n-1} \nu \otimes \rho)
= \int_{\R^{n-1}} H \big( \pi_n \nu^{(x_1, \dots, x_{n-1})} \vert {(\rho^n)}^{(x_1, \dots, x_{n-1})} \big) \pi_{n-1} \nu (d(x_1, \dots, x_{n-1}))
\end{equation}
for all $n \in \N \backslash \{1\}.$ But applying the quoted theorem to the left-hand side of the
previous equation we obtain
$$
H(\pi_n \nu \vert \pi_{n-1} \nu \otimes \rho)
= \int_{\R^{n-1}} H \big( \pi_n \nu^{(x_1, \dots, x_{n-1})} \vert (\pi_{n-1} \nu \otimes \rho)^{(x_1, \dots, x_{n-1})} \big)
\pi_{n-1} \nu (d(x_1, \dots, x_{n-1})),
$$
and since
$
(\pi_{n-1}\nu \otimes \rho)^{(x_1, \dots, x_{n-1})} = \delta_{(x_1, \dots, x_{n-1})} \otimes \rho
= {(\rho^n)}^{(x_1, \dots, x_{n-1})},
$
(\ref{substRelEntrEq}) follows.
\end{proof}

\begin{prop} \label{affineI}
The function ${\cal I}_{\vert {\cal M}_1^s(\Sigma_b^+)}$ is affine and for $\nu \in {\cal M}_1^{s}(\Sigma_b^+)$ we have
\begin{equation} \label{entEq}
H(\pi_n \nu \vert \pi_{n-1} \nu \otimes \eta) \uparrow {\cal I}(\nu)
\end{equation}
as $n \to \infty.$
\end{prop}
\begin{proof}
We know (cf. Lemma 6.5.16, Corollary 6.5.17 and the preceding discussion in \cite{DeZe-98})
that for $\nu \in {\cal M}_1^s(\Sigma_b^+)$
the value of ${\cal I}(\nu)$ is given as the limit
of the nondecreasing sequence $H(\pi_n \nu \vert \pi_{n-1} \nu \otimes \eta)$ of relative entropies.

To show that ${\cal I}$ restricted to ${\cal M}_1^s(\Sigma_b^+)$
is affine, let $\beta \in (0,1)$ and $\mu, \nu \in {\cal M}_1^s(\Sigma_b^+).$
We distinguish cases:

$(i)$
Assume ${\cal I}(\nu), {\cal I}(\mu) < \infty.$ Then (\ref{entEq}) applies and using
Lemma \ref{relEntropRep} we deduce $\pi_n \nu \ll \eta^n$ and $\pi_n \mu \ll \eta^n$ for all $n \in \N.$
The convexity of relative entropy yields
\begin{align}\label{suppEntComp}
\begin{split}
\beta H(\pi_n \nu \vert \eta^n) &+ (1-\beta) H(\pi_n \mu \vert \eta^n)\\
&\geq H(\beta \pi_n \nu + (1-\beta) \pi_n \mu \vert \eta^n) \\
&\geq \int_{[b,0]^n} \Big( \beta \frac{d \pi_n \nu}{d \eta^n} \log \big( \beta \frac{d \pi_n \nu}{d \eta^n} \big)
+ (1-\beta) \frac{d \pi_n \mu}{d \eta^n} \log \big( (1-\beta) \frac{d \pi_n \mu}{d \eta^n} \big) \Big) d \eta^n \\
&= \beta H(\pi_n \nu \vert \eta^n) + \beta \log \beta
+ (1-\beta) H(\pi_n \mu \vert \eta^n) + (1-\beta) \log (1-\beta).
\end{split}
\end{align}
Dividing by $n$ and taking $n \to \infty$ we obtain in combination with Lemma \ref{relEntropRep} and (\ref{entEq}) that
\begin{equation} \label{affineEq}
{\cal I}(\beta \nu + (1-\beta) \mu) = \beta {\cal I}(\nu) + (1-\beta) {\cal I}(\mu).
\end{equation}

$(ii)$
It remains to consider the case where at least one of the
terms ${\cal I}(\mu), {\cal I}(\nu)$ equals infinity. In this case we want to have
${\cal I}(\beta \nu + (1-\beta) \mu) = \infty,$ and in consideration of (\ref{entEq})
the only nontrivial situation can occur if we have
$H(\pi_n(\beta \nu + (1-\beta) \mu) \vert \pi_{n-1} (\beta \nu + (1-\beta) \mu) \otimes \eta) < \infty$
for all $n \in \N.$ Then
$\pi_n(\beta \nu + (1-\beta) \mu) \ll \pi_{n-1} (\beta \nu + (1-\beta) \mu) \otimes \eta$
and iteratively we deduce
$\pi_n(\beta \nu + (1-\beta) \mu) \ll \eta^n$ and thus
$\pi_n \nu \ll \eta^n$ as well as $\pi_n \mu \ll \eta^n$ for all $n \in \N.$
The same reasoning as in (\ref{suppEntComp}) and (\ref{affineEq}) then yields the desired result.

\end{proof}

\begin{corollary} \label{processLevelIZero}
 The only zero of ${\cal I}$ is given by $\eta^{\N_0}.$
\end{corollary}
\begin{proof}
Proposition \ref{affineI} in combination with
Lemma (\ref{relEntropLem}) shows that for $\nu$ such that ${\cal I}(\nu)$ is finite,
${\cal I}(\nu)$ is given as the limit of the non-decreasing
sequence
$
(H(\pi_n \nu \vert \eta^n)/n)_{n \in \N}.
$
Now since the only zero of $H(\cdot \vert \eta^n)$ is given by $\eta^n,$
we have $H(\pi_n \nu \vert \eta^n) = 0$ for all $n \in \N$ if and
only if $\pi_n \nu = \eta^n$ for all $n \in \N.$
This, however, is equivalent to $\nu = \eta^{\N_0}$ by Kolmogorov's consistency theorem,
which finishes the proof.
\end{proof}

The next lemma is standard.
\begin{lem} \label{extremalPointsLemma}
The set of extremal points of
${\cal M}_1^s(\Sigma_b^+)$ is given by ${\cal M}_1^e(\Sigma_b^+)$.
\end{lem}
\begin{proof}
The proof proceeds analogously to Theorem A.9.10 of \cite{El-85} and is omitted here.
\end{proof}

The following result is closely connected to Proposition \ref{clambda} (b) and shows that the critical
value $\beta_{cr}$ also applies to the constant zero-potential. It is crucial for proving the finiteness
of $L_p^{sup}$ on $(-\infty, \beta_{cr})$ (cf. Lemma \ref{LpsupProps}) and as such in the transition from the variational
formula of Corollary \ref{annealedLyapExpVarFormula} to the representation of the annealed Lyapunov exponents given
in Theorem \ref{annealedLyapExp}.

\begin{lem} \label{betaCrBeta}
Assume (\ref{potIID}) and $h \in (0,1).$ Then
$$
\E_1 \exp\{\beta T_0\} = \infty
$$
for all $\beta > \beta_{cr},$
while
$$
\E_1 \exp\{\beta_{cr} T_0\} < \infty.
$$
\end{lem}
\begin{proof}
The first part of the result follows from the definition of $\beta_{cr}.$ To prove the second
part, we start with
showing that
$\E_1 \exp\{\beta T_0\}$ is finite for all $\beta < \beta_{cr}.$
For this purpose choose such $\beta.$
We now assume
\begin{equation} \label{infinityAssump}
\E_1 \exp \{ \beta T_0\} = \infty
\end{equation}
and lead this assumption to a contradiction. Indeed, setting
$\varepsilon := \beta_{cr}-\beta > 0,$ due to Proposition \ref{clambda} (b) there exists a finite constant
$C_{\beta + \varepsilon/2}$ such that
\begin{equation} \label{asBounded}
\E_1 \exp \Big\{ \int_0^{T_0} (\xi(Y_s) + \beta + \varepsilon/2) \, ds \Big\} \leq C_{\beta+\varepsilon/2} \quad \text{a.s.}
\end{equation}
But
\begin{align} \label{thinSetEst}
\begin{split}
 \ind{1}_{\xi(m) \geq -\varepsilon/2 \, \forall m \in \{1, \dots, n\}}
&\E_1 \exp \Big\{ \int_0^{T_0} (\xi(Y_s) + \beta+ \varepsilon/2) \, ds \Big\} \ind{1}_{Y_s \in \{ 1, \dots, n\} \, \forall s \in [0,T_0)}\\
&\geq \ind{1}_{\xi(m) \geq -\varepsilon/2 \, \forall m \in \{1, \dots, n\}}
\E_1 \exp \{\beta T_0\} \ind{1}_{Y_s \in \{ 1, \dots, n\} \, \forall s \in [0,T_0)}.
\end{split}
\end{align}
With (\ref{infinityAssump}), we deduce
\begin{equation} \label{monConv}
\E_1 \exp \{ \beta T_0\} \ind{1}_{Y_s \in \{1, \dots, n\} \, \forall s \in [0,T_0)} \to \infty
\end{equation}
as $n \to \infty;$ furthermore, due to (\ref{potEssSup}) and (\ref{potIID}), 
$
(\{\xi(m) \geq -\varepsilon /2 \, \forall m \in \{1, \dots, n\}\})_m
$
is a decreasing sequence of sets with positive probability each, and therefore, (\ref{monConv}) in
combination with (\ref{thinSetEst}) yields a contradiction the a.s. boundedness given in (\ref{asBounded}).
Hence, (\ref{infinityAssump}) cannot hold true.
To finish the proof we decompose for $\beta < \beta_{cr}:$
\begin{align*}
\E_1 \exp \{\beta T_0\}
&= \frac{1+h}{2} \frac{\kappa}{\kappa - \beta} + \frac{1-h}{2} \frac{\kappa}{\kappa-\beta} \E_2 \exp\{\beta T_0\}\\
&\geq \frac{(1-h)\kappa}{2(\kappa - \beta)} (\E_1 \exp \{ \beta T_0\})^2.
\end{align*}
Consequently,
$$
\E_1 \exp\{ \beta T_0\} \leq \frac{2(\kappa-\beta)}{\kappa(1-h)}
$$
and monotone convergence yields
$
\E_1 \exp\{ \beta_{cr} T_0\} \leq \frac{2(\kappa-\beta_{cr})}{\kappa(1-h)} < \infty,
$
which finishes the proof.
\end{proof}

\begin{lem} \label{bddCont}
For fixed $\beta \in (-\infty,\beta_{cr}),$
\begin{enumerate}
\item
there exist constants $0 < c < C < \infty$ such that
$$
\E_1 \exp \Big\{ \int_0^{T_0} (\zeta(Y_s) + \beta) \, ds \Big\} \in [c,C]
$$
for all $\zeta \in \Sigma_b^+.$

\item
the mapping
$$
\Sigma_b^+ \ni \zeta \mapsto \log \E_1 \exp \Big\{ \int_0^{T_0} (\zeta (Y_s) + \beta) \, ds \Big\}
$$
is continuous.
\end{enumerate}
\end{lem}
\begin{proof}
$(a)$
For any $\zeta \in \Sigma_b^+$ we have
$$
\E_1 \exp \Big\{ \int_0^{T_0} (\zeta(Y_s) + \beta) \, ds \Big\}
\geq \E_1 \exp \{ (b + \beta) T_0 \} =: c > 0.
$$ 

% as well as
% $$
% \E_1 \exp \Big\{ \int_0^{T_0} ( \zeta(Y_s) + \beta) \, ds \Big\}
% \leq \E_1 \exp \{\beta T_0 \} =: C < \infty,
% $$
% where the last inequality follows from Proposition \ref{clambda} (b) for the constant zero-potential $\xi(x) \equiv 0.$

The upper bound follows from Lemma \ref{betaCrBeta}.

$(b)$
This follows using part (a) and dominated convergence.
\end{proof}
\begin{corollary} \label{VaradhanCond}
For fixed $\beta \in (-\infty, \beta_{cr}),$ the mapping
$$
{\cal M}_1(\Sigma_b^+) \ni \nu \mapsto L (\beta, \cdot)
$$
is continuous and bounded.
\end{corollary}
\begin{proof}
$(a)$
This follows directly from the previous lemma.
\end{proof}

For technical reasons we will need the following two lemmas in the proof of the lower annealed bound;
they can be considered refinements of the corresponding results in the quenched case (cf. Lemma \ref{lDiff} (b) and (c)).
\begin{lem} \label{lVarDiff}
For fixed $\nu \in {\cal M}_1(\Sigma_b^+),$ the mapping
$$
(-\infty, \beta_{cr}) \ni \beta \mapsto L(\beta, \nu)
$$
is continuously differentiable with derivative
\begin{equation}
\frac{\partial L}{\partial \beta}(\beta, \nu)
= \int_{\Sigma_b^+} \frac{\E_1 T_0 \exp \{ \int_0^{T_0} (\zeta(Y_s) + \beta) \, ds \}}
{\E_1 \exp \{ \int_0^{T_0} (\zeta(Y_s) + \beta) \,ds \}} \nu(d \zeta).
\end{equation}
\end{lem}
\begin{proof}
The proof proceeds in analogy to the proof of Lemma \ref{lDiff} (b) and takes advantage of
Lemma \ref{betaCrBeta}.
\end{proof}

\begin{lem} \label{yFindBeta}
\begin{enumerate}
 \item 
For arbitrary $y \in (0, \infty),$ $\nu \in {\cal M}_1(\Sigma_b^+)$
and large enough $M \in (0, \infty),$ there exists $\beta_M(y) \in \R$ such that
$$
y = \int_{\Sigma_b^+} \frac{ \E_{M,1}^\zeta T_0 \exp \{ \beta_M(y) T_0 \}}
{\E_{M,1}^\zeta \exp \{ \beta_M(y) T_0 \}} \nu(d\zeta).
$$

\item
For all $b \in (-\infty, 0),$
$$
\lim_{\beta \to -\infty} \sup_{\nu \in {\cal M}^1(\Sigma_b^+)} \frac{\partial L}{\partial \beta}(\beta, \nu) = 0.
$$
\end{enumerate}
\end{lem}

\begin{proof}
$(a)$

It suffices to show the assertions that
\begin{equation} \label{intInfty}
\int_{\Sigma_b^+} \frac{ \E_{M,1}^\zeta T_0 \exp \{ \beta T_0 \}} { \E_{M,1}^\zeta \exp \{ \beta T_0 \}} \nu(d\zeta) \to \infty
\end{equation}
for $\beta$ large enough
and $M \to \infty$
as well as that 
its integrand
tends to $0$ a.s. for $M$ fixed as $\beta \to -\infty.$
The result then follows from the continuity of
this integrand
and the Intermediate Value theorem in combination with dominated convergence.

The second of these assertions follows as in
the proof of part $(c)$ of Lemma \ref{lDiff},
replacing $\P_{1}^\xi$ by $\P_{M,1}^\zeta.$ For the first assertion,
let $\beta$ be large enough such that $\E^\zeta_1 \exp \{ \beta T_0 \} = \infty$ on a set of positive measure.
It follows for such $\beta$ that
the above integrand tends to infinity
as $M \to \infty$ on the corresponding set, from which we infer
(\ref{intInfty})
for $M \to \infty.$

$(b)$
For this purpose, due to Lemma \ref{lVarDiff}, it suffices to show that
$$
\frac{\partial L}{\partial \beta}(\beta, \delta_\zeta) \to 0
$$
uniformly in $\zeta \in \Sigma_b^+$ as $\beta \to -\infty.$ As in the above we obtain the estimate
(\ref{densityIntEst}), but now $c>0$ can be chosen not to depend on $\zeta \in \Sigma_b^+$ due
to the uniform boundedness of $\zeta.$ Proceeding as in part (a), the claim follows.
\end{proof}

The next lemma states that the supremum in the definition of $L_p^{sup}$ is actually a
maximum when.
% taken over ${\cal M}_1^s(\Sigma_b^+)$ instead of ${\cal M}_1^e(\Sigma_b^+)$ (cf. Lemma \ref{modChoquet}).
\begin{lem} \label{supIsMax}
 For $\beta \in (-\infty, \beta_{cr}),$ there exists $\nu \in {\cal M}_1^s (\Sigma_b^+)$ such that
$
L_p^{sup} (\beta) = L(\beta, \nu) - {\cal I}(\nu)/p.
$
\end{lem}
\begin{proof}
Fix $\beta \in (-\infty, \beta_{cr}).$
Since ${\cal M}_1^s (\Sigma_b^+)$
can be considered
a compact metric space, we find
a converging sequence $(\nu_n)_{n \in \N_0} \subset {\cal M}_1^s(\Sigma_b^+)$ such that
$
L(\beta, \nu_n) - {\cal I}(\nu_n)/p \to L_p^{sup}(\beta)
$ for $n \to \infty.$
As $L(\beta, \cdot)$ is continuous (Corollary \ref{VaradhanCond}) and ${\cal I}$ is lower semi-continuous,
we deduce $L(\beta, \nu) - {\cal I}(\nu)/p = L_p^{sup}(\beta)$ for
$\nu := \lim_{n \to \infty} \nu_n \in {\cal M}_1^s (\Sigma_b^+).$
\end{proof}

In order to deduce the representation for $\lambda_p$ given in Theorem \ref{annealedLyapExp}, we will need the following
lemma.

\begin{lem} \label{LpsupProps}
 The function $\beta \mapsto L_p^{sup} (\beta)$ is finite, strictly increasing, convex and continuous on $(-\infty, \beta_{cr}).$
\end{lem}
\begin{proof}
Lemma \ref{betaCrBeta} implies that $L_p^{sup}$ is finite on $(-\infty, \beta_{cr}).$
 With respect to the strict monotonicty, choose $\nu \in {\cal M}_1^s(\Sigma_b^+)$ such that
$L_p^{sup} (\beta) = L(\beta, \nu) - {\cal I}(\nu)/p,$ which is possible due to Lemma \ref{supIsMax}. The fact
that $L(\cdot, \nu)$ is strictly increasing and $L_p^{sup} \geq L(\cdot, \nu) - {\cal I}(\nu)/p$ now imply that
$L_p^{sup}$ is strictly increasing.

The convexity follows since $L(\cdot, \nu)$ is convex and thus $L_p^{sup}$ as a supremum of
convex functions is convex;
continuity is implied by convexity.
\end{proof}

\section{Proofs for the annealed regime} \label{proofsAnnealed}

The aim of this section is to prove Theorem \ref{annealedLyapExp}.
Similarly to the quenched case we derive upper and lower bounds for
$
t^{-1} \log \langle u(t,0)^p \rangle^{1/p}
$
as $t \to \infty$ (cf. Lemmas \ref{annealedUpperBd} and \ref{annealedLowerBd}).
The additional techniques needed here are Varadhan's Lemma
(see proof of Lemma \ref{annealedUpperBd}) as well as an exponential change of measure (in the
proof of Lemma \ref{annealedLowerBd}), both applied to the sequence $(R_n \circ \xi)_{n \in \N}$ of empirical
measures. Further estimates similar to the quenched regime (Lemmas \ref{annealedLeftSmallerRightExp}
and \ref{annealedFarAwaySummands}) lead to a variational formula for $\lambda_p$ given in
Corollary \ref{annealedLyapExpVarFormula}. Results on the properties of $L_p^{sup}$ (Lemmas \ref{supIsMax} and
\ref{LpsupProps}) then complete of the proof of Theorem \ref{annealedLyapExp}.

As in Theorem \ref{annealedLyapExp}, we assume (\ref{potEssSup}) and (\ref{potIIDBdd})
for the rest of this section.
Notice that since the potential is bounded, 
$L$ is well-defined on $\R.$

\begin{lem} \label{annealedUpperBd}
Let $I \subset [0, \infty)$ be a compact interval and $p \in (0, \infty).$ Then
\begin{equation} \label{annealedUpperBdFormula}
 \limsup_{t \to \infty} \frac1t \log \Big \langle \Big( \sum_{n \in tI} \E_n \exp \Big\{ \int_0^t \xi(Y_s) \, ds \Big\} \ind{1}_{Y_t = 0}
\Big)^p \Big \rangle^\frac1p
\leq \sup_{\alpha \in I} \inf_{\beta < \beta_{cr}}
(-\beta + \alpha L_p^{sup}(\beta)).
\end{equation}
\end{lem}

\begin{proof}

$(i)$
With the same notations as in the quenched case we
first assume $I = [\varepsilon, \gamma]$ with $\varepsilon > 0$ and deduce
% $$
% \limsup_{t \to \infty} \frac1t \log \langle u(t,0)^p \rangle
% \leq -p \beta_{cr} \vee  \sup_{\alpha \in [\gamma_1, \gamma_2]} \limsup_{t \to \infty} \frac1t \log \Big \langle
% \Big ( \E_\floor{\alpha t} \exp \Big \{ \int_0^t \xi(Y_s) \, ds \Big \} \ind{1}_{Y_t = 0} \Big)^p \Big \rangle.
% $$
% Indeed, consider $\alpha \in [\gamma_1, \gamma_2]$ and $\delta>0$ small.
% From (\ref{sumEstimate}) we obtain for $\floor{(\alpha - \delta)t} \leq n \leq \floor{(\alpha + \delta)t}$
using the exponential Chebyshev inequality:
\begin{align}
\Big \langle & \Big ( \sum_{n \in t[\alpha_k^\delta, \alpha_{k+1}^\delta]}
\E_n \exp \Big \{ \int_0^t \xi(Y_s) \, ds \Big \} \ind{1}_{Y_t = 0} \Big)^p \Big \rangle\\
&\leq \Big \langle \Big( \sum_{n \in t[\alpha_k^\delta, \alpha_{k+1}^\delta]}
\E_n \exp \Big\{ \int_0^{T_0} \xi(Y_s) \, ds \Big\} 
\ind{1}_{\frac{T_0}{n} \leq \frac{1}{\alpha_k^\delta}}
\exp \{ c^*(t-T_0)\} \Big)^p \Big \rangle \nonumber \\
&\leq \inf_{\beta > 0} \Big \langle \Big( \sum_{n \in t[\alpha_k^\delta, \alpha_{k+1}^\delta]}
 \E_n \exp \Big\{ \int_0^{T_0} (\xi(Y_s)-\beta) \, ds \Big\}
\exp \{ -\beta_{cr} (t-T_0)\} \Big)^p \Big \rangle
\exp \Big \{\frac{\beta p n }{\alpha_k^\delta} \Big\} \nonumber \\
&\leq \inf_{\beta > 0}
\exp \Big \{ pt \Big( \frac{\beta \alpha_{k+1}^\delta}{\alpha_k^\delta} - \beta_{cr} \Big) \Big\}
\Big \langle \Big( \sum_{n \in t[\alpha_k^\delta, \alpha_{k+1}^\delta]}
\exp \big\{n p L(-\beta + \beta_{cr}, R_n \circ \xi) \big \} \Big)^p \Big \rangle, \label{upperChebyParamEst}
\end{align}
where
to obtain the penultimate line we used (\ref{cStarBetaEq}).
Recall that at the beginning of subsection \ref{annealedRegimeSubSec},
$R_n$ was defined as the empirical measure of a shifted sequence.
Consequently, we conclude
\begin{align}
&\limsup_{t \to \infty} \frac1t \log 
\Big \langle \Big ( \sum_{n \in t[\alpha_k^\delta, \alpha_{k+1}^\delta]}
\E_n \exp \Big \{ \int_0^t \xi(Y_s) \, ds \Big \} \ind{1}_{Y_t = 0} \Big)^p \Big \rangle \nonumber \\
&\leq \inf_{\beta > 0}
\Big[ p \Big(\frac{\beta \alpha_{k+1}^\delta}{\alpha_k^\delta} - \beta_{cr} \Big) 
+ \limsup_{t \to \infty} \frac1t \log
\big \langle \exp \big\{ \alpha_{j_k}^\delta pt L(-\beta + \beta_{cr}, R_\floor{\alpha_{j_k}^\delta t} \circ \xi) \big\} \big \rangle \Big]
\label{annealedFirstUpperSplit}
\end{align}
with $j_k = k+1$ if the second summand in (\ref{annealedFirstUpperSplit}) is positive 
in that case and $j_k = k$ otherwise
(note that this decision depends on $\beta$ but not on the choice of $j_k$).
Corollary \ref{VaradhanCond} tells us that the conditions concerning $L(\beta, \cdot)$
with respect to the upper bound of Varadhan's lemma
(Lemma 4.3.6 in \cite{DeZe-98}) are fulfilled. Thus, bearing in mind the large deviations principle for $R_n \circ \xi$
given in Corollary 6.5.15 of \cite{DeZe-98}
with rate function ${\cal I}$ (cf. (\ref{processLevelRateFunction})),
we can estimate the right-hand side of (\ref{annealedFirstUpperSplit}) by
% \begin{align*}
% \limsup_{t \to \infty} \frac1t &\log
% \big \langle \exp \big\{ \alpha_{j_k}^\delta pt L(-\beta + \beta_{cr}, R_\floor{\alpha_{j_k}^\delta t}) \big\} \big \rangle\\
% &\leq \alpha_{j_k}^\delta \sup_{\nu \in {\cal M}_1(\Sigma_b^+)} \big( pL(-\beta + \beta_{cr}, \nu) -{\cal I}(\nu) \big).
% \end{align*}
% Since ${\cal I}$ is finite for shift-invariant $\nu$ only, we deduce
\begin{align*}
% \limsup_{t \to \infty} &\frac1t \log \Big \langle  \Big ( \sum_{n \in t[\alpha_k^\delta, \alpha_{k+1}^\delta]}
% \E_n \exp \Big \{ \int_0^t \xi(Y_s) \, ds \Big \} \ind{1}_{Y_t = 0} \Big)^p \Big \rangle \nonumber \\
\inf_{\beta > - \beta_{cr}} \Big[ \Big( p \frac{(\beta + \beta_{cr}) \alpha_{k+1}^\delta}{\alpha_k^\delta} -\beta_{cr} \Big)
+ \alpha_{j_k}^\delta \sup_{\nu \in {\cal M}_1^s(\Sigma_b^+)} \big( p L(-\beta, \nu) - {\cal I}(\nu)\big) \Big].
\end{align*}
% Now ${\cal I}$ is lower semi-continuous, whence in combination with
% Corollary \ref{VaradhanCond} the function
% $
% \nu \mapsto L (\beta, \nu) - {\cal I}(\nu)
% $
% is upper semi-continuous. Furthermore, $L(\beta, \cdot)$ is affine and so is ${\cal I}$ on
% ${\cal M}_1^s(\Sigma_b^+)$ (cf. Proposition \ref{affineI}).
% Employing Lemma \ref{extremalPointsLemma},
% as well as Choquet's theorem (see e.g. \cite{Ph-01}, p. 14) and
% Lemma 10.7 of the same reference, we conclude that:
% Employing Lemma \ref{modChoquet} we may restrict the supremum to ergodic measures to obtain in combination
% with the previous:
% \begin{align*}
%  &\limsup_{t \to \infty} \frac1t \log \Big \langle \Big( \sum_{n \in tI} \E_n \exp \Big\{ \int_0^t \xi(Y_s) \, ds \Big\}
% \ind{1}_{Y_t = 0} \Big)^p \Big \rangle\\
% &\leq \max_{k \in \{1, \dots, n-1\}} \inf_{\beta < \beta_{cr}}
% \Big[ p \Big( \frac{(-\beta + \beta_{cr}) \alpha^\delta_{k+1}}{\alpha_k^\delta} -\beta_{cr} \Big)
% + \alpha_{j_k}^\delta \sup_{\nu \in {\cal M}_1^e(\Sigma_b^+)} \big( p L(\beta, \nu) - {\cal I}(\nu) \big) \Big].
% \end{align*}
Since $\inf I > 0$ by assumption, the ratios $\alpha_{k+1}^\delta / \alpha_k^\delta$ are bounded from above and
similarly to the quenched case (proof of Lemma \ref{quenchedUpperBd}) we obtain in combination
with the previous and taking $\delta \downarrow 0:$
\begin{align*}
\limsup_{t \to \infty} \frac1t \log &\Big \langle \Big( \sum_{n \in tI}
\E_n \exp \Big \{ \int_0^t \xi(Y_s) \, ds \Big \} \ind{1}_{Y_t = 0} \Big)^p \Big \rangle^\frac1p
\leq \sup_{\alpha \in I} \inf_{\beta < \beta_{cr}} ( -\beta + \alpha L_p^{sup}(\beta) ).
% = \sup_{\alpha \in I} \inf_{\beta < \beta_{cr}} ( -\beta + \alpha L_p^{sup}(\beta) ).
\end{align*}
% Tatsächlich genügt es hier, zu zeigen, dass zu $\varepsilon > 0$ ein $\beta^*$ klein genug existiert, sodass
% für alle $beta < beta^*$ und $\nu \in {\cal M}(\Sigma_b^+)$ schon $L'(beta, nu) < \varepsilon$ gilt, denn
% dann können wir auch von $\varepsilon-$Supremum sagen, dass beim Betrachten des $\inf_{\beta < \beta_{cr}^0$
% die $\beta < \beta^*$ keine Rolle spielen.
% Here, the equality follows since $\inf I > 0$ and $L_p^{sup} (\beta) = \infty$ for
% $\beta \in (\beta_{cr}^0, \beta_{cr}).$

$(ii)$
Now assume $\inf I = 0.$
We proceed similarly to the quenched part (cf. proof of Lemma \ref{quenchedUpperBd}) and use Proposition
\ref{clambda} (b) to estimate the first factor on the right-hand side of (\ref{quenchedDegenUpperBd}) uniformly by a constant.
% In analogy to the quenched case we then arrive at
% \begin{align*}
% \limsup_{t \to \infty} \frac1t  \log \Big \langle \Big( \sum_{n \in tI} 
% \E_n \exp \Big \{ \int_0^t \xi(Y_s) \, ds \Big \} \ind{1}_{Y_t = 0} \Big )^p \Big \rangle \leq -p\beta_{cr}.
% \end{align*}
% Since the right-hand side of (\ref{annealedUpperBdFormula}) just evaluates to $-\beta_{cr}$ when replacing
% $\sup_{\alpha \in I}$ by $\alpha = 0,$ this finishes the proof.
\end{proof}

For proving the lower bound, we need the following technical lemma.
\begin{lem} \label{technicalExpDecay}
For $\beta \in \R,$ $M > 0,$ $\zeta \in \Sigma_b^+$ and $n \in \N$
denote
$$
\P_{n, \beta, M}^{\zeta} (A) :=
\frac{\E_n \exp \big\{ \int_0^{T_0} (\zeta(Y_s) + \beta) \, ds \big\} \ind{1}_{\tau_k \leq M \, \forall k \in \{1, \dots, n\}}
\ind{1}_{A}}
{\E_n \exp \big\{ \int_0^{T_0} (\zeta(Y_s) + \beta) \, ds \big\} \ind{1}_{\tau_k \leq M \, \forall k \in \{1, \dots, n\}}}
$$
with $A \in {\cal F}.$
Then, for $\varepsilon > 0$ and $\nu \in {\cal M}_1^s(\Sigma_b^+)$ there exist
$\delta > 0$ and a neighbourhood $U(\nu)$ of $\nu$ such that for all $\zeta \in \Sigma_b^+,$
\begin{equation*}
\P_{n, \beta, M}^{\zeta} (T_0 / n \notin (y-\varepsilon, y + \varepsilon) ) \ind{1}_{R_n(\zeta) \in U(\nu)}
\leq 2 \exp\{-n \delta\} \ind{1}_{R_n(\zeta) \in U(\nu)}
\end{equation*}
holds for all $n \in \N,$ where
$$
y := \int_{\Sigma_b^+} \E_{n, \beta, M}^\zeta (T_0) \, d\nu.
$$
\end{lem}
\begin{proof}
Define the function
$$
G_M(\beta): 
\Sigma_b^+ \ni \zeta \mapsto 
\log \E_1 \exp \Big\{ \int_0^{T_0} (\zeta(Y_s) + \beta) \, ds \Big\} \ind{1}_{T_0 \leq M}.
$$
Using the exponential Chebyshev inequality for $\alpha \geq 0$, we compute with $\delta > 0$
\begin{align}\label{QEst}
\begin{split}
&\P_{n, \beta, M}^{\zeta} (T_0 / n \geq y + \varepsilon )
\ind{1}_{R_n(\zeta) \in U(\nu)}
\leq \E^{\zeta}_{n, \beta, M}
\exp \{\alpha T_0 - n\alpha(y+\varepsilon) \}
\ind{1}_{R_n(\zeta) \in U(\nu)}\\
&\leq \exp\{-n \alpha (y+\varepsilon)\} \exp \Big\{
n \Big( \int_{\Sigma_b^+} (G_M(\beta + \alpha) - G_M (\beta)) \, d\nu + \delta \Big)
\Big\}
 \ind{1}_{R_n(\zeta) \in U(\nu)}
\end{split}
\end{align}
for some neighbourhood $U(\nu)$ of $\nu$ (depending on $\alpha$ also) such that
$$
\Big \vert \int_{\Sigma_b^+} (G_M(\beta + \alpha) - G_M(\beta)) \, d\nu
- \int_{\Sigma_b^+} (G_M(\beta + \alpha) -G_M(\beta)) \, d\mu \Big \vert 
\leq \delta
$$
holds for all $\mu \in U(\nu).$
Writing
$$
g(\alpha) := -\alpha(y+ \varepsilon) + \int_{\Sigma_b^+} (G_M(\beta + \alpha) - G_M(\beta)) \, d\nu,
$$
the right-hand side of (\ref{QEst}) equals $\exp\{n(g(\alpha) + \delta)\} \ind{1}_{R_n(\zeta) \in U(\nu)}.$
We observe $g(0) = 0$ and $g'(0) = -y - \varepsilon + y < 0.$
Hence, there exists $\underline{\alpha} > 0$ such that $g(\underline{\alpha}) < 0$.
Setting $\delta := -g(\underline{\alpha})/2$ and refining
$U(\nu)$ such that
\begin{equation*}
\Big \vert \int_{\Sigma_b^+} G_M(\underline{\alpha}) \,d\nu - \int_{\Sigma_b^+} G_M(\underline{\alpha}) \,d\mu \Big \vert < \delta
\end{equation*}
holds for all
$\mu \in U(\nu),$ we deduce from (\ref{QEst}) with $\alpha = \underline{\alpha}$ that
\begin{equation*}
\P_{n, \beta, M}^{\zeta} (T_0 / n \geq y + \varepsilon ) \ind{1}_{R_n(\zeta) \in U(\nu)} \leq \exp\{-n \delta\} \ind{1}_{R_n(\zeta) \in U(\nu)}.
\end{equation*}
% with appropriate adaptations of $\delta$ and $U(\nu)$ for $\alpha \geq 0$:
% \begin{align*}
% Q_{M,n}^{\zeta} (T_0 / n \leq y - \varepsilon )
% \ind{1}_{R_n(\zeta) \in U(\nu)}
% &\leq \exp \Big\{n \Big( \int_{\Sigma_b^+} (G_M(\beta - \alpha) - G_M(\beta)) \, d\nu + \delta \Big)\Big\}\\
% &\quad \times \exp\{n\alpha(y -\varepsilon)\}
% \ind{1}_{R_n(\zeta) \in U(\nu)}
% \end{align*}
% and
% as before
In complete analogy we obtain
\begin{equation*}
\P_{n, \beta, M}^{\zeta} (T_0 / n \leq y - \varepsilon ) \ind{1}_{R_n(\zeta) \in U(\nu)} \leq \exp\{-n \delta\} \ind{1}_{R_n(\zeta) \in U(\nu)},
\end{equation*}
where possibly $\delta > 0$ is even smaller and $U(\nu)$ even more refined; the result then follows.
\end{proof}
We can now proceed to prove the lower annealed bound.
\begin{lem} \label{annealedLowerBd}
Let $I \subset [0, \infty)$ be an interval and $p \in (0, \infty).$ Then
\begin{align} \label{annealedLowerBdFormula}
 \begin{split}
 \liminf_{t \to \infty} \frac1t &\log \Big \langle 
\Big( \sum_{n \in tI} \E_n \exp \Big\{ \int_0^t \xi(Y_s) \, ds \Big\} \ind{1}_{Y_t = 0} \Big)^p \Big \rangle^\frac1p \\
&\geq \sup_{\alpha \in \stackrel{\circ}{I} \cup (I \cap \{0\})} \inf_{\beta < \beta_{cr}} (-\beta + \alpha L_p^{sup}(\beta)).
\end{split}
\end{align}
\end{lem}

\begin{proof}
Observe that for $\alpha \in \stackrel{\circ}{I},$ $y \in (0, 1/\alpha)$ and $\varepsilon > 0$ small enough we have
due to the independence of the medium
\begin{align}
\Big \langle \E_\floor{\alpha t} &\exp \Big\{ \int_0^{T_0} \xi(Y_s) \, ds \Big\} \ind{1}_{Y_t = 0} \Big)^p \Big \rangle^\frac{1}{p}\\
&\geq \Big \langle \Big( \E_\floor{\alpha t} \exp \Big\{ \int_0^{T_0} \xi(Y_s) \, ds \Big\}
\ind{1}_{\frac{T_0}{\floor{\alpha t}} \in ( y-\varepsilon, y + \varepsilon)}
\Big(\E_0 \exp \Big\{ \int_0^{t-r} \xi(Y_s) \,ds \Big\}
\ind{1}_{Y_{t-r} = 0} \Big)_{r=T_0} \Big)^p \Big \rangle \nonumber \\
&\geq \Big \langle \Big( \E_\floor{\alpha t} \exp \Big\{ \int_0^{T_0} \xi(Y_s) \, ds \Big\}
\ind{1}_{\frac{T_0}{\floor{\alpha t}} \in ( y-\varepsilon, y + \varepsilon)} \Big)^p \Big\rangle \nonumber \\
&\quad \times \Big \langle \min_{r \in (y-\varepsilon, y + \varepsilon)}
\Big( \E_0 \exp \Big\{ \int_0^{t(1-\alpha r)} \xi(Y_s) \, ds \Big\} 
\ind{1}_{Y_{t(1-\alpha r)} = 0, T_1 > t(1-\alpha r)} \Big)^p \Big\rangle. \label{divFirstSecPart}
\end{align}
$(i)$ 
To deal with the second of the factors on the right-hand side we observe
\begin{equation} \label{fatouMod}
\frac1t \log \E_0 \exp \Big\{ \int_0^t \xi(Y_s) \, ds \Big\} \ind{1}_{Y_t = 0, \, T_1 > t }
\geq b + \frac1t \log \P_0(Y_s = 0 \, \forall s \in [0,t]) = b-\kappa
\end{equation}
a.s.
% Since
% $
% t \mapsto \log \P_0( Y_t = 0, T_1 > t)
% $
% is super-additive,
% $
% \lim_{t \to \infty} t^{-1} \log \P_0( Y_t = 0, \, T_1 > t)
% $
% exists in $(-\infty, 0]$. Thus, there exist a large real $T$ and some constant $c > -\infty$
% such that a.s. we have
% \begin{equation} \label{fatouMod}
% \frac{1}{t} \log \E_0 \exp \Big\{ \int_0^t \xi(Y_s) \, ds \Big\} \ind{1}_{Y_t = 0,\, T_1 > t} \geq c
% \end{equation}
% for all $t \geq T.$
Hence, Jensen's inequality and a modified version of Fatou's lemma (taking advantage of
(\ref{fatouMod})) apply to yield the inequality in
\begin{align}
\liminf_{t \to \infty} & \frac{1}{t} \log \Big \langle \Big (\E_0 \exp \Big\{ \int_0^{t(1-\alpha y)} \xi(Y_s) \, ds \Big\}
\ind{1}_{Y_{t(1-\alpha y)} = 0,\, T_1 > t(1-\alpha y)} \Big)^p \Big \rangle \nonumber \\
&\geq p \Big \langle \liminf_{t \to \infty} \frac{1}{t} \log \E_0 \exp \Big\{ \int_0^{t(1-\alpha y)} \xi(Y_s) \, ds \Big\}
\ind{1}_{Y_{t(1-\alpha y)} = 0,\, T_1 > t(1-\alpha y)} \Big \rangle \nonumber \\
&= (1-\alpha y)pc^*. \label{divSecPartEst}
\end{align}
To obtain the equality, we used the modification of (\ref{cMEquality}) for $M = 1,$ cf. line thereafter.

$(ii)$
With respect to the first factor on the right-hand side of
(\ref{divFirstSecPart}), we aim to perform an exponential change of measure and
introduce for $\beta \in \R$, $n \in \N$ as well as $M > 0$ the function
$$
G_M^{(n)} (\beta): 
\Sigma_b^+ \ni \zeta \mapsto 
\log \E_n \exp \Big\{ \int_0^{T_0} (\zeta(Y_s) + \beta) \, ds \Big\} \ind{1}_{\tau_k \leq M \, \forall k \in \{1, \dots,n\}}.
$$
Choosing $\nu \in {\cal M}_1^s (\Sigma_b^+)$ we fix $\beta_M \in \R$ such that
\begin{equation} \label{betaMChoice}
\int_{\Sigma_b^+} \frac{\E_1 T_0 \exp \big\{ \int_0^{T_0} (\zeta(Y_s) + \beta_M) \, ds \big\} \ind{1}_{T_0 \leq M}}
{\E_1 \exp \big\{ \int_0^{T_0} (\zeta(Y_s) + \beta_M) \, ds \big\} \ind{1}_{T_0 \leq M}} \, \nu(d\zeta)  = y,
\end{equation}
which is possible due to part (a) of Lemma \ref{yFindBeta}.
Since for $\beta$ fixed,
$
G_M(\beta) := G_M^{(1)}(\beta)
$
is bounded
as a function on $\Sigma_b^+$, we find for each $\delta > 0$
a neighbourhood $U(\nu)$ of $\nu$ in ${\cal M}_1(\Sigma_b^+)$ such that
\begin{equation} \label{neighbourhoodCond}
\Big \vert \int_{\Sigma_b^+} G_M(\beta_M) \,d\nu
- \int_{\Sigma_b^+} G_M(\beta_M) \, d\mu \Big \vert \leq \delta/2
\end{equation}
for all $\mu \in U(\nu)$.
We obtain
\begin{align}
\Big \langle \Big( \E_n & \exp \Big\{ \int_0^{T_0} \xi(Y_s) \, ds \Big\}
\ind{1}_{\frac{T_0}{n} \in (y- \varepsilon, y+\varepsilon)} \Big)^p \Big \rangle
\nonumber \\
&= \Big \langle \Big( \E_n \exp 
\Big\{ \int_0^{T_0} \xi(Y_s) \, ds + \beta_M T_0 - G_M^{(n)}(\beta_M) - \beta_M T_0 +G_M^{(n)}(\beta_M) \Big\}
\times \ind{1}_{\frac{T_0}{n} \in (y- \varepsilon, y+\varepsilon) } \Big)^p \Big \rangle \nonumber \\
&\geq \Big \langle
\Big( \frac{\E_n \exp \big\{ \int_0^{T_0} (\xi(Y_s) + \beta_M) \, ds \big\}
\ind{1}_{\tau_k \leq M \, \forall k \in \{1, \dots, n\}}
\ind{1}_{\frac{T_0}{n} \in (y- \varepsilon, y+\varepsilon)}}
{\E_n \exp \big\{ \int_0^{T_0} (\xi(Y_s) + \beta_M) \, ds \big\} \ind{1}_{\tau_k \leq M \, \forall k \in \{1, \dots, n\}}}\Big)^p \nonumber \\
& \quad \times \exp \big\{ p G_M^{(n)}(\beta_M) \big\}
\ind{1}_{R_n \circ \xi \in U(\nu)} \Big \rangle \exp\{-n \beta_M p (y \pm \varepsilon)\} \nonumber \\
&\geq \exp\{-n \beta_M p (y \pm \varepsilon)\} \exp \Big \{ p n \Big ( \int_{\Sigma_b^+} G_M(\beta_M) \, d\nu - \delta /2 \Big ) \Big \} \nonumber \\
& \quad \times \big \langle \big(1 - \P_{n, \beta, M}^{\zeta} (T_0 /n \notin (y-\varepsilon, y + \varepsilon) ) \big)^p \ind{1}_{R_n \circ \xi \in U(\nu)}
\big \rangle, \label{annealedFirstLower}
\end{align}
where
% $$
% Q_{M,n}^{\xi} (A) :=
% \frac{\E_n \exp \big\{ \int_0^{T_0} (\xi(Y_s) + \beta_M) \, ds \big\} \ind{1}_{\tau_k \leq M \, \forall k \in \{1, \dots, n\}}
% \ind{1}_{A}}
% {\E_n \exp \big\{ \int_0^{T_0} (\xi(Y_s) + \beta_M) \, ds \big\} \ind{1}_{\tau_k \leq M \, \forall k \in \{1, \dots, n\}}}
% $$
% for $A \in \sigma \{ T_0 \}$ and $n \in \N$. Furthermore, 
$y \pm \varepsilon$ is supposed to denote $y+\varepsilon$ if
$\beta_M > 0$ and $y-\varepsilon$ otherwise.

Therefore, choosing $\delta >0$ and $U(\nu)$ according to Lemma \ref{technicalExpDecay}
we infer that
$
\P_{n, \beta, M}^{\zeta} (T_0 /n \notin (y-\varepsilon, y + \varepsilon)) \ind{1}_{R_n \circ \xi \in U(\nu)}
$
decays exponentially in $n.$
Thus,
% small enough such that the derivations leading to (\ref{upperQEst1}) and (\ref{upperQEst2}) hold, 
in combination with the large deviations principle for $R_n \circ \xi$ given in
Corollary 6.5.15 of \cite{DeZe-98} we obtain
% in combination with (\ref{upperQEst1}) and (\ref{upperQEst2})
$$
\liminf_{n \to \infty} \frac{1}{n} \log \big \langle \big( \P_{n, \beta, M}^{\zeta} (T_0 / n \in (y -\varepsilon,y+\varepsilon)) \big)^p
\ind{1}_{R_n \circ \xi \in U(\nu)} \big \rangle \geq -{\cal I}(\nu).
$$
Continuing (\ref{annealedFirstLower}) we get taking $\varepsilon \downarrow 0$ on the right-hand side
\begin{align}\label{firstLowerPart} 
\begin{split}
\liminf_{n \to \infty} &\frac{1}{n} \log \Big \langle \Big( \E_n \exp \Big\{ \int_0^{T_0} \xi(Y_s) \, ds \Big\}
\ind{1}_{\frac{T_0}{n} \in (y- \varepsilon, y+\varepsilon)} \Big)^p \Big \rangle\\
&\geq -\beta_M p y
+ p \Big( \int_{\Sigma_b^+} G_M(\beta_M) d\nu - \delta/2 \Big ) - {\cal I} (\nu).
\end{split}
\end{align}
% One observes that for $M$ large enough, $\beta_M$ as a function in $M$ is non-increasing and from part (b) of Lemma \ref{lambdaChoice}
% we infer that $\beta_M$ as a function in $M$ is bounded from below as $M \to \infty.$
% Hence,
% $(\beta_n)_{n}$ has a converging subsequence $(\beta_{l(n)})_{n \in \N}$ with
% $
% \lim_{n \to \infty} \beta_{l(n)} =: \beta \in (-\infty, \kappa - b].
% $
% Note here, that $U(\nu)$ implicitly depends on $M$ as well. Hence, we proceed as follows: We choose $\beta_\delta < \beta$ such that
% $
% \langle G_n (\beta_\delta) \rangle_\nu + \delta > \langle G_n (\beta) \rangle_\nu;
% $
% further we once again restrict $U(\nu)$, now in such a way that in addition we have
% $
% \langle G_n(\beta_\delta) \rangle_\mu + \delta > \langle G_n (\beta_\delta) \rangle_\nu
% $
% for all $\mu \in U(\nu).$
% The dominated convergence therefore yields that (\ref{firstLowerPart}) is estimated from below by
We observe that
$
M \mapsto G_M(\beta)
$
is nondecreasing, whence the sets
$$
\Big \{ \beta \in \R: -\beta y + \int_{\Sigma_b^+} G_M(\beta) \, d\nu
\leq \lim_{M \to \infty} \inf_{\beta \in \R} \Big( -\beta y + \int_{\Sigma_b^+} G_M(\beta) \, d\nu \Big) \Big \}
$$
are nonincreasing in $M.$ Furthermore, they are non-empty since 
differentiation of 
\begin{equation} \label{infFunc}
\beta \mapsto -\beta y + \int_{\Sigma_b^+} G_M(\beta) \, d\nu
\end{equation}
with respect to $\beta$
yields that this function takes its infimum in $\beta = \beta_M$ (cf. (\ref{betaMChoice})).
From this in combination with the strict convexity of the function in (\ref{infFunc})
we infer the boundedness of the above sets.
Furthermore, these sets are closed since the map in (\ref{infFunc})
is continuous. We therefore conclude that the intersection over all $M > 0$ of these sets contains some
$\beta_\nu \in (-\infty, \beta_{cr}-b]$ and
in combination with (\ref{firstLowerPart})
we deduce
\begin{align}
\liminf_{n \to \infty} &\frac{1}{n} \log \Big \langle \Big( \E_n \exp \Big\{ \int_0^{T_0} \xi(Y_s) \, ds \Big\}
\ind{1}_{\frac{T_0}{n} \in (y- \varepsilon, y+\varepsilon)} \Big)^p \Big \rangle \nonumber \\
&\geq -\beta_\nu p y + p ( L(\beta_\nu, \nu) - \delta/2 ) - {\cal I} (\nu). \label{divFirstFirstPart}
\end{align}
Taking $\delta \downarrow 0,$ (\ref{divFirstSecPart}), (\ref{divSecPartEst}) and (\ref{divFirstFirstPart})
therefore supply us with
\begin{align} \label{annealedLowerInfBetaInserted}
\begin{split} 
\liminf_{t \to \infty} &\frac1t \log 
\Big \langle \Big( \E_\floor{\alpha t} \exp \Big\{ \int_0^t \xi(Y_s) \, ds \Big\} \ind{1}_{Y_t = 0} \Big)^p \Big \rangle^\frac1p\\
&\geq
\alpha \sup_{\nu \in {\cal M}_1^s(\Sigma_b^+)}
\big( -\beta_\nu y + L(\beta_\nu, \nu) - {\cal I}(\nu)/p \big) + (1-\alpha (y-\varepsilon)) c^*\\
&\geq
\alpha \sup_{\nu \in {\cal M}_1^s(\Sigma_b^+)} \inf_{c \leq \beta \leq \beta_{cr}-b}
\big( -\beta y + L(\beta, \nu) - {\cal I}(\nu)/p \big) + (1-\alpha (y-\varepsilon)) c^*,
\end{split}
\end{align}
where the last line holds for $c \in (-\infty, 0)$ small enough due to Lemma \ref{yFindBeta} (b).
Using Sion's minimax theorem (cf. e.g. \cite{Ko-88}),
the first summand of the right-hand side of (\ref{annealedLowerInfBetaInserted})
equals
$$
\alpha \inf_{c \leq \beta \leq \beta_{cr}-b} \sup_{\nu \in {\cal M}_1^s(\Sigma_b^+)}
\big( -\beta y + L(\beta, \nu) - {\cal I}(\nu)/p \big),
$$
where we note that,
mainly due to (\ref{potEssSup}) and Lemma \ref{betaCrBeta},
the infimum over $\beta$ can be restricted to $[c, \beta_{cr}].$
Taking $\varepsilon \downarrow 0$ in (\ref{annealedLowerInfBetaInserted})
and the suprema in $y$ and $\alpha$ we therefore get
\begin{align*}
\liminf_{t \to \infty} &\frac1t \log \Big \langle \sum_{n \in tI} \E_n \exp \Big\{ \int_0^t \xi(Y_s) \,ds \Big\}
\ind{1}_{Y_t = 0} \Big)^p \Big \rangle^\frac1p\\
&\geq \sup_{\alpha \in \stackrel{\circ}{I}} \sup_{y \in (0, 1/\alpha)}
\Big( \alpha \inf_{\beta \leq \beta_{cr}} \sup_{\nu \in {\cal M}_1^s(\Sigma_b^+)}
\big( -\beta y + L(\beta, \nu) - {\cal I}(\nu)/p \big) + (1-\alpha y) c^* \Big)\\
&\geq \sup_{\alpha \in \stackrel{\circ}{I}} \inf_{\beta < \beta_{cr}} ( -\beta + \alpha L_p^{sup}(\beta) )
= \sup_{\alpha \in \stackrel{\circ}{I}} \inf_{\beta < \beta_{cr}} ( -\beta + \alpha L_p^{sup}(\beta) ).
\end{align*}
To obtain the last inequality we used Lemma \ref{supIsMax}, while
the last equality follows
since for $\beta > \beta_{cr}$ we have $L_p^{sup} (\beta) = \infty.$

For the case $0 \in I$ it remains to estimate
$$
\liminf_{t \to \infty} \frac1t \log \Big \langle \Big( \E_0 \exp \Big\{ \int_0^t \xi(Y_s) \, ds \Big\} \ind{1}_{Y_t = 0} \Big)^p \Big \rangle^\frac1p.
$$
Analogously as for the second factor of (\ref{divFirstSecPart}) we obtain
$$
\liminf_{t \to \infty} \frac1t \log \Big \langle 
\Big( \E_0 \exp \Big\{ \int_0^t \xi(Y_s) \, ds \Big\} \ind{1}_{Y_t = 0} \Big)^p \Big \rangle^\frac1p
\geq c^* = - \beta_{cr}.
$$
This finishes the proof.

\end{proof}

Similarly to the quenched case, we have the following two results.
\begin{lem} \label{annealedLeftSmallerRightExp}
For all $\delta \geq 0$ and $\gamma \geq \delta$
we have
 \begin{align}
\limsup_{t \to \infty} &\frac1t \log \Big\langle \Big( \sum_{n \in t[-\gamma , -\delta ]} \E_n \exp \Big\{ \int_0^t \xi(Y_s) ds\Big\}
\ind{1}_{Y_t = 0} \Big)^p \Big \rangle^\frac1p \nonumber \\
&\leq \delta \log \frac{1-h}{1+h}
+ \liminf_{t \to \infty} \frac1t \log \Big\langle \Big( \sum_{n \in t[\delta, \gamma]} \E_n \exp \Big\{ \int_0^t \xi(Y_s) ds \Big\}
\ind{1}_{Y_t = 0} \Big)^p \Big \rangle. \label{annealedOneSidedConc}
\end{align}
\end{lem}
\begin{proof}
The proof proceeds similarly to the proof of Lemma \ref{quenchedLeftSmallerRightExp} and is omitted here.
\end{proof}

\begin{lem} \label{annealedFarAwaySummands}
We have
\begin{equation}
\limsup_{t \to \infty} \frac1t \log
\Big \langle \Big(\sum_{n \notin t[-\gamma, \gamma]} \E_n \exp \Big \{ \int_0^t \xi(Y_s) \, ds \Big\} \ind{1}_{Y_t = 0} \Big)^p\Big \rangle
\to -\infty
\end{equation}
as $\gamma \to \infty$.
\end{lem}
\begin{proof}
The proof is similar to that of Lemma \ref{quenchedFarAwaySummands} and is omitted here.
\end{proof}

We are now ready to prove the existence of $\lambda_p$ and give a variational formula.

\begin{corollary} \label{annealedLyapExpVarFormula}
For $p \in (0, \infty),$ the annealed Lyapunov exponent $\lambda_p$ exists and is given by
 \begin{equation} \label{annealedVarFormula}
\lambda_p = \sup_{\alpha \in [0, \gamma]} \inf_{\beta < \beta_{cr}}
(-\beta + \alpha L_p^{sup}(\beta))
\end{equation}
for all $\gamma > 0$ large enough. 
\end{corollary}
\begin{proof}
Using Lemmas \ref{annealedUpperBd}, \ref{annealedLowerBd}, \ref{annealedLeftSmallerRightExp} and \ref{annealedFarAwaySummands}
one may proceed similarly to the quenched case.
\end{proof}

\begin{proof} [{\bf Proof of Theorem \ref{annealedLyapExp}}]
In order to derive the representation of Theorem \ref{annealedLyapExp},
we distinguish two cases: First assume that $L_p^{sup}$ does not have a zero in $(0, \beta_{cr}).$
From the fact that $L_p^{sup}$ is continuous and increasing (cf. Lemma \ref{LpsupProps}),
we infer that $L_p^{sup} (\beta) < 0$ for all $\beta \in (-\infty, \beta_{cr}),$ whence
taking $\alpha \downarrow 0$ in Corollary
\ref{annealedLyapExpVarFormula}
yields $\lambda_p = -\beta_{cr}.$

Otherwise, if such a zero exists, the properties of $L_p^{sup}$ derived in Lemma \ref{LpsupProps} first imply the
uniqueness of such a zero and then, in an analogous way to the proof of Theorem \ref{quenchedLyapExp}
and in combination with (\ref{annealedVarFormula}), that $\lambda_p$ equals the zero
of $L_p^{sup}(-\cdot).$
% Tats\"achlich, w\"ahle $\alpha^{-1} \in (\frac{d^- L_p^{sup}}{d\beta}(\beta_{zero}), \frac{d^+ L_p^{sup}}{d\beta}(\beta_{zero})),$
% wobei $\beta_{zero}$ die Nullstelle von $L_p^{sup}$ bezeichnet.
\end{proof}

\section{Further results} \label{furtherResults}
While in sections \ref{proofsQuenched} and \ref{proofsAnnealed} we derived the existence
of the corresponding Lyapunov exponents and gave formulae for them, we now concentrate on
their further properties and the proof of Proposition \ref{pIntermittency} and Theorem \ref{lyapExpCont}.
\subsection{Quenched regime}
 Writing $\lambda_0(\kappa)$ to denote the dependence
of $\lambda_0$ on $\kappa,$ we get the following result.

\begin{prop} \label{quenchedLyapConvexKappa}
\begin{enumerate}
\item
$(0, \infty) \ni \kappa \mapsto \lambda_0(\kappa)$ is convex and nonincreasing.

\item
$
\lim_{\kappa \downarrow 0} \lambda_0(\kappa) = 0.
$

\item
The limits
$
\lim_{\kappa \to \infty} \kappa^{-1} \lambda_0(\kappa)
$
and
$
\lim_{\kappa \downarrow 0} \kappa^{-1} \lambda_0(\kappa)
$
exist and are given by
\begin{equation} \label{zeroAsymp}
\lim_{x \downarrow 0} \lim_{t \to \infty} \frac1t \log \E_0 \exp\Big\{ x \int_0^t \xi(X_s) \, ds \Big\} \in (-1,0]
\end{equation}
and
\begin{equation} \label{inftyAsymp}
\lim_{x \to \infty} \lim_{t \to \infty} \frac1t \log \E_0 \exp\Big\{ x \int_0^t \xi(X_s) \, ds \Big\} \in (-1,0],
\end{equation}
respectively. For a non-degenerate potential, both limits are contained in $(-1,0).$
\end{enumerate}
\end{prop}

\begin{proof}
$(a)$
We first show convexity. Writing $u_\kappa (t,x)$ to emphasise the dependence of the solution to (\ref{PAMArb})
on $\kappa,$ we have for 
a random walk $(X_t)_{t \in \R_+}$ with generator $\Delta_h$:
\begin{align}
\lambda_0(\kappa) &= \lim_{t \to \infty} \frac1t \log u_\kappa(t,x)
= \lim_{t \to \infty} \frac1t \log \E_0 \exp \Big\{ \int_0^t \xi(X_{\kappa s}) \, ds \Big\} \nonumber \\
&= \lim_{t \to \infty} \frac1t \log \E_0 \exp \Big\{ \frac{1}{\kappa} \int_0^{\kappa t} \xi(X_s) \, ds \Big\} \nonumber \\
&= \kappa \lim_{t \to \infty} \frac1t \log \E_0 \exp \Big\{ \frac{1}{\kappa} \int_0^t \xi(X_s) \, ds \Big\}
= \kappa \Psi(1/\kappa), \label{quenchedLyapExpAltRep}
\end{align}
where 
$$
\Psi(x) := \lim_{t \to \infty} \frac1t \log \E_0 \exp\Big\{ x \int_0^t \xi(X_s) \, ds \Big\}
$$
for $x \geq 0$.
Note that the limit defining $\Psi(x)$ exists in $[-\kappa,0]$ for all $x \in \R_+$ due to Theorem \ref{quenchedLyapExp}
and (\ref{betaCrEst}).
H\"older's inequality now tells us that $\Psi$ is convex
% AUSFÜHRUNG
% Indeed, for $\beta \in (0,1)$ and $x,y \geq 0$ we
% compute
% \begin{align*}
% \Psi(\beta x + (1-\beta) y) &= \lim_{t \to \infty} \frac1t \log \E_0 \exp\Big\{ \beta x \int_0^t \xi(X_s) \, ds 
% + (1-\beta) y \int_0^t \xi(X_s) \, ds \Big\} \\
% &\leq \beta \lim_{t \to \infty} \frac1t \log \E_0 \exp \Big\{ x \int_0^t \xi(X_s) \, ds \Big\} \\
% &\quad + (1-\beta) \lim_{t \to \infty} \frac1t \log \E_0 \exp \Big\{ y \int_0^t \xi(X_s) \, ds \Big\}\\
% &= \beta \Psi(x) + (1-\beta) \Psi(y).
% \end{align*}
and choosing $\alpha := \frac{\beta x}{\beta x + (1-\beta)y}$ and
$
\gamma := \frac{(1-\beta) y}{\beta x + (1-\beta)y},
$
we obtain the convexity of $x \Psi(1/x)$ in a similar manner:
\begin{align*}
(\beta x + (1-\beta) y) &\Psi \Big( \frac{1}{\beta x + (1-\beta)y} \Big)
= (\beta x + (1-\beta) y) \Psi \Big( \frac{\alpha}{x} + \frac{\gamma}{y} \Big) \\
&\leq \alpha (\beta x + (1-\beta)y) \Psi(1/x) + \gamma (\beta x + (1-\beta) y) \Psi(1/y) \\
&= \beta x \Psi(1/x) + (1-\beta) y \Psi(1/y).
\end{align*}
In combination with (\ref{quenchedLyapExpAltRep})
the convexity of $\kappa \mapsto \lambda_0(\kappa)$ follows.

To show that $\lambda_0(\kappa)$ is nonincreasing in $\kappa \in (0, \infty)$ assume to the contrary that
there is $0 < \kappa_1 < \kappa_2$ such that $\lambda_0(\kappa_1) < \lambda_0(\kappa_2)$. The convexity
of $\lambda_0$ would then imply $\lim_{\kappa \to \infty} \lambda_0(\kappa) = \infty$ which is impossible
since we clearly have $\lambda_0 (\kappa) \leq 0$ for all $\kappa \in (0, \infty)$.

$(b)$
From Theorem \ref{quenchedLyapExp} we deduce 
\begin{equation} \label{quenchedLyapExpBounds}
\lambda_0(\kappa) \in [-\beta_{cr}, 0],
\end{equation}
and using (\ref{betaCrEst}) the claim follows.

$(c)$
The existence and representation
of both limits follows directly from (\ref{quenchedLyapExpAltRep}) and the existence of $\lim_{x \to \infty} \Psi(x)$
as well as $\lim_{x \downarrow 0} \Psi(x),$ which again is due to the monotonicity of $\Psi.$
The fact that the limits are contained in $[-1,0]$ is a consequence of $\Psi(x) \in [-1, 0]$ for all $x \in (0, \infty),$ which
itself is due to Theorem \ref{quenchedLyapExp} for $\kappa = 1.$ 
The fact that for a non-degenerate potential the left-hand side of
(\ref{zeroAsymp}) is different from $0$ can be deduced from
$\lim_{\kappa \downarrow 0} \lambda_0(\kappa) =0,$ the convexity of $\kappa \mapsto \lambda_0(\kappa)$
and the fact that $\lambda_0(1) < 0.$
The last inequality is implied by Theorem \ref{quenchedLyapExp} due to $L(0) < 0.$

As for (\ref{inftyAsymp}), using $\lambda_0(1) < 0,$ $\Psi(1) = \lambda_0(1)$ and the monotonicity of
$\Psi$ we deduce $\lim_{x \to \infty} \Psi(x) < 0$ which finishes the proof.
\end{proof}

\begin{remark}
The result of (a) can be interpreted as follows: The larger $\kappa,$ the harder it gets for the random walk
$X$ appearing in the Feynman-Kac formula (\ref{feynmanKac})
to remain at islands of high peaks of $\xi.$
\end{remark}

\subsection{Annealed regime}
In this subsection we primarily deal with results concerning the annealed Lyapunov exponents. As a special case,
the corresponding results may apply to the quenched regime also.

\begin{prop} \label{lyapExpGenProp}
\begin{enumerate}
\item
The function $p \mapsto \lambda_p$ is nondecreasing in $p \in [0, \infty)$.

\item
The function $p \mapsto p \lambda_p$ is convex in $p \in (0, \infty)$.

\item
For any $p \in [0, \infty)$, $\kappa \mapsto \lambda_p(\kappa)$ is convex in $\kappa \in (0, \infty).$

\item
If $u$ is $p$-intermittent for some $p \in (0,\infty)$, then it is $q$-intermittent for all $q > p$.

\end{enumerate}
\end{prop}

\begin{proof}
$(a)$
For $p > 0$ we directly obtain $\lambda_0 \leq \lambda_p$
from the corresponding formulae given in Theorems $\ref{quenchedLyapExp}$ and $\ref{annealedLyapExp}$. If $0 < p < q$, then Jensen's 
inequality supplies us with
$
\langle u(t,0)^p \rangle^{\frac{1}{p}} \leq \langle u(t,0)^q \rangle^{\frac{1}{q}}
$
and the statement follows from the definition of $\lambda_p$.

$(b)$
For $\beta \in (0,1)$ and $0 < p < q$ we get
$$
\langle u(t,0)^{\beta p + (1-\beta) q} \rangle \leq \langle u(t,0)^p \rangle^\beta \langle u(t,0)^q \rangle^{1-\beta}
$$
by H\"older's inequality,
which implies the desired convexity on $(0, \infty)$.

$(c)$
For $p = 0$ this follows from Proposition \ref{quenchedLyapConvexKappa} (a); for $p \in (0,\infty)$ the proof proceeds in complete analogy
to the corresponding part of the proof of Proposition \ref{quenchedLyapConvexKappa} (a).

$(d)$
Assume to the contrary that $u$ is $p$-intermittent but not $q$-intermittent for some $q > p$.
Then, by the definition of $p$-intermittency,
we have $\lambda_p < \lambda_{p+\varepsilon}$ for all $\varepsilon >0$ and there exists
$\varepsilon^* > 0$ such that $\lambda_q = \lambda_{q+\varepsilon^*}$.
Fixing $\varepsilon := (q-p)/2 \wedge \varepsilon^*,$ we get
using the convexity statement of part (b) and $\lambda_p < \lambda_q:$
\begin{align*}
q \lambda_q &\leq \frac{\varepsilon}{q + \varepsilon - p} p \lambda_p 
+ \frac{q-p}{q+\varepsilon-p}(q+\varepsilon) \lambda_{q+\varepsilon}
< \frac{\varepsilon}{q+\varepsilon-p} p \lambda_q + \frac{q-p}{q+\varepsilon-p}(q+\varepsilon) \lambda_{q+\varepsilon}\\
&= \frac{\varepsilon p / q}{q+\varepsilon-p} q \lambda_q + \frac{(q-p)(q+\varepsilon)/q}{q+\varepsilon-p} q \lambda_q
= q\lambda_q,
\end{align*}
a contradiction. Hence, $u$ must be $q$-intermittent as well.
\end{proof}

\begin{proof}[{\bf Proof of Proposition \ref{pIntermittency}}]

% 
% NICHT LÖSCHEN NICHT LÖSCHEN NICHT LÖSCHEN NICHT LÖSCHEN NICHT LÖSCHEN NICHT LÖSCHEN NICHT LÖSCHEN
% 
% 
% 
% First assume (\ref{intermittencyCrit}) to hold and choose $\nu \in {\cal M}_1^e (\Sigma_b^+)$ such that
% ${\cal I}(\nu) < \infty$ and $\lim_{\beta \uparrow \beta_{cr}^0} L(\beta, \nu) > 0.$ Then for $p$ large enough
% we have $\lim_{\beta \uparrow \beta_{cr}^0} L_p^{sup} (\beta) > 0$ and Theorem \ref{annealedLyapExp} implies
% $L_p^{sup}(-\lambda_p) = 0.$ Due to Lemma \ref{supIsMax} there exists $\nu_p \in {\cal M}_1^e(\Sigma_b^+)$
% such that $L_p^{sup} (-\lambda_p) = L(-\lambda_p, \nu_p) - {\cal I}(\nu_p)/p.$ For $\varepsilon > 0$ we get
% in particular
% $$
% L_{p+\varepsilon}^{sup} (-\lambda_p) \geq L (-\lambda_p \nu_p) - {\cal I}(\nu_p)/(p+\varepsilon) >
% L_p^{sup}(-\lambda_p) = 0,
% $$
% whence Theorem \ref{annealedLyapExp} yields $\lambda_{p+\varepsilon} > \lambda_p.$ Thus, $u$ is $p$-intermittent.
% 
% Now assume (\ref{intermittencyCrit}) is not fulfilled. Then
% $
% \lim_{\beta \uparrow \beta_{cr}^0} L_p^{sup}(\beta) \leq 0
% $
% for all $p > 0$, whence $\lambda_p = -\beta_{cr}$ for all $p > 0$ by Theorem \ref{annealedLyapExp}.
% % In particular, $L$ has no zero in $(0, \beta_{cr}),$ whence due to Theorem \ref{quenchedLyapExp}
% % we have $\lambda_0 = -\beta_{cr}.$
% Hence, $\lambda_p = -\beta_{cr}$ for all
% $p \in (0, \infty)$
% and $u$ is not $p$-intermittent for
% any $p > 0$.
% 
% 
%  NICHT LÖSCHEN NICHT LÖSCHEN NICHT LÖSCHEN NICHT LÖSCHEN NICHT LÖSCHEN NICHT LÖSCHEN NICHT LÖSCHEN
% 
We first show that $L_p^{sup}$ has a zero in $(0, \beta_{cr})$ for $p > 0$ large enough
and then invoke Lemma \ref{supIsMax} to conclude the proof.

To show the existence of such a zero,
let $\mu \in {\cal M}_1 ([b,0])$ such that $H(\mu \vert \eta) < \infty$ and $\mu ([-\beta_{cr}/3,0]) = 1.$
Then due to (\ref{relEntropRep}) and Proposition \ref{affineI} we have
$$
{\cal I}(\mu^{\N_0}) = \lim_{n \to \infty} H(\mu^{n} \vert \mu^{n-1} \otimes \eta) = H(\mu \vert \eta) < \infty
$$
as well as
\begin{align*}
L(\beta_{cr} / 2, \mu^{\N_0}) 
&= \int_{\Sigma_b^+} \log \E_1 \exp \Big\{ \int_0^{T_0} (\zeta(Y_s) + \beta_{cr}/2) \, ds \Big\} \mu^{\N_0}(d\zeta)\\
&\geq \log \E_1 \exp \{ (-\beta_{cr} / 3 + \beta_{cr} /2) T_0 \} > 0.
\end{align*}
We deduce
$$
L_p^{sup} (\beta_{cr}/2) \geq L(\beta_{cr} /2, \mu^{\N_0}) - {\cal I}(\mu^{\N_0})/p > 0
$$
for $p > 0$ large enough, in which case $L_p^{sup}$ has zero $-\lambda_p \in (0, \beta_{cr}),$ cf. Theorem \ref{annealedLyapExp}.

Lemma \ref{supIsMax} now tells us that we find $\nu_p \in {\cal M}_1^s(\Sigma_b^+)$ with
$
L_p^{sup}(-\lambda_p) = L(-\lambda_p, \nu_p) - {\cal I}(\nu_p) / p.
$
Since $\Prob$ can be assumed to be non-degenerate, one can show that for $p$ large enough we have
$\nu_p \not= \Prob.$
We then have ${\cal I}(\nu_p) \in (0,\infty)$ and for $\varepsilon > 0$ we obtain
\begin{align*}
L_{p+\varepsilon}^{sup} (-\lambda_p) \geq L(-\lambda_p, \nu_p) - {\cal I}(\nu_p)/(p+\varepsilon)
> L(-\lambda_p, \nu_p) - {\cal I}(\nu_p)/p = L_p^{sup}(-\lambda_p) = 0.
\end{align*}
Therefore, $L_{p+\varepsilon}^{sup}$ has a zero in $(0, -\lambda_p),$ whence due to Theorem
\ref{annealedLyapExp}
we have $\lambda_{p+\varepsilon} > \lambda_p$
and $u$ is $p$-intermittent.
\end{proof}

The following claim is employed in the proof of Theorem \ref{lyapExpCont}.
\begin{claim}
For each neighbourhood $U$ of $\Prob= \eta^{\N_0}$ in ${\cal M}_1(\Sigma_b^+),$ there exists $\varepsilon >0$ such that
$\{{\cal I} \leq \varepsilon\} \subseteq U$.
\end{claim}
\begin{proof}
Indeed, if this was not the case, we would find an open neighbourhood $U$ of $\Prob$ such that
$
\{{\cal I} \leq \varepsilon\} \not\subseteq U
$
for all $\varepsilon > 0$.
Now since ${\cal I}$ is a good rate function (cf. Corollary 6.5.15 in \cite{DeZe-98}) $\{{\cal I} \leq \varepsilon\} \cap U^c$ is compact and
non-empty whence there exists $\nu \in {\cal M}_1(\Sigma_b^+)$ with ${\cal I} (\nu) = 0$ and $\nu \not\in U$. 
% But
% $H(\pi_n \nu \vert \pi_{n-1} \nu \otimes \eta) \uparrow {\cal I}(\nu)$ as $n \to \infty$ and the only zero of $H( \cdot \vert \pi_{n-1} \otimes \eta)$
% is $\pi_{n-1} \otimes \eta$ itself; it therefore ensues that
% $
% \pi_n \nu = \pi_{n-1} \nu \otimes \eta
% $
% for all $n \in \N$. From this we infer that $\pi_n \nu = \eta^n$ for all $n \in \N$ and thus $\nu = \eta^\Z$, a contradiction
% to $\nu \notin U$. Hence, the claim must hold.
But due to Corollary \ref{processLevelIZero}, $\eta^{\N_0}$ is the only zero of ${\cal I},$ contradicting $\nu \notin U.$
\end{proof}

\begin{proof}[{\bf Proof of Theorem \ref{lyapExpCont}}]
The continuity on $(0,\infty)$ follows from Proposition \ref{lyapExpGenProp} (b). It therefore remains
to show the continuity in $0.$

For this purpose, we first show that $L_p^{sup} \downarrow L$ pointwise as $p \downarrow 0$ on $(0, \beta_{cr}).$

Fix $\beta \in (0, \beta_{cr}).$ Then
$M:= \sup_{\nu \in {\cal M}_1^s(\Sigma_b^+)} L(\beta, \nu) < \infty$ due to Corollary \ref{VaradhanCond}
and for $\varepsilon > 0$
we may therefore find a neighbourhood $U(\Prob)$ of $\Prob$ such that
$
\vert L(\beta, \nu)- L(\beta) \vert < \varepsilon
$
for all $\nu \in U(\Prob).$
Choosing $\delta > 0$ small enough such that $\{ {\cal I} \leq \delta\} \subset U(\Prob)$
(which is possible due to the above claim), we set
$p_\varepsilon := \delta/(M-L(\beta)).$ Then for $p \in (0, p_\varepsilon),$ we have
$\vert L_p^{sup} (\beta) - L(\beta) \vert \leq \varepsilon.$ This proves the above convergence.

The continuity of $p \mapsto \lambda_p$ in zero now follows from Theorems \ref{quenchedLyapExp}
and \ref{annealedLyapExp} where we may distinguish the cases that $L$ does or does not have a zero
in $(0, \beta_{cr}).$
\end{proof}

\section{The case of maximal drift} \label{maxDrift}

In subsection \ref{modProofMaxDrift} we will give the modifications necessary to adapt the proofs leading
to the results of section \ref{mainResults} to the case $h=1.$

Subsequently, in subsection \ref{anMaxDriftCase} we will provide an alternative approach to establish
the existence of the first annealed Lyapunov exponent using a modified subadditivity argument.
By means of the Laplace transform we will then retrieve
an easy formula for the $p$-th annealed Lyapunov exponent for $p \in \N.$

Note that there have been some initial investigations of the first annealed Lyapunov exponent in the case
$h=1$ using a large deviations approach to establish its existence (cf.
\cite{Sc-05}).

\subsection{Modifications in proofs for maximal drift} \label{modProofMaxDrift}

As one may have noticed, some of the results and proofs given so far depended on $h$
being strictly smaller than $1.$
Already Proposition \ref{clambda} does not hold true anymore in the case of maximal drift. Indeed, with the previous definitions
one computes
\begin{equation} \label{cStarBetaRel}
 \beta_{cr} = \kappa \leq \kappa - \xi(0) = -c^*;
\end{equation}
in particular, $c^*$ is in general a non-degenerate random variable.
On the other hand, in the case $h= 1$ we have the simple representations
$$
L(\beta) = \Big \langle \log \frac{\kappa}{\kappa - \xi(1) - \beta} \Big\rangle
\quad
\text{and}
\quad
\Lambda(\beta) = \Big \langle \log \frac{\kappa - \xi(1)}{\kappa - \xi(1) - \beta} \Big\rangle,
\quad \beta \in (-\infty, \kappa).
$$
Notwithstanding these differences between the cases of $h=1$ and $h\in (0,1),$
our main results are still valid in the case $h=1.$ To verify this, we make use of the identity
\begin{equation} \label{zeroSummDecay}
 \E_0 \exp \Big\{ \int_0^t \xi(Y_s) \, ds \Big\} \ind{1}_{Y_t = 0} = \exp \{(-\kappa + \xi(0))t\}.
\end{equation}
We will now exhibit the modifications necessary to derive the results of section \ref{mainResults}.

The proof of Lemma \ref{quenchedLowerBd} is as follows:
\begin{proof}
 For $\alpha > 0$ and bearing in mind (\ref{cStarBetaRel}) and (\ref{zeroSummDecay}),
the supremum on the right-hand side of (\ref{quenchedLowerBdFormula})
is obtained as in the case $h \in (0,1).$ For $\alpha = 0$ it evaluates to $-\beta_{cr} = -\kappa$ and in this
case, choosing for arbitrary $\varepsilon > 0$ an $n \in \N$ such that $\xi(n) > -\varepsilon$ yields
in combination with the Markov property applied at time $t-1$:
\begin{align} \label{favourableBehaviour}
 \begin{split}
&\liminf_{t \to \infty} \frac1t \log \E_n \exp \Big\{ \int_0^t \xi(Y_s) \, ds \Big\} \ind{1}_{Y_t = 0}\\
&\geq \liminf_{t \to \infty} \frac1t \log \Big( \E_n \exp \{ \xi(n)t \} \ind{1}_{Y_{t-1} = n} 
\big( \min_{k \in \{0, \dots, n\}} \exp\{ \xi(k)\} \P_n (Y_0 = n, Y_1 = 0) \big) \Big)\\
&= - \kappa - \varepsilon.
\end{split}
\end{align}
Since $\varepsilon > 0$ was chosen arbitrarily, this finishes the proof.
\end{proof}

Bearing in mind (\ref{cStarBetaRel}) and (\ref{zeroSummDecay}) again, the proof of Lemma \ref{quenchedUpperBd}
proceeds very similarly to the case $h \in (0,1);$ note that, as it will frequently be the case, the proof facilitates
lightly since for $h = 1$ we do not have to consider the negative summands appearing in (\ref{feynmanKac}).
This is also the reason why Lemma \ref{quenchedLeftSmallerRightExp} is not required for $h = 1.$
The proof of Lemma
\ref{quenchedFarAwaySummands} does not depend on $h$ at all, whence no modifications are required.
With these results at hand, Corollary
\ref{quenchedLyapExpVar} is proven as before and the same applies to Theorem \ref{quenchedLyapExp}
and Corollary \ref{optimalSpeed}.

When turning to section \ref{auxiliaryAnnealed},
we note that Lemma \ref{betaCrBeta} is not needed in the case $h=1.$ Furthermore,
for $h= 1$ we note that $\beta_{cr} = \kappa,$ whence Lemma \ref{bddCont}
can be easily verified
to hold true using
$$
\E_1 \exp \Big\{ \int_0^{T_0} (\zeta(Y_s) + \beta) \, ds \Big\} = \frac{\kappa}{\kappa - \zeta(1) - \beta},
\quad \beta < \kappa, \, \zeta \in \Sigma_b^+.
$$

With respect to section \ref{proofsAnnealed}, we note that to derive Lemma \ref{annealedUpperBd}
we just have to employ the relations (\ref{cStarBetaRel}) and (\ref{zeroSummDecay}) in the proof
to obtain the same result.
When it comes to Lemmas \ref{technicalExpDecay} and
\ref{annealedLowerBd}, we observe that the proof goes along similar lines but
facilitates at different steps. But note that e.g. in (\ref{annealedLowerInfBetaInserted})
the infimum in $\beta$ should be taken over $[c, \beta_{cr}-\delta]$ for some $\delta > 0$ small enough
since $L(\beta_{cr})$ might be infinite (whereas the quoted minimax theorem is applicable to
real-valued functions only).

Corollary \ref{annealedLyapExpVarFormula} and Lemma
\ref{LpsupProps} are proven analogously, whence the same applies to Theorem \ref{annealedLyapExp}.

\subsection{Analysis of the maximal drift case} \label{anMaxDriftCase}
When considering annealed Lyapunov exponents for an i.i.d. medium, the
situation that $h=1$ is much easier to analyse
than the case of $h \in (0,1).$ This is the case since in this setting the independence of the medium
yields a product structure for expressions such as
$$
\Big \langle \E_0 \exp \Big\{ \int_0^{T_n} \xi(X_s) \, ds \Big\} \Big \rangle,
$$
which evaluates to $\langle \kappa / (\kappa - \xi(0))\rangle^n.$

\subsubsection{Additional derivations for the annealed regime}
While in general even showing the mere existence of the Lyapunov exponents requires quite some effort, in the case of maximal
drift and an i.i.d. potential, the existence of $\lambda_1$ can be retrieved by a modified subadditivity argument.
\begin{lem} \label{modSubAddPropLemma}
Let
$
f: \R_+ \to \R
$
be a continuous function fulfilling the following property:
For all $\delta > 0$ there exists $K_\delta > 0$ such that for all $s, t \in \R_+$ we have
\begin{equation} \label{modSubAddProp}
f(s+t) \leq K_\delta + \delta s + f(s) + f(t).
\end{equation}
Then $\lim_{t \to \infty} f(t)/t$ exists in $[-\infty, \infty).$
% and we have
% $$
% \lim_{t \to \infty} \frac{f(t)}{t} = \inf_{t \in (0,\infty)} \frac{f(t)}{t} \in [-\infty, \infty).
% $$
\end{lem}
\begin{proof}
For $t$ and $T$ such that $0 < t < T$ choose $n \in \N$ and $r \in [0, t)$ such that $T = nt + r.$ We infer using
(\ref{modSubAddProp}) that
\begin{align*}
\frac{f(T)}{T} &\leq \frac1T \big( (K_\delta + \delta t + f(t))n + f(r) \big)\\
&\leq \frac 1t (K_\delta + \delta t + f(t)) + \frac{f(r)}{T}.
\end{align*}
It follows that
$$
\limsup_{T \to \infty} \frac{f(T)}{T} \leq \liminf_{t \to \infty} \frac{f(t)}{t} + \delta < \infty
$$
for all $\delta > 0$ and thus $\lim_{t \to \infty} f(t)/t$ exists
$
[-\infty, \infty).
$
\end{proof}
When trying to apply this lemma to the function $t \to \log \langle u(t,0)\rangle$ we compute writing
$H(t) := \log \langle e^{t \xi(0)} \rangle,$ denoting by $T_n$ the first hitting time of $n$ by $X,$
and employing the strong Markov property:
\begin{align}
&\langle u(s+t,0) \rangle
= \sum_{n \in \N_0} \Big \langle \E_0 \exp \Big\{ \int_0^{s+t} \xi(X_r) \, dr \Big \}
\ind{1}_{ T_n \leq s < T_{n+1}} \Big \rangle \nonumber\\
&\leq \sum_{n \in \N_0} \E_0 \Big( \Big \langle \exp \Big\{ \int_0^{T_n} \xi(X_r) \, dr \Big\} \Big \rangle
\Big \langle \ind{1}_{T_n \leq s < T_{n+1}} \E_0 \Big( \exp \Big \{ \int_s^{s+t} \xi(X_u) \, du \Big\} \Big \vert {\cal F}_s \Big) \Big \rangle \Big)
\nonumber \\
&= \sum_{n \in \N_0} \E_0 \Big( \Big \langle \exp \Big \{ \int_0^s \xi(X_r) \, dr \Big \} \Big \rangle e^{-H(s-T_n)}
\ind{1}_{T_n \leq s < T_{n+1}} \Big \langle \E_n \exp \Big \{ \int_0^t \xi(X_u) \, du \Big\} \Big \rangle \Big) \label{testMittig}\\
&\leq \sum_{n \in \N_0} \E_0 \Big( \Big \langle \exp \Big \{ \int_0^s \xi(X_r) \, dr \Big \} \Big \rangle
\ind{1}_{ T_n \leq s < T_{n+1}} \Big) e^{-H(s)} \langle u(t,0) \rangle \nonumber\\
& \leq e^{K_\delta + \delta s} \langle u(s,0) \rangle \langle u(t,0) \rangle, \nonumber
\end{align}
where to obtain the last line we used
$
0 \geq \frac{H(t)}{t} \to 0
$
as $t \to \infty$, which implies that for all $\delta >0$ there exists $K_\delta > 0$ such that
$
-H(t) \leq K_\delta + \delta t
$
for all $t > 0.$
Taking logarithms on both sides of (\ref{testMittig}), Lemma \ref{modSubAddPropLemma} is applicable and yields the existence of
$\lambda_1.$ 
It is now promising to consider the Laplace transform
\begin{equation} \label{LaplaceTransform}
\R \ni \beta \mapsto \int_0^\infty e^{-\beta t} \langle u(t,0) \rangle \, dt;
\end{equation}
observe that $\lambda_1$ is given as the critical value of $\beta$ for the divergence of this integral.
By direct computation,  the integral in (\ref{LaplaceTransform}) can be shown to equal
$$
\frac{1}{\kappa} \sum_{n \in \N} \Big \langle \frac{\kappa}{\kappa + \beta - \xi(0)} \Big \rangle^n
$$
for $\beta \geq -\kappa,$ see also Lemma 3.2 in \cite{Sc-05}. Thus,
given the existence of $\lambda_1$ and using (\ref{LaplaceTransform}), we observe that $\lambda_1$ is given
as the zero of
\begin{equation} \label{maxDriftZeroRep}
\beta \mapsto \log \Big \langle \frac{\kappa}{\kappa + \beta - \xi(0)} \Big \rangle
\end{equation}
in $(-\kappa, 0)$ if this zero exists; otherwise, we conclude $\lambda_1 \leq - \kappa$ and 
by considering realisations of $X$ in (\ref{feynmanKac}) which stay at sites $n$ with
$\vert \xi(n) \vert$ small for nearly all the time, we may conclude
$\lambda_1 \geq - \kappa,$ cf. (\ref{favourableBehaviour}). Thus, we get $\lambda_1 = -\kappa$ in this situation.
We have therefore proven the following proposition for $p=1$:
\begin{prop} \label{maxDriftRep}
Assume (\ref{potEssSup}) as well as (\ref{potIIDBdd}) to hold. Then for $h = 1$ and $p \in \N$,
the $p$-th annealed Lyapunov exponent $\lambda_p$ is given as the zero of
$$
\beta \mapsto \log \Big \langle \Big( \frac{\kappa}{\kappa + \beta - \xi(0)} \Big)^p \Big \rangle
$$
in $(-\kappa, 0)$ if this zero exists and $-\kappa$ otherwise.
\end{prop}

\begin{remark}
 While Theorem \ref{annealedLyapExp} yields the existence and implicit formulae
for all $\lambda_p$, $p \in (0, \infty)$, simultaneously,
Proposition \ref{maxDriftRep} provides a nicer representation in the cases $p \in \N$ with $h=1.$
\end{remark}

\begin{proof}
While for $p=1$ we
showed how to employ a subadditivity argument to infer existence of $\lambda_1$,
for general $p \in \N$ we now refer to Theorem
\ref{annealedLyapExp} for this purpose. We can then use the Laplace transform again to deduce
a more convenient representation of $\lambda_p$. For the sake of simplicity, we prove the proposition
for $p=2$ and give corresponding remarks where generalisations to arbitrary $p \in \N$ are not straightforward.

Denote by $X^{(1)}$ and $X^{(2)}$ two independent copies of $X$ and by $\P_{0,0}$ and $\E_{0,0}$ we denote
the probability and expectation, respectively, of these processes both starting in $0$. Note that since $h= 1,$ these are Poisson processes
with intensity $\kappa$. We set $\tau_0^{(j)} := 0$, $\tau_{k}^{(j)} := T_k^{(j)} - T_{k-1}^{(j)}$
for $k \in \N$
and $j \in \{1, 2\},$ where by $T_k^{j}$ we denote the first hitting time of $k$ by $X^{(j)}.$
Note that the $\tau_{k}^{(j)}$ are i.i.d. exponentially distributed with parameter $\kappa.$
We distinguish cases:

$(i)$
Assume 
that (\ref{maxDriftZeroRep}) has a zero in $(-\kappa, 0).$

Using H\"older's inequality\footnote{For arbitrary $p$ we retreat to the generalised H\"older inequality with
the $p$ exponents $1/p, \dots, 1/p.$}
we estimate for
\begin{equation} \label{lambdaCond}
\beta > -2\kappa:
\end{equation}
\begin{align} \label{maxDriftUpperEstI}
\begin{split}
&\int_0^\infty e^{-\beta t} \langle u(t,0)^2 \rangle \, dt\\
&= \int_0^\infty e^{-\beta t} \sum_{m,n \in \N_0} \Big \langle \E_{0,0} \Big( \exp \Big\{ \sum_{k=1}^m \tau_k^{(1)} \xi(k-1) + (t-T_m^{(1)}) \xi(m) \Big\}\\
&\quad \times \exp \Big \{ \sum_{k=1}^n \tau_k^{(2)} \xi(k-1) + (t-T_n^{(2)}) \xi(n) \Big\}  \ind{1}_{X_t^{(1)}=m} \ind{1}_{X_t^{(2)} = n} \Big)
\Big \rangle \, dt\\
&\leq \sum_{m,n \in \N_0} \Big \langle \int_0^\infty e^{-\beta t} \Big( \E_{0,0} \exp \Big\{ \sum_{k=1}^m \tau_k^{(1)} \xi(k-1) + (t-T_m^{(1)}) \xi(m) \Big \}
\ind{1}_{X_t^{(1)} = m} \Big)^2 \, dt \Big \rangle^{\frac12}\\
&\quad \times \Big \langle \int_0^\infty e^{-\beta t} \Big( \E_{0,0} \exp \Big \{ \sum_{k=1}^n \tau_k^{(2)} \xi(k-1) + (t-T_n^{(2)}) \xi(n) \Big \}
\ind{1}_{X_t^{(2)} = n} \Big)^2 \, dt \Big \rangle^{\frac12}.
\end{split}
\end{align}
% Since we have
% \begin{align*}
% \Big( &\E_{0,0} \exp \Big \{ \sum_{k=1}^n \tau_k^{(j)} \xi(k-1) + (t-T_n^{(j)}) \xi(n) \Big \}
% \ind{1}_{X_t^{(j)} = n} \Big)^2 \\
% &=\E_{0,0} \exp \Big \{ \sum_{k=1}^n \tau_k^{(1)} \xi(k-1) + (t-T_n^{(1)}) \xi(n) \Big \}
% \ind{1}_{X_t^{(1)} = n} \\
% & \quad \times  \E_{0,0} \exp \Big \{ \sum_{k=1}^n \tau_k^{(2)} \xi(k-1) + (t-T_n^{(2)}) \xi(n) \Big \}
% \ind{1}_{X_t^{(2)} = n}
% \end{align*}
% for $j \in \{1,2\}$,
We now can
estimate the diagonal summands as follows:
\begin{align} \label{maxDriftUpperEstII}
\begin{split}
&\Big \langle \int_0^\infty e^{-\beta t} \E_{0,0} \Big( \exp \Big\{ \sum_{k=1}^m (\tau_k^{(1)} + \tau_k^{(2)}) \xi(k-1)\Big\}
\exp \{ (2t - T_m^{(1)} - T_m^{(2)}) \xi(m) \} \ind{1}_{X_t^{(1)} = X_t^{(2)} = m} \Big) \, dt \Big \rangle\\
&= \Big \langle \int_0^\infty e^{-\beta t} \E_{0,0} \Big( \exp \Big \{ \sum_{k=1}^m (\tau_k^{(1)} + \tau_k^{(2)} ) \xi(k-1)
+ (2t - T_m^{(1)} - T_m^{(2)}) \xi(m) \Big\}\\
&\quad \times \exp \Big \{ - \kappa ( 2t - T_m^{(1)} - T_m^{(2)} ) \Big\}
\ind{1}_{T_m^{(1)} \leq t} \ind{1}_{T_m^{(2)} \leq t } \Big) \, dt \Big \rangle\\
&\leq \E_{0,0} \Big \langle \int_0^\infty \exp \Big \{ \sum_{k=1}^m (\tau_k^{(1)} + \tau_k^{(2)}) (\xi(k-1) - \beta/2) \Big \}\\
&\quad \times \exp \big \{(2t - T_m^{(1)} - T_m^{(2)}) (\xi(m) - \kappa - \beta/2) \big\}
\ind{1}_{T_m^{(1)} + T_m^{(2)} \leq 2t } \, dt \Big \rangle\\
&\hspace{-2.4em}\stackrel{t \mapsto t + \frac{T_m^{(1)} + T_m^{(2)}}{2}}{=}
\E_{0,0} \Big \langle \exp \Big \{ \sum_{k=1}^m (\tau_k^{(1)} + \tau_k^{(2)})(\xi(k-1) - \beta/2) \Big\} \Big \rangle
\underbrace{
\Big \langle \int_0^\infty \exp\{2t(\xi(0) - \kappa - \beta/2)\} \, dt \Big \rangle
}_{=: C< \infty, \text{ since }
\beta > -2\kappa}\\
&= C \Big \langle \Big( \frac{\kappa}{\kappa +\beta/2 - \xi(0)} \Big)^2 \Big \rangle^m 
\end{split}
\end{align}

Hence, combining (\ref{maxDriftUpperEstI}) and (\ref{maxDriftUpperEstII}) we have
\begin{align}
\int_0^\infty e^{-\beta t} \langle u(t,0)^2 \rangle \, dt
&\leq C^2 \sum_{m,n \in \N_0} \Big \langle \Big( \frac{\kappa}{\kappa + \beta/2 - \xi(0)} \Big)^2 \Big 
\rangle^{\frac{m+n}{2}} \nonumber\\
&= C^2 \Big( \sum_{m \in \N_0} \Big \langle \Big( \frac{\kappa}{\kappa + \beta/2 - \xi(0)} \Big)^2 \Big 
\rangle^{\frac{m}{2}} \Big)^2. \label{annQuadUpperBd}
\end{align}

For the lower bound we compute
\begin{align}
\int_0^\infty & e^{-\beta t} \langle u(t,0)^2 \rangle \, dt 
\geq \sum_{m \in \N_0} \int_0^\infty \Big \langle \E_{0,0} \exp \Big \{ \sum_{k=1}^m (\tau_k^{(1)} + \tau_k^{(2)})(\xi(k-1) - \beta/2)\Big\}
\nonumber \\
&\quad \times \exp \big\{ (2t - T_m^{(1)} - T_m^{(2)}) (\xi(m) - \beta/2) \big\} \ind{1}_{X_t^{(1)} = X_t^{(2)} = m} \Big \rangle \, dt \nonumber \\
&= \sum_{m \in \N_0} \Big \langle \E_{0,0} \Big( \int_0^\infty
\exp \Big \{ \sum_{k=1}^m (\tau_k^{(1)} + \tau_k^{(2)})(\xi(k-1) - \beta/2) \Big \} \nonumber \\
&\quad \times \exp \big \{ (2t - T_m^{(1)} - T_m^{(2)})(\xi(m) - \kappa - \beta/2) \big\} \ind{1}_{T_m^{(1)} \leq t}
\ind{1}_{T_m^{(2)} \leq t} \, dt \Big) \Big \rangle \nonumber \\
&\hspace{-2.4em} \stackrel{t \mapsto t + \frac{T_m^{(1)} + T_m^{(2)}}{2}}{=} \sum_{m \in \N_0} \Big \langle \E_{0,0} \Big( \int_0^\infty
\exp \Big \{ \sum_{k=1}^m (\tau_k^{(1)} + \tau_k^{(2)})(\xi(k-1) - \beta /2) \Big\} \nonumber \\
&\quad \times \exp \{ 2t(\xi(m) - \kappa - \beta /2) \} \ind{1}_{\frac{T_m^{(1)} - T_m^{(2)}}{2} \leq t}
\ind{1}_{\frac{T_m^{(2)} - T_m^{(1)}}{2} \leq t} \, dt \Big) \Big \rangle \label{indFuncs}\\
&\hspace{-2.6em} \stackrel{t \mapsto t + \frac{ \vert T_m^{(1)} - T_m^{(2)} \vert }{2}}{=} \sum_{m \in \N_0} 
\Big \langle \E_{0,0} \Big( \int_0^\infty \exp \Big \{ \sum_{k=1}^m (\tau_k^{(1)} + \tau_k^{(2)})(\xi(k-1) - \beta /2) \Big\} \nonumber \\
&\quad \times \exp \big \{ (2t + \vert T_m^{(1)} - T_m^{(2)} \vert )(\xi(m) - \kappa - \beta/2) \big \} \, dt \Big) \Big \rangle \nonumber \\
\begin{split}&\geq \sum_{m \in \N_0} \Big \langle \E_{0,0} \Big( \ind{1}_{ \vert T_m^{(1)} - T_m^{(2)} \vert \leq m\delta } 
\exp \Big \{ \sum_{k=1}^m (\tau_k^{(1)} + \tau_k^{(2)})(\xi(k-1) - \beta/2) \Big\} \Big) \Big \rangle \\
&\quad \times \Big \langle \exp\{ m\delta (\xi(m) - \kappa - \beta/2) \} 
\int_0^\infty \exp \{2t(\xi(m) - \kappa - \beta /2)\} \, dt \Big \rangle.\end{split} \label{tempLastLine}
\end{align}
Note here that for arbitrary $p \in \N$ the indicators appearing in (\ref{indFuncs}) are replaced by
$$
\prod_{j=1}^p \ind{1}_{T_m^{(j)} - \frac{\sum_{1\leq k\leq p} T_m^{(k)}}{p} \leq t}
$$
which can be estimated from below by
$$
\ind{1}_{\max_{1\leq j,k \leq p} \vert T_m^{(j)} - T_m^{(k)}\vert \leq t}.
$$
The subsequent substitution can duely be replaced by
$$
t \mapsto t + \frac{\max_{1\leq j,k \leq p} \vert T_m^{(j)} - T_m^{(k)}\vert }{p},
$$
and the remaining steps are analogous to $p=2.$

Now we continue (\ref{tempLastLine}) with $p=2$ and
bearing in mind that $\beta > -2\kappa$, we estimate the right-hand side factor using Jensen's inequality to get
\begin{align} \label{maxDriftSecondFactorEst}
 \begin{split}
\Big \langle \exp &\{ m\delta (\xi(m) - \kappa - \beta/2) \} \int_0^\infty \exp\{ 2t(\xi(m) - \kappa - \beta/2)\} \, dt \Big \rangle\\
&\geq \frac12 \Big \langle
\underbrace{
\exp\{\xi(m) - \kappa - \beta/2\}^\delta
}_{\to 1 \Prob\text{-a.s. as } \delta \downarrow 0}
\underbrace{
\Big( \frac{1}{\kappa + \beta/2 - \xi(m)} \Big)^\frac1m
}_{\to 1 \Prob\text{-a.s. as } m \to \infty}
\Big \rangle^m\\
&\geq \frac{(1-\varepsilon)^m}{2},
\end{split}
\end{align}
where the last inequality holds for arbitrary $\varepsilon >0$ and all $m \geq m_{\delta, \varepsilon}$ large enough.
Writing $H(t) := \log \langle e^{t \xi(0)} \rangle$ again, we obtain combining (\ref{tempLastLine}) and
(\ref{maxDriftSecondFactorEst}):
\begin{align*}
\int_0^\infty e^{-\beta t} \langle u(t,0)^2 \rangle \, dt
&\geq \frac12 \sum_{m \geq m_\varepsilon} (1-\varepsilon)^m\E_{0,0} \Big( \ind{1}_{\vert T_m^{(1)} - T_m^{(2)} \vert \leq m\delta}\\
&\quad \times \exp \Big \{ \sum_{k=1}^m  \big( H(\tau_k^{(1)} + \tau_k^{(2)}) - \beta/2(\tau_k^{(1)} + \tau_k^{(2)}) \big) \Big\} \Big).
\end{align*}
Next we define
$$
\hat{\P}_m(A):= \frac{\E_{0,0} (\exp\{ \sum_{k=1}^m H(\tau_k^{(1)} + \tau_k^{(2)}) - \beta/2(\tau_k^{(1)} + \tau_k^{(2)})\}\ind{1}_A)}
{
\underbrace{
\E_{0,0} \big( \exp \big \{ \sum_{k=1}^m H(\tau_k^{(1)} + \tau_k^{(2)}) - \beta/2(\tau_k^{(1)} + \tau_k^{(2)}) \big \} \big)
}_{= \langle (\frac{\kappa}{\kappa + \beta/2 - \xi(0)})^2 \rangle^m < \infty \text{ as } \beta > -2\kappa,\text{ cf. } (\ref{lambdaCond})}
}
$$
for $m \in \N$ and measurable $A.$
Now $(\tau_k^{(1)} - \tau_k^{(2)})_{k \in \{1, \dots, m\}}$ have mean $0$ and are square integrable and
i.i.d. with respect to $\hat{\P}_m.$
Thus, a weak law of large numbers supplies us with
\begin{align}
\int_0^\infty e^{-\beta t} \langle u(t,0)^2 \rangle \, dt 
&\geq
\sum_{m \geq m_{\delta, \varepsilon}} \hat{\P}_m (\vert T_m^{(1)} - T_m^{(2)} \vert \leq \delta m )
\Big \langle \Big( \frac{\kappa}{\kappa + \beta/2 - \xi(0)} \Big)^2 \Big \rangle^m\\
&\quad \times \frac{(1-\varepsilon)^m}{2} \notag \\
&\geq \frac14 \sum_{m \geq m_{\delta, \varepsilon}} \Big ( \Big \langle \Big( \frac{\kappa}{\kappa + \beta/2 - \xi(0)} \Big)^2
\Big \rangle (1-\varepsilon) \Big)^m, \label{annQuadLowerBd}
\end{align}
where we choose $m_{\delta, \varepsilon}$ large enough such that
$
\hat{\P}_m (\vert T_m^{(1)} - T_m^{(2)} \vert \leq \delta m) \geq 1/2
$
for all $m \geq m_{\delta, \varepsilon}$ due to the law of large numbers.

Since $\varepsilon > 0$ was chosen arbitrarily, we infer combining
(\ref{annQuadUpperBd}) and (\ref{annQuadLowerBd})
that $\lambda_2$ equals the zero of
$$
\beta \mapsto \log \Big \langle \Big( \frac{\kappa}{\kappa + \beta - \xi(0)} \Big)^2 \Big \rangle.
$$
$(ii)$
Now assume 
that (\ref{maxDriftZeroRep}) does not have a zero in $(-\kappa, 0).$

Again, considering realisations of $X$ in (\ref{feynmanKac}) which stay at sites $n$ with $\vert \xi(n) \vert$ small
for nearly all the time, we arrive at
$
\lambda_2 \geq - \kappa,
$
cf. (\ref{favourableBehaviour}).
But from (\ref{annQuadUpperBd}) we infer
$
\lambda_2 \leq -\kappa$ if (\ref{maxDriftZeroRep}) does not have a zero in $(-\kappa, 0),$
which finishes the proof.
\end{proof}

% However, regarding Theorem \ref{annealedLyapExp} and the dependence of $\lambda_p$ on $p$ therein, we do not expect
% Proposition \ref{maxDriftRep} to hold for all $p \in \N$.
It is inherent to the approach along which we proved Proposition \ref{maxDriftRep} that it applies
to natural $p$ only. Nevertheless, we expect the formula to hold true for general $p \in (0, \infty).$
\begin{conjecture}
Assume (\ref{potEssSup}) as well as (\ref{potIIDBdd}) to hold. Then for $h = 1$ and $p \in (0,\infty)$,
the $p$-th annealed Lyapunov exponent $\lambda_p$ is given as the zero of
$$
\beta \mapsto \log \Big \langle \Big( \frac{\kappa}{\kappa + \beta - \xi(0)} \Big)^p \Big \rangle
$$
in $(-\kappa, 0)$ if this zero exists and $-\kappa$ otherwise.
\end{conjecture}

{\bf Acknowledgment.}
I would like to thank J\"urgen G\"artner and Rongfeng Sun for various fruitful discussions. Moreover, I would like to thank
Wolfgang K\"onig for a helpful suggestion on the computation of the annealed Lyapunov exponents.

\newpage

\bibliographystyle{alpha}
\bibliography{pamBib}

\end{document}